\documentclass[10pt]{amsart}
\theoremstyle{plain}
\usepackage{amsmath, amssymb, amscd,graphicx,jrsymb,color}
\usepackage[all]{xy}
\newtheorem{lem}{Lemma}[section]
\newtheorem{prop}[lem]{Proposition}
\newtheorem{thrm}[lem]{Theorem}
\newtheorem{cor}[lem]{Corollary}

\newtheorem*{defn*}{Definition}

\newtheorem{tm}{Theorem}

\newtheorem{crrr}{Corollary}

\theoremstyle{definition}
\newtheorem{defn}[lem]{Definition}

\newcommand{\R}{{\mathbb{R}}}

\newcommand{\Z}{{\mathbb{Z}}}
\newcommand{\Q}{{\mathbb{Q}}}

\newcommand{\wt}{\widetilde}

\newcommand{\im}{{{\text {im} \ }}}

\newcommand{\gr}{{{\mathrm{gr} }}}

\newcommand{\ts}{{{\thinspace}}}

\renewcommand{\H}{\overline{H}}
\newcommand{\Hpm}{H^\pm}
\newcommand{\Hm}{H}
\newcommand{\Hu}{\widetilde{H}}
\newcommand{\Cu}{\widetilde{C}}
\newcommand{\Cr}{\overline{C}}

\newcommand{\cc}{\mathcal{C}}
\newcommand{\ck}{\mathcal{K}}
\newcommand{\ccp}{\mathcal{C}^+}
\newcommand{\ckp}{\mathcal{K}^+}
\newcommand{\od}{\overline{D}}
\renewcommand{\wt}{\widetilde}

\title[Some differentials on Khovanov-Rozansky homology]{Some differentials on
  Khovanov-Rozansky homology}
\author{Jacob Rasmussen}
\address{Princeton University Dept. of Mathematics, Princeton, NJ 08544}
\email{jrasmus@math.princeton.edu}
\thanks{The author was partially supported by an NSF Postdoctoral fellowship.}

\begin{document}

\begin{abstract}
We study the relationship between the HOMFLY and \(sl(N)\) knot
homologies introduced by Khovanov and Rozansky. For each \(N>0\), we
show there is a spectral sequence which starts at the HOMFLY homology
and converges to the \(sl(N)\) homology.  
As an application, we determine the KR-homology of knots with 9
crossings or fewer.
\end{abstract}

\maketitle

\section{Introduction}

In \cite{KRI,KRII}, Khovanov and Rozansky introduced a new class of
 homological knot invariants which generalize the original
construction of the Khovanov homology \cite{Khovanov}. 
In this paper, we investigate these {\it
  KR-homologies} and the relations between them.  Our
motivation was to give some substance to the conjectures
made in  \cite{superpolynomial} about the behavior of these
theories and their relation to the knot Floer homology. Although we
are unable to say anything about the latter problem, we hope that we can at
least shed some light on the structure of  KR-homology. 

In order to state our results, we briefly recall the form of 
these homologies, restricting for the moment to the case of a knot \(K \subset
S^3\). To such 
\(K \),  
the theory of \cite{KRII} assigns a triply-graded homology group
\(\H^{i,j,k}(K)\) whose graded Euler characteristic is the HOMFLY
polynomial. To be precise, we denote by \(P_K(a,q)\) the
HOMFLY polynomial of \(K\) normalized to satisfy the skein relation
\begin{equation*}
a P(\undercrossing) - a^{-1} P(\overcrossing) = (q-q^{-1})
P(\asmoothing), 
\end{equation*} 
and so that \(P\) of the unknot is
equal to \(1\).
Then with an appropriate choice of gradings, 
\begin{equation*}
\sum_{i,j,k} (-1)^{(k-j)/2}a^jq^i \dim \H^{i,j,k}(K) = P_K(a,q). 
\end{equation*}

The definition of \(\H\) is closely related to that of
another family of homology theories \(\H_N^{I,J}(K)\) (\(N>0\)) 
introduced by
Khovanov and Rozansky in \cite{KRI}. Their graded Euler
characteristics give the \(sl(N)\) polynomials:
\begin{equation*}
\sum_{I,J}(-1)^J q^I \dim \H_N^{I,J} (K) = P_K(q^N,q).
\end{equation*}
The large \(N\) behavior of these theories was studied first by Gukov,
Schwartz, and Vafa in \cite{GSV}, and then later in 
\cite{superpolynomial}, where it was conjectured that 
the limit of \(\H_N(K)\) as \(N\to \infty\) should be a triply graded
homology theory \(\mathcal{H}(K)\) with Euler characteristic \(P_K(a,q)\). 
In fact, 
the limiting theory is a regraded version of \(\H(K)\).

\begin{tm}
\label{Thm:Limit}
For all sufficiently large \(N\),
\( \displaystyle
\H_N^{I,J}(K) \cong \bigoplus_{\substack{i+Nj=I \\ (k-j)/2=J}} \H^{i,j,k}(K)\).
\end{tm}
\noindent We remark
 that \(\H(K)\) is finite dimensional, so when \(N\) is large, 
there will be at most one nontrivial summand on the right-hand
side. The exact value of \(N\) needed for the theorem to hold depends
on \(K\), but it need not be especially big. There are
many knots for which \(N>1\) is enough.


Theorem~\ref{Thm:Limit}
 is a special case of the  following more general relation
between \(\H\) and \(\H_N\).

\begin{tm}
\label{Thm:d+}
 For each \(N>0\), there is a spectral sequence \(E_k(N)\) which
 starts at \(\H(K)\) and converges to \(\H_N(K)\).  
The higher terms in this sequence are  invariants of \(K\).   
\end{tm} 

In some sense, these sequences are all generalizations of Lee's
original spectral sequence \cite{ESL2} for the Khovanov homology. As
described in \cite{superpolynomial}, the idea that they should exist
arose from Gornik's work on the \(sl(N)\) homology \cite{Gornik}. The
exact method by which they are constructed is rather different from
that envisioned in \cite{superpolynomial}, but we expect that their
content is the same.

Strictly speaking, the statement of the theorem is
weaker than the conjecture made in
\cite{superpolynomial}, which says that \(\H_N(K)\) should be the
homology of a differential \(d_N:\mathcal{H}(K) \to \mathcal{H}(K)\). 
The first differential in the sequence \(E_k(N)\) provides 
us with a map  \(d:\H(K) \to
\H(K)\) whose behavior  with respect to the triple grading on
\(\H\) matches that predicted for \(d_N\). Thus if we knew that the
spectral sequence converged after the first differential, this
part of the conjecture would hold. 
In all the examples we have considered, \(E_k(N)\) does indeed
converge after the first differential, but we see no 
 {\it a priori} reason why this should always be the case. 

More generally, 
 \cite{superpolynomial} conjectured that there should be
 differentials \(d_N:\mathcal{H}(K) \to \mathcal{H}(K)\)
  not just for \(N>0\), but for all \(N \in \Z\). Furthermore, 
  \(d_N\) and \(d_{-N}\)  should be
 interchanged by an involution \(\phi: \mathcal{H} (K) \to \mathcal{H}(K)\)
which generalizes the well-known symmetry of the HOMFLY polynomial:
 \(P_K(a,q) = P_K(a,q^{-1})\). So far, we are unable to
 explain either this symmetry or the differentials \(d_N \) \((N \leq
 0)\) in terms of \(\H\). However, there is one surprising exception.
The symmetry \(\phi\) should exchange \(d_1\) and \(d_{-1}\), so the
 conjecture implies that 
\begin{equation*}
H(\mathcal{H}(K),d_{-1}) \cong H_*(\mathcal{H}
(K),d_1) \cong \H_1(K).
\end{equation*}
The latter group is always isomorphic to \(\Q\), so we
expect that \(H(\mathcal{H}(K),d_{-1}) \cong \Q\) as well.
 In fact, we have

\begin{tm}
\label{Thm:d-1}
There is a spectral sequence \(E_k(-1)\) which starts at
\(\H(K)\) and converges to \(\Q\).   
\end{tm}

\noindent
 The grading behavior of the first differential \(d: \H(K) \to \H(K)
 \) matches the expected behavior of \(d_{-1}\), so again,
if the sequence converged after
 this differential, 
we would be in the situation of the conjecture. 
The construction of the sequence \(E_k(-1)\), while simple, is 
 unlike anything familiar from Khovanov homology.  It
 certainly behaves as if it should be dual to \(E_k(1)\) under the symmetry
 \(\phi\), but it is not clear how this duality might be realized.

\vskip0.05in

Although the KR-homologies are entirely combinatorial in nature, they
have been surprisingly difficult to compute. 
As an application of the theorems above, we  determine the
KR-homology of some  
simple knots. For example, combining Theorem~\ref{Thm:Limit} with
the main result of \cite{khthin} gives 

\begin{crrr} 
\label{Cor:TwoBridge}
If \(K\) is a two-bridge knot, then \(\H^{i,j,k}(K) = 0
  \) unless \(i+j+k = \sigma(K)\). 
\end{crrr}
\noindent
This condition is  similar to the usual notion of thinness in
Khovanov homology \cite{DBN,Khovanov2}. We call knots 
 which satisfy it {\it KR-thin}. The KR-homology of such a knot is completely
determined by its HOMFLY polynomial and signature. 
Many other small knots are KR-thin, and Theorems~\ref{Thm:d+} and
\ref{Thm:d-1} provide strong constraints on the homology of those
which are not. Using them, it is not difficult to determine the
KR-homology of all knots with 9 crossings or fewer.

\vskip0.05in 
The rest of the paper is organized as follows. In the  first three sections
we review (and in some cases, sharpen) various notions introduced by
Khovanov and Rozansky, starting with the 
definitions of the different  KR-homologies in  section~\ref{Sec:Defs}.
Section~\ref{Sec:Koszul} contains 
material  related to the theory of matrix factorizations,
 while section~\ref{Sec:Singular} describes the relation between
 KR-homology and the Murakami-Ohtsuki-Yamada state model. 
In sections~\ref{Sec:d+} and \ref{Sec:d-1} we construct the spectral
sequences of  Theorems~\ref{Thm:d+} and \ref{Thm:d-1}, respectively.
Finally, in section~\ref{Sec:Misc}, we  explain 
how these sequences can be applied 
to the problem of computing the KR-homology. 
\vskip0.05in
In writing, we have aimed to give a reasonably
self-contained treatment of  the KR-homology. In particular, we do not
assume that the reader is familiar with \cite{KRI,KRII}, and much of
the  the first three sections is devoted to a review of those
papers. The reasons
for this are both technical and expository. On the technical side, the
 proof of Theorem~\ref{Thm:d+} rests on results which
are very similar, but unfortunately not quite identical, to those in
\cite{KRI,KRII}. In order to give a complete treatment of these facts,
it seemed best to begin at the beginning. From the expository point of
view, we hope that readers who are unfamiliar with KR-homology will
find it convenient to have the definitions and normalization
conventions for the different theories housed under one roof. 
\vskip0.05in

\vskip0.05in
\noindent {\it Acknowledgements:} The author would like to thank
Dror Bar-Natan, 
Matt Hedden, Mikhail Khovanov, Marco Mackaay, Ciprian Manolescu, Peter Ozsv{\'a}th,
and Zolt{\'a}n Szab{\'o} for many helpful conversations during the
course of this work.

\section{Definitions}
\label{Sec:Defs}

Our goal in this section is to give a concise (but still
self-contained) definition of the various Khovanov-Rozansky
homologies. The material here is  all drawn
from \cite{KRI}, \cite{KRII}, and \cite{Gornik}, but we have slightly
modified some of the definitions.  In particular, the reader should be
aware that our grading conventions for the HOMFLY
homology are different from the ones introduced in \cite{KRII}.

\subsection{Matrix factorizations}
\label{SubSec:GMF} We begin by describing a class of  algebraic objects
known as matrix factorizations. These objects  first appeared 
in the context of algebraic geometry. Their application  to knot theory 
was one of the seminal advances of \cite{KRI}.

\begin{defn}
Suppose \(R\) is a commutative ring, and  that \(w \in R\). A {\em
  \(\Z\)--graded matrix factorization with potential} \(w\) consists of a free
 graded
  \(R\)-module \(C^*\) \((*\in \Z)\), together with a pair of {\em
  differentials} \(d_{\pm}: C^* \to C^{*\pm1}\) with the property
  that \((d_++d_-)^2 = w \cdot \text{Id}_{C}\).
\end{defn}
\noindent{\bf Remark:} We have included the phrase \(\Z\)--graded to
distinguish this definition  from the one used in 
in \cite{KRI} and \cite{KRII}, where matrix factorizations
 are \(\Z/2\)--graded. Unless we're trying to emphasize the
 distinction, we'll generally be careless and call a \(\Z\)--graded
 matrix factorization a matrix factorization. 
 \vskip0.05in 
 The  \(\Z\)--grading implies that the
condition \((d_++d_-)^2 = w \cdot \text{Id}_{C}\) is equivalent to
\begin{equation*}
d_+^2 = d_-^2 = 0 \quad \text{and} \quad d_+d_-+d_-d_+ = w \cdot \text{Id}_{C}. 
\end{equation*}
Thus a \(\Z\)--graded matrix factorization  gives rise
to two different chain complexes \(C^*_\pm\) with underlying group
\(C^*\) and differentials \(d_{\pm}\). If it happens that \(w=0\), we
get a third, \(\Z/2\)--graded
 chain complex structure \(C^*_{tot}\) on \(C^*\), with
differential \(d_{tot} = d_++d_-\). 

A morphism between two  matrix factorizations \(C^*\) and
\(D^*\) is a homomorphism of graded modules \(f: C^* \to D^*\)
which commutes with both differentials. We denote
the category of matrix
factorizations over a fixed ring \(R\) by \(GMF(R)\) and the
subcategory of factorizations with
fixed potential \(w\) by  \(GMF_w(R)\). 

%
%



The tensor product construction plays an important role in the
definition of the KR-homology.  If \(C^*\) and \(D^*\) are two 
matrix factorizations over \(R\), we endow the graded group
 \(C^* \otimes_R D^*\) with differentials
\(d_\pm\) defined by the requirement that \(d_-\) is the
differential on the chain complex \(C_-^* \otimes D_-^*\), and
similarly for \(d_+\). The reader can  easily verify

\begin{lem}
\label{Lem:Tensor}
If \(C^*\) and \(D^*\) are matrix factorizations  with
potentials \(w_1\) and \(w_2\), then \(C^*\otimes D^*\) is a  matrix
factorization with potential \(w_1+w_2\). 
\end{lem}


The final notion we need is that of a {\it complex of matrix
  factorizations with potential \(w\).} This is a \(\Z\)--graded chain
  complex defined over the category \(GMF_w(R)\). (Recall that the
  definition of a chain complex makes sense over any additive
  category.) More prosaically, such a complex consists of 
 a doubly graded group \(C^{*,*}\) equipped with 
 differentials 
\begin{equation*}
d_\pm:C^{i,j} \to C^{i\pm1,j} \quad \text{and} \quad d_v: C^{i,j}
 \to C^{i,j+1}
\end{equation*}
such that \((d_++d_-)^2 = w \cdot \text{Id}_C\), \(d_v^2 = 0\), and
  \(d_v\) commutes with both \(d_+\) and \(d_-\). Often, it is more
  convenient to have \(d_v\) {anticommute} with \(d_\pm\). This
  can be arranged by replacing \(d_v\) with \((-1)^id_v\). 

It is  helpful to think of \(C^{*,*}\) as being a sort of
generalized double complex. We envision the group \(C^{i,j}\) as
sitting over the point \((i,j)\) in the \(xy\)--plane, so that the
differentials \(d_\pm\) carry us
 one unit to the right and left, respectively, and \(d_v\) carries us
 one unit up. In keeping with this picture, we refer to \(i\) and
 \(j\) as the {\it horizontal} and {\it vertical} gradings on
 \(C^{*,*}\), and denote them by  \(\gr_h\) and \(\gr_v\), respectively.
In addition to these gradings, it 
is also natural to consider the quantities 
\(\gr_{\pm} = gr_v\pm \gr_h\), which are the  total gradings on the
double complexes \(C^{*,*}_\pm\).

In the sequel, we will frequently take the tensor product of complexes
of matrix factorizations. Since we know how to take  tensor products
of chain complexes and  of matrix factorizations, it's clear how this
is to be done.  
From Lemma~\ref{Lem:Tensor}, we see that the tensor product of a
complex of matrix factorizations with potential \(w_1\) with a complex
of matrix factorizations with potential \(w_2\) is 
a complex of matrix factorizations with potential \(w_1+w_2\). 


 \subsection{Tangle diagrams}
\begin{figure}
\includegraphics{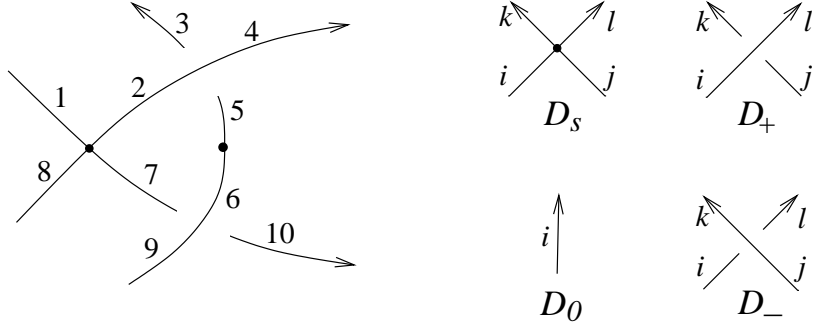}
\caption{\label{Fig:Tangles}Some singular tangles. The left-hand side
  shows a labeled singular tangle, including a mark and a crossing of
  each type. On the right hand side are diagrams of the four
  elementary tangles}
\end{figure}

KR-homology is most naturally defined in the context
 of {\it singular oriented planar
tangles}. These are oriented planar diagrams which in addition to the usual
over- and undercrossings may also contain some singular points, as
illustrated in Figure~\ref{Fig:Tangles}. (In the notation of
\cite{KRI} and \cite{KRII}, singular points correspond to wide edges.) 
From now on, we will just refer to them as {\it tangle diagrams.}

More formally, a tangle diagram is an oriented planar graph, all of
whose vertices have valence \(1\), \(2\), or \(4\).
 The 4-valent vertices are {\it crossings}, and come with an
 additional decoration indicating whether they are {\it positive},
 {\it negative}, or {\it singular}, as represented by the  diagrams
 \(D_+\), \(D_-\), and \(D_s\) shown in the figure. Bivalent vertices
 are called {\it marks}, and must have one incoming and one outgoing
 edge. Univalent vertices are  called {\it free ends}. An edge
adjacent to such an end  is called  {\it external}; all other edges
are {\it internal}. 
The {\it connected components} of a diagram are the connected
components of the underlying graph (not the
connected components of the associated tangle). A component with no
free ends is {\it closed}; other components are {\it open}. 
We keep track of the edges in a tangle diagram by labeling them
by integers \(1,2, \ldots, n\), where \(n\) is the number of edges in
the diagram. A free end is identified by the label of its adjacent
edge. 

We now describe some operations for building new tangle diagrams out of
old ones. First, if \(D_1\) and \(D_2\) are two tangle diagrams, we
can take their disjoint union \(D_1\sqcup D_2\). Second, if \(i\) is
an edge of \(D\), we can form an new diagram \(D(i)\) by inserting a
bivalent vertex into \(i\). We can also perform the inverse
operation, which is known as {\it mark removal}. Finally, 
suppose that \(D\) is a tangle diagram with incoming and outgoing 
free ends labeled \(i\) and \(j\), and \(i\) is {\it adjacent} to \(j\) in
the sense that they can be isotoped onto each other without hitting
the rest of the graph. Then we can form a new diagram \(D\vert_{i=j}\) 
by identifying \(i\) and \(j\) to form a single bivalent vertex. 
Any tangle diagram can be built up from the elementary  diagrams
\(D_+\), \(D_-\), \(D_s\), and \(D_0\) shown in the figure
by the operations of disjoint union, identifying free ends, and mark
removal. 

\subsection{Edge rings}
\label{SubSec:Rings}
Suppose \(D\) is  a tangle diagram with edges labeled \(1,\ldots,n\),
 and let \(R'(D)\) be the ring \(\Q[X_1,\ldots,X_n]\). 
 To an internal vertex \(v\) of \(D\), we assign a linear relator 
\(\rho(v)\) in \(R'(D)\). \(\rho(v)\)  is the sum of the variables corresponding
 to outgoing edges of \(v\) minus the sum of the variables
 corresponding to ingoing edges. In other words, if \(v\)
 is a mark with incoming edge \(i\) and
 outgoing edge  \(j\), \(\rho(v) = X_j-X_i\), and if
 \(v\) is a crossing with incoming edges \(i\) and \(j\) and outgoing
 edges \(k\) and \(l\), \(\rho(v) = X_k+X_l-X_i-X_j\).
\begin{defn}
 The {\it  edge ring} \(R(D)\) is  the graded ring
\(R'(D)/(\rho(v_j))\) where \(j\) runs over all internal vertices
of \(D\).  The grading on \(R(D)\) is denoted by \(q\); it 
is determined by the requirement  that \(q(X_i) = 2\) for all \(i\).
\end{defn}

The edge ring behaves nicely under the operations of
disjoint union, mark removal, and identifying free ends. The reader
can easily verify that
\begin{align*}
R(D_1 \sqcup D_2) & \cong R(D_1) \otimes_{\Q} R(D_2) \\
R(D(i)) &\cong R(D) \\
R(D\vert_{i=j}) & \cong R(D)\vert_{X_i=X_j}.
\end{align*}
More generally, suppose that  \(D\) is obtained from
diagrams \(D_1\) and \(D_2\)  by first taking their disjoint union
 and then identifying ends \((i_1,i_2,\ldots i_m)\)
of \(D_1\) with ends \((j_1,j_2,\ldots j_m)\) of \(D_2\). Applying the
relations above, we see that
\begin{equation*}
R(D) \cong R(D_1) \otimes_{\Q[y_1,\ldots,y_m]} R(D_2)
\end{equation*}
where \(y_k\) acts as \(X_{i_{k}}\) on \(R(D_1)\) and \(X_{j_k}\) on
\(R(D_2)\).

\begin{lem}
\label{Lem:EdgeRing} 
Suppose \(D\) is a tangle diagram with \(V\) internal
vertices, \(E\) edges, and \(C\)
closed components. Then \(R(D)\) is isomorphic to a polynomial ring on 
\(E-V+C\) variables. 
\end{lem}

\begin{proof} It is enough to prove the claim when
  \(D\) is connected; the general result then follows from 
the tensor product formula for disjoint unions. If \(D\) is open and
connected, the statement amounts to saying that the relations
\(\rho(v_j)\) 
are linearly independent in the vector space spanned by the \(X_i\). 
Suppose \(\sum_j\alpha_j \rho(v_j) = 0 \). Then for any internal edge
\(i\), the coefficient of \(X_i\) in the sum vanishes, which means
that the values of \(\alpha\) on its two ends must be
equal. Since \(D\) is connected, it follows that all of the
\(\alpha_j\)'s are equal. If \(D\) is open, considering an external
edge shows that \(\alpha_j \equiv 0\), while if \(D\) is closed
and connected,
 there is a unique linear relation between the
\(\rho(v_j)\). 
\end{proof}

We will also use  two  subrings of the edge
ring. These are the {\it external ring} \(R_e(D)\), which is the
subring  generated by the \(X_j\), where \(j\) runs over the
external edges of \(D\), and the {\it reduced ring} \(R_{r}(D)\) which
is generated by the differences \(U_{ij}=X_i-X_j\), where \(i\) and
\(j\) run over all edges of \(D\). 

More explicitly, the external ring can be described as follows. 
We assign a sign \(\epsilon_j\) to each free end of \(D\) 
according to the rule that
\(\epsilon_j=1\) if \(j\) is an outgoing end, and \(\epsilon_j=-1\) if
it is incoming. If \(C\) is a connected component of \(D\), we assign
to it the polynomial \(\rho(C) =  \sum \epsilon_jX_j\), where
\(j\) runs over the free ends of \(C\). (Note that if \(C\) is an
elementary diagram, this reduces to the previous definition.) Then we
have

\begin{lem} \(R_e(D) \cong \Q[X_j]/(\rho(C))\), where \(j\) runs over
  the external edges of \(D\) and \(C\) runs over the set of connected
  components of \(D\). 
\end{lem}

\begin{proof} Consider the vector space \(V = \langle X_i \ts \vert
 \ts i \text{ is an edge of \(D\)}\rangle\), along with its linear subspaces
 \(V_e =\langle X_j \ts \vert
 \ts j \text{ is an exterior edge of \(D\)}\rangle\) and \(V_c = 
\langle \rho(c) \ts \vert
 \ts c \text{ is a crossing  of \(D\)}\rangle\). Let \(I \subset \Q[X_j]
 \) be the ideal generated by \(V_e \cap V_c\). Then
 \(R_e(D) \cong \Q[X_j]/I\), so it suffices to show that \(V_e \cap
 V_c \) is generated by the elements
\(\rho(C)\), where \(C\) runs over the components of \(D\). Now if
\( \rho = \sum _{c} \alpha_c \rho(c) \in V_e\), the component of \(\rho\) along
 each internal edge must vanish, which means that \(\alpha_c\) has the
 same value at the two ends of the edge. Thus
 \(\alpha_c\) is constant on connected components, and the claim is proved.  
\end{proof}

\noindent The
edge ring and the reduced ring are related by

\begin{lem}
\label{Lem:ReducedRing}
\(R(D) \cong R_r(D)[x].\)
\end{lem}

\begin{proof}
Let \(R'_r(D)\) be the subring of \(R'(D) \) generated
by the \(X_i-X_j\). Then the map which sends \(x\) to \(X_1\)
defines an isomorphism from  \( R_r'(D)[x]\) to \(R'(D)\). Since the relations
\(\rho(v_j)\) are all contained in \(R_r'(D)\), this descends to an
isomorphism  \( R_r(D)[x] \cong R(D)\). 
\end{proof}

 \subsection{The KR-complex}
\label{SubSec:Complex}

The key step in the definition of KR-homology is a process
which assigns to a tangle a triply-graded complex of matrix factorizations.
More precisely, let 
\(D\) be a tangle diagram, and fix as an auxiliary
parameter a polynomial \(p(x) \in \Q[x]\). Then the {\it KR-complex}
\(C_p(D)\) associated to the pair \((D,p)\) is a complex of matrix
factorizations over the ring \(R(D)\) with potential 
\begin{equation*}
w_p(D) = \sum_{j} \epsilon_j p(X_j),
\end{equation*}
where the sum runs over the external edges of \(D\).

\(C_p(D)\) is a graded module over the graded ring
\(R(D)\). This grading corresponds to the power of \(q\) in the HOMFLY
polynomial, and will be referred to as the {\it \(q\)--grading.}  The
other two gradings on \(C_p(D)\) are the homological gradings
\(\gr_h\) and \(\gr_v\) coming from its structure as a complex of
matrix factorizations. The differentials on \(C_p(D)\) interact with
the \(q\)-grading as follows: \(d_v\) preserves the \(q\)-grading,
while \(d_+\) increases it by \(2\). \(d_-\) is usually not homogenous
with respect to the \(q\)-grading, but if
\(p(x)=x^n\), \(d_-\) raises the \(q\)-grading by \(2n-2\). 
We summarize our conventions regarding the various gradings in the following
\begin{defn}
\label{Def:Gradings}
We say that \(x
\in C^{i,j,k}_p(D)\) if \(x\) is homogenous with respect to
all three gradings, and  \((i,j,k) = (q(x),
2\gr_h(x), 2\gr_v(x))\). With respect to this grading, 
\(d_v\) is homogenous of degree
\((0,0,2)\) and \(d_{+}\) is homogenous of degree \((2,2,0)\). If
\(p(x) = x^{n}\), \(d_-\) is homogenous of degree \((2n-2,-2,0)\). 
\end{defn}
\noindent{\bf Remark:} At first sight, the fact that we have chosen
to double the homological gradings may seem rather strange. In fact,
there are two good reasons for this choice of normalization. First, as
we will explain in section~\ref{Sec:Singular},
the quantity \(2 \gr_h\) is naturally related to the power of \(a\) in
the HOMFLY polynomial. Second, with this normalization, \(i,j\), and
\(k\) all have the same parity when \(D\) is an ordinary diagram ({\it
  i.e.} one with 
no singular crossings.) 

\subsection{Elementary tangles}
\label{SubSec:ETangle}

Before  defining the KR-complex in general, we describe  it for the 
elementary diagrams \(D_s\), \(D_+\), and \(D_-\) shown in
Figure~\ref{Fig:Tangles}. In each case, \(C_p(D)\) will be
a complex of matrix factorizations over the ring  
\begin{equation*}
R  =
\Q[X_i,X_j,X_k,X_l]/(X_k+X_l-X_i-X_j)\cong \Q[X_i,X_j,X_k]
\end{equation*}
 with
potential 
\begin{align*}
W_p(X_i,X_j,X_k,X_l) & = p(X_k)+p(X_l) - p(X_i)-p(X_j) \\ & = p(X_k) + 
p(X_i+X_j-X_k) - p(X_i) - p(X_j).
\end{align*}

The complex \(C_p(D_s)\) is a free \(R\)--module of
rank \(2\). Since the potential is nonvanishing, the map \(d_+\) must take one
copy of \(R\) to the other. It is given by multiplication by
 \(X_kX_l-X_iX_j\), which is equal (in \(R\)) to 
\(-(X_k-X_i)(X_k-X_j)\). The  map  \(d_-\)  
 takes the first copy of
\(R\) back to the second, and must be given by multiplication by  
\begin{equation*}
p_{ij} = -W_p/{(X_k-X_i)(X_k-X_j)}.
\end{equation*}
Note  that if we substitute either \(X_k=X_i\) or
\(X_k=X_j\) into \(W_p\), the result vanishes, so the quotient
\(p_{ij}\) really is an element of \( R\). 
Finally, the two copies of \(R\) have the same vertical grading, so
\(d_v\) is necessarily trivial. 
More succinctly, we can represent \(C_p(D_s)\)  by the diagram
$$\xymatrixcolsep{7pc} C_p(D_s) = 
\xymatrix{ R\{1,-2,0\}  \ar@<0.4ex>[r]^{-(X_k-X_i)(X_k-X_j)}&
  R\{-1,0,0\} \ar@<0.4ex>[l]^{p_{ij}}}.$$
Following \cite{KRI}, \cite{KRII}, we use the  notation \(R\{i,j,k\}\) to
indicate a free \(R\)--module of rank one  with 
gradings shifted so that if \(1\) is a
generator of \(R\{i,j,k\}\), then \(1\in C^{i,j,k}_p(D_s)\).

Using the same notation, the complexes \(C_p(D_+)\) and \(C_p(D_-)\) 
 are  given by diagrams
$$\xymatrixcolsep{2.0pc} 
\xymatrixrowsep{1.7pc}
\xymatrix{ 
& R\{0,-2,0\}  \ar@<0.4ex>[rrr]^{(X_k-X_i)}& & &
  R\{0,0,0\} \ar@<0.4ex>[lll]^{p_{i}} \\
C_p(D_+) = & & & & \\
& R\{2,-2,-2\}  \ar@<0.4ex>[rrr]^{-(X_k-X_i)(X_k-X_j)} 
\ar[uu]^{(X_j-X_k)}& & &
  R\{0,0,-2\} \ar@<0.4ex>[lll]^{p_{ij}} \ar[uu]_{1}} $$
 and 
$$\xymatrixcolsep{2pc}  
\xymatrixrowsep{1.7pc}
\xymatrix{ 
& R\{0,-2,2\}  \ar@<0.4ex>[rrr]^{-(X_k-X_i)(X_k-X_j)}& & &
  R\{-2,0,2\} \ar@<0.4ex>[lll]^{p_{ij}} \\
C_p(D_-) = & & & & \\
& R\{0,-2,0\}  \ar@<0.4ex>[rrr]^{(X_k-X_i)} 
\ar[uu]^{1}& &  &
  R\{0,0,0\}. \ar@<0.4ex>[lll]^{p_{i}} \ar[uu]_{(X_j-X_k)}} $$
Here \(p_i=W_p/(X_k-X_i)\),
and the vertical arrows represent components of
the map  \(d_v\). 
The reader
can easily verify that in all three complexes, \(d_+\) and   \(d_v\)  are
homogenous of degree \((2,2,0)\) and \((0,0,2)\), respectively.

\subsection{General tangles}
For an arbitrary tangle diagram \(D\), \(C_p(D)\) 
is defined to be a tensor product of smaller complexes, one for each
crossing in \(D\). More precisely, if  \(c\) is a crossing of \(D\),
let \(D_c\) be the subdiagram composed of the four edges of \(D\)
adjacent to \(c\). \(D_c\) is an elementary diagram, so the complex
\(C_p(D_c)\) was defined in the previous section. 
 It is a complex of matrix factorizations over the ring
\(R_c = \Q[X_i,X_j,X_k,X_l]/(X_k+X_l-X_i-X_j)\) with potential
\(w_p(c) = p(X_k)+p(X_l)-p(X_i)-p(X_j)\).

Next, we consider the complex  \(C_p(D_c) \otimes_{R_c} R(D)\), which
 is obtained by replacing
each copy of \(R_c\) in \(C_p(D_c)\) with a copy of \(R(D)\). 
It is a complex of matrix factorizations over \(R(D)\). 
The global complex  \(C_p(D)\) is defined to be the tensor product
 over the ring \(R(D)\)
\begin{equation*}
C_p(D) = \bigotimes_{c} \Big( C_p(D_c) \otimes_{R_c} R(D)\Big). 
\end{equation*}
where the product runs over all crossings of \(D\). In particular, if
there are no crossings, \(C_p(D) = R(D)\).

We can now  verify that \(C_p(D)\) has the properties advertised
in section~\ref{SubSec:Complex}. First, it is clearly defined over the
ring \(R(D)\). 
 Second, it is easy to see that \(\sum_cw_p(c) =
w_p(D)\), so it follows from Lemma~\ref{Lem:Tensor} that \(C_p(D)\)
has potential \(w_p(D)\). Third, the differentials on each individual
factor satisfy the grading conventions established in
definition~\ref{Def:Gradings}, so the same is true for \(C_p(D)\). 

An important (indeed, the defining) property of \(C_p(D)\) is that it is
 {\it local} in the following sense:

\begin{lem}
\label{Lem:Gluing}
Suppose \(D\) is obtained from
diagrams \(D_1\) and \(D_2\)  by first taking their disjoint union
 and then identifying ends \((i_1,i_2,\ldots i_m)\)
of \(D_1\) with ends \((j_1,j_2,\ldots j_m)\) of \(D_2\). Then
\begin{equation*}
C_p(D) \cong C_p(D_1) \otimes_{\Q[y_1,\ldots,y_m]} C_p(D_2)
\end{equation*}
where \(y_k\) acts as \(X_{i_{k}}\) on \(C_p(D_1)\) and \(X_{j_k}\) on
\(C_p(D_2)\). 
\end{lem}
 
\begin{proof}
This follows from the fact that the set of crossings for \(D\) is the
union of the sets of crossings for \(D_1\) and \(D_2\), together with the
relation
\begin{equation*}
R(D) \cong R(D_1) \otimes_{\Q[y_1,\ldots,y_m]} R(D_2)
\end{equation*}
observed in section~\ref{SubSec:Rings}. 
\end{proof}

\subsection{The HOMFLY homology}
\label{SubSec:MiddleH}
We  now define the various KR--homologies, starting with the HOMFLY
homology of
\cite{KRII}. There are
 several  ways of normalizing this invariant, all of which
 contain the same information. In addition to the {\it reduced theory}
 used in the introduction, there is also an {\it unreduced
 theory} which appears naturally in the context of the \(sl(N)\)
 homology. We start with a third variant, which interpolates
 between these two and is closest to the version of the theory
 described in \cite{KRII}. 

For the next few sections, we assume that \(L\) is an oriented link in
\(S^3\), and that \(L\) is represented by a connected tangle diagram
\(D\) which is the closure of a braid. (The restriction that \(D\) be
connected is simply for ease of exposition. The necessary
modifications for disconnected diagrams are described in
section~\ref{SubSec:Disconnected}.) 

\begin{defn} 
\label{Def:MiddleH}
The middle HOMFLY homology of \(L\) is the group
\begin{equation*}
\Hm(L)= H(H(C_p(D),d_+),d_v^*)\{-w+b, w+b-1, w-b+1\}, 
\end{equation*}
where 
\(w\) and \(b\) are the writhe and number of strands of the braid
diagram \(D\). 
\end{defn}
\noindent {\bf Remarks:} There are several aspects of this definition
which are worth pointing out. First, observe that 
we have taken homology {\it twice}: first with respect to
 \(d_+\), and then with respect to \(d_v^*\), which is
the map induced on \(H({C}_p(D),d_+)\) by \(d_v\). Second, note
that \(d_+\) and \(d_v\) are homogenous with respect to 
all three gradings, so the 
 triple grading on \(C_p(D)\) descends to a triple grading on
 \(H(L)\). 
 Finally, since \(d_-\) does not appear in the definition, \(\Hm(L)\)
is independent of the parameter  \(p\). 

In \cite{KRII}, Khovanov and Rozansky proved 
\begin{thrm}
\label{Thm:Invariance}
 \cite{KRII} \(\Hm(L)\) is an invariant of \(L\). 
\end{thrm}

{\it A priori}, there is nothing stopping us 
from considering the
homology  \(H(H(C_p(D),d_+),d_v^*)\) for an arbitary 
diagram \(D\) representing \(L\), but the restriction to diagrams 
which are braid closures plays an important role in the proof of
Theorem~\ref{Thm:Invariance}. Indeed, Khovanov and Rozansky
prove the invariance of \(\Hm(L)\) under {braidlike}
Reidemeister moves and then use the fact that any two braid diagrams
of  \(L\) are related by such moves to conclude that 
\(H(L)\) is a link invariant. 

The second major result of \cite{KRII} is the relation between
\(\Hm(L)\) and the HOMFLY polynomial:
\begin{thrm}
\label{Thm:Chi} 
\cite{KRII} For any  \(L \subset S^3\), we have  
\begin{equation*}
\sum_{i,j,k} (-1)^{(k-j)/2}a^jq^i  \dim \Hm^{i,j,k}(L)
= - \frac{P(L)}{q-q^{-1}}. 
\end{equation*}
\end{thrm}
\noindent Here, both sides of the equation should be interpreted as 
 Laurent series in \(q\). 

\subsection{Reduced and unreduced complexes}
\label{SubSec:UnRed}
If \(i\) is an edge of \(D\), we define the reduced KR--complex
\(\Cr_p(D,i)\) to be the quotient 
\( C_p(D)/(X_i)\). In \cite{KRII}, Khovanov and Rozansky observe that 
when \(p=0\), this definition is actually
independent of \(i\). 
To see this, recall that \(C_0(D)\) is a direct sum of copies of
\(R(D)\). Define \(C_r(D) \subset C_0(D)\) to be the subgroup obtained
by replacing each copy of \(R(D)\) with a copy of the reduced ring
\(R_r(D)\). Inspecting the coefficients of  \(d_+\) and \(d_v\) in 
  \(C_0(D_s)\), \(C_0(D_+)\), and \(C_0(D_-)\), we see that
they  are all contained in \(R_r(D)\). It follows that \(C_r(D)\) is a
subcomplex of \(C_0(D)\). 

\begin{lem}
\label{Lem:GlobalReduced}
In the category \(GMF(\Q)\), there are  isomorphisms 
\(C_0(D) \cong C_r(D)\otimes_\Q \Q[x]\) and  \(\Cr_0(D,i)\cong
C_r(D)\). 
\end{lem}

\begin{proof}
The first claim follows immediately from
Lemma~\ref{Lem:ReducedRing}. For the second, consider the map
\(\phi: R_r(D)  \to R(D)/(X_i)\) which is the composition of the
inclusion \(R_r(D) \to R(D)\) and the projection \(R(D) \to
R(D)/(X_i)\). It's easy to see that \(\phi\) is an isomorphism of
vector spaces.  Since \(C_0(D)\) is free over \(R(D)\), 
 the induced map \(\phi:C_r(D) \to C_0(D,i)\) is also an isomorphism. 
\end{proof}

\begin{defn}
\label{Def:Reduced}
The reduced HOMFLY homology \(\H(L)\) is
defined to be
\begin{equation*}
\H(L) = H(H(C_r(D),d_+),d_v^*)\{-w+b-1,w+b-1,w-b+1\}  
\end{equation*}
where as before, \(w\) and \(b\) are the writhe and number of strands
in the braid diagram \(D\). 
\end{defn}

From the first part of 
Lemma~\ref{Lem:GlobalReduced}, we see  that 
 \(\Hm(L) \cong \H(L)\otimes_\Q \Q[x]\). It follows that the graded
 Euler characteristic of \(\H(L)\) is given by the HOMFLY
 polynomial:
\begin{equation*}
\sum_{i,j,k} (-1)^{(k-j)/2}a^jq^i  \dim \H^{i,j,k}(L)
= P(L). 
\end{equation*}

There is also  an {unreduced} version of the KR--complex.
If \(i\) is an edge of \(D\), we let \(U_p(i)\) be the matrix factorization
$$\xymatrixcolsep{5pc} C_p(D_s) = 
\xymatrix{ \Q[X_i]\{0,-2,0\}  \ar@<0.4ex>[r]^{0}&
  \Q[X_i] \{0,0,0\} \ar@<0.4ex>[l]^{p'(X_i)}}.$$
 The unreduced complex \(\Cu_p(D,i)\) is defined to be
\(C_p(D) \otimes_{\Q[X_i]} U_p(i)\).
\begin{defn}
\label{Def:Unreduced}
The unreduced HOMFLY homology \(\Hu(L)\) is given by
\begin{equation*}
\Hu(L) = H(H(\Cu_p(D,i),d_+),d_v^*)\{-w+b,w+b,w-b\}  
\end{equation*}
where  \(w\) and \(b\) are the writhe and number of strands
in the braid diagram \(D\). 
\end{defn}

Since both \(d_+\) and \(d_v\) are trivial on \(U_p(i)\), we see that
 \(\Hu(L) \cong \Hm(L) \otimes H^*(S^1)\). Its graded Euler
 characteristic  is the unnormalized HOMFLY
polynomial of \(L\):
\begin{equation*}
 \sum_{i,j,k} (-1)^{(k-j)/2}a^jq^i  \dim 
\Hu^{i,j,k}(L) = \frac{a-a^{-1}}{q-q^{-1}} P(L) = \wt{P}(L).
\end{equation*}
\noindent{\bf Remark:} The quantity \(w+b\) always has the same
parity as the number of components of  \(L\). Since all the grading
shifts in the complexes \(C_p(D_+)\) and \(C_p(D_-)\) are even, it
follows that all three gradings of \(\Hu(L)\) have the same parity as
the number of components of \(L\), and  all three gradings of
\(\H(L)\) have the opposite parity. 

As an example, we describe \(H\), \(\H\), and \(\Hu\) for the
unknot. The unknot can be represented by a braid diagram \(D\)
consisting of a single edge (labeled \(1\)), a single mark, and no
crossings. The relation associated to the mark is \(X_1-X_1 = 0\), so
\(R(D) = \Q[X_1]/(0) \cong \Q[X_1]\), and \(R_r(D) \cong \Q\). 
Since there are no crossings, \(C_p(D) \cong R(D)\). 
It follows that 
\(H(U) \cong \Q[X]\), where \(1 \in \Q[X]\) has triple grading
\((1,0,0)\); \(\H(U) \cong \Q\), with triple grading \((0,0,0)\); and
\(\Hu(U) \cong \Q[X] \oplus \Q[X]\), where the generators have
gradings \((1,1,-1)\) and \((1,-1,-1)\).

\subsection{The \(sl(N)\) homologies} 
\label{SubSec:SlN}


To define the
KR--homologies corresponding to the \(sl(N)\) polynomial, we
add the differential \(d_-\) into the mix. Suppose that 
 \(D\) is a connected tangle diagram --- not necessarily in braid form
 --- representing the link \(L\). Then \(D\) is closed, so the potential
  \(w_p(D) = 0 \), and
 the differential \(d_{tot} = d_++d_-\)  makes 
 \(\Cr_p(D,i)\) and \(\Cu_p(D,i)\) into chain complexes. 

\begin{defn}
For \(p(x) \in \Q[x]\), the reduced and unreduced
\(p\)--homologies are defined by
\begin{align*}
 \H_p(L,i) & = H(H(\Cr_p(D,i),d_{tot}),d_v^*) \\
\Hu_p(L) & = H(H(\Cu_p(D,i),d_{tot}),d_v^*) 
\end{align*} 
\end{defn}

\noindent
When \(p(x) = x^{N+1}\), this definition was
introduced  by Khovanov and Rozansky  in \cite{KRI}. The
fact that the definition is interesting for other values of \(p\) was
observed by Gornik  \cite{Gornik}. 

 For the definition to make sense, we should check that 
\(\Hu_p(L)\) depends only on \(L\), and not on the choice of the
diagram \(D\) or the marked edge \(i\). This is done in
section~\ref{Sec:d+}. Similarly, the reduced homology \(\H_p(L,i)\)
depends only on \(L\) and the component of \(L\) containing \(i\).  
Unlike the HOMFLY homology, \(\H_p(L,i)\) really does depend on
the marked component.
 However, in the special case when \(L=K\) is a knot,  there is  
only one component to choose from, so it makes sense to talk about the
 reduced homology \(\H_p(K)\). 

Next, we consider the grading on these homology groups. 
For a general polynomial \(p\), 
 \(d_{tot}\) will not be homogenous with respect to any linear combination of
the gradings \(q\) and \(\gr_h\) on \(C_p(D)\), so \(\H_p\) and \(\Hu_p\)
 will  have only the single grading coming from \(\gr_v\). 
However, when \(p(x)=X^{N+1}\),
 \(d_{tot}\) is homogenous with respect
to the grading 
\begin{equation*}
\gr_N = q + (N-1) \gr_h = i+\frac{N-1}{2}j, 
\end{equation*}
 so we can view  \(\H_p(L,i)\) and \( \Hu_p(L)\) as being
 doubly graded, with gradings \((\gr_N,\gr_v)\). An additional global
 shift is needed to make the first grading into a link invariant. We
 put
\begin{align*}
\H_N(L,i) & = \H_{x^{N+1}}(D,i)\{(N-1)w,0\} \\
\Hu_N(L) & = \Hu_{x^{N+1}}(D)\{(N-1)w,0\},
\end{align*}
where \(w\) is the writhe of the diagram \(D\). 
In section~\ref{SubSec:DefCom}, we verify that
\(\H_N\) and \(\Hu_N\) are the \(sl(N)\) homology groups defined by
Khovanov and Rozansky in \cite{KRI}. 
Their graded Euler characteristic is given by the
\(sl(N)\) polynomial:
\begin{thrm}
\cite{KRI} \(\Hu_N(L)\) is an invariant of the link \(L\), while
\(\H(L,i)\) is an invariant of the link \(L\) and the marked component
\(i\). They satisfy
\begin{align*}
\sum_{I,J} (-1)^Jq^I \dim \H_N^{I,J}(L,i) = P_L(q^N,q) \\
 \sum_{I,J} (-1)^Jq^I \dim \Hu_N^{I,J}(L) = \wt{P}_L(q^N,q) 
\end{align*}
\end{thrm}

As an example, we again consider the homology of the unknot. The
reduced complex \(\Cr_P(U) \cong \Q[X_1]/(X_1) \cong \Q\), so
\(\H_p(U) \cong \Q\), for any \(p\). The complex \(\Cu_p(U)\) is more
complicated. It is composed of two copies of \(\Q[X]\), situated in
gradings \((0,0,0)\) and \((0,-2,0)\). The differential \(d_-\)
takes a generator of the first summand to \(p'(X)\) times the
generator of the second. Thus \(\Hu_p(U) \cong \Q[X]/(p'(X))\)
supported in homological grading \(0\). When \(p(X) = X^{N+1}\), we
see that \(\Hu_N(U) \cong \Q[X]/(X^N)\). The generator \(1 \in 
\Q[X]/(X^N)\) has polynomial grading \(\gr_N = 1-N\). 

\subsection{Disconnected diagrams} 
\label{SubSec:Disconnected}
We conclude our discussion of KR--homology by describing what
happens when  the diagram \(D\) is disconnected. In this case, 
we must modify the definition
of the complexes \(\Cu_p(D)\) and \(\Cr_p(D)\). The unreduced complex
\(\Cu_p(D)\) is the tensor product
\begin{equation*}
\Cu_p(D) = \bigotimes_j \Cu_p(D_j,i_j)
\end{equation*}
where \(j\) runs over the connected components of 
\(D\). The definition requires that we specify a collection of edges
\(i_j\) --- one for each component of \(D\). 
In section~\ref{SubSec:DefCom}, we will show that 
\(\Cu_p(D)\) 
is essentially independent of the choice of \(i_j\). 
From Lemma~\ref{Lem:Gluing}, we see that
\begin{equation*}
C_p(D) = \bigotimes_j C_p(D_j),
\end{equation*}
 so from the point of view of the HOMFLY homology, the extra factors
 \(\otimes U_p(i_j)\)  just 
add a factor of \(H^*(S^1)\) for each component of \(D\). 

To define the reduced KR--complex, assume that the special marked edge
\(i\) is in the component \(D_1\). Then
\begin{equation*}
\Cr_p(D) = \Cr_p(D_1,i) \otimes \bigotimes_{j>1} \Cu_p(D_j). 
\end{equation*}
The definitions of the various KR--homologies now proceed exactly as
they did in the case when \(D\) had only one component.

%

\section{Matrix Factorizations}
\label{Sec:Koszul}

In this section, we develop some ideas about \(\Z\)--graded
matrix factorizations which will be needed in the rest of the paper. The
main difficulty with such factorizations, as compared to the
\(\Z/2\)--graded factorizations used in \cite{KRI} and \cite{KRII}, is
that they lack a good notion of homotopy equivalence. Our first task is
to develop an appropriate substitute --- the notion of a
quasi-isomorphism. After that, we discuss the class of Koszul
factorizations introduced by Khovanov and Rozansky in
\cite{KRI} and adapt some of their results to the \(\Z\)--graded
context. We conclude  by verifying
that the definitions of the various KR--groups given in section 2
coincide with the original definitions in \cite{KRI} and
\cite{KRII}. 

\subsection{Positive homology} Given a \(\Z\)--graded  matrix factorization
 \(C^*\), we define its   {\it positive
    homology} to be the group 
\begin{equation*}
H^+(C^*) = H(C^*,d_+). 
\end{equation*}
If it happens that \(C^* = C_p(D)\),
we abbreviate still further and  write  \(H^+(D)\)  in place of
  \(H^+(C_p(D))\), and similarly for   \(\Hu^+(D) =
    H^+(\Cu_p(D))\) and  \(\H^+(D) = H^+(\Cr_p(D))\). The operation of
    taking the positive homology gives a covariant functor
    \(\mathcal{H}^+\) from the category \(GMF(R)\) to the category of
    graded \(R\)--modules. This  naturally extends to a functor from 
\(Kom(GMF_w(R))\) to \(Kom (R)\). For example,
 Definition~\ref{Def:Unreduced} can be rewritten as
\begin{equation*}
\Hu(L) = H(\Hu^+(D),d_v^*)\{-w+b,w+b,w-b\}
\end{equation*}
in this notation.

 When the factorization has potential \(0\), we can
say more:
\begin{lem} There are functors 
\begin{align*}
\mathcal{H}^+&: GMF_0(R) \to Kom(R) \\
\mathcal{H}^+&: Kom(GMF_0(R)) \to Kom(Kom(R)). 
\end{align*}
\end{lem}

\begin{proof} If \(C^*\) has zero potential, the differentials \(d_-\) and
  \(d_+\) anticommute. The induced map \(d_-^*:H^+(C^*) \to
  H^+(C^*)\) makes \(H^+(C^*)\) into a chain complex. 
\end{proof}

\subsection{Quasi-Isomorphisms}
Roughly speaking, we want to think of two matrix factorizations
as being equivalent if their positive homologies are isomorphic as chain
complexes. When the potential is nonzero, however, the positive
homology isn't a chain complex. 
To get around this problem, we adopt the following

\begin{defn}
Suppose \(C^*,D^*\) are objects of \(GMF_w(R)\), and that \(f:C_+^*
\to D_+^*\) is a chain map. We say that \(f\) is a {\it
  quasi-isomorphism} if for every object \(E^*\) of \(GMF_{-w}(R)\)
the induced map 
\((f\otimes1)^*:H^+(C^*\otimes E^*) \to H^+(D^*\otimes E^*)\)
 is an isomorphism which commutes with \(d_-^*\). More generally, we
 say that \(C^*\) and \(D^*\) are {\it quasi-isomorphic} and write
 \(C^* \sim D^*\) if they can
 be joined  by a chain of quasi-isomorphisms. 
\end{defn}

Note that \(f\) is not required to be a morphism of matrix
factorizations, but only a map on the positive chain complexes which 
``looks like'' such a morphism when we pass
to homology, in the sense that it commutes with \(d_-^*\).

In practice, many of the quasi-isomorphisms we will consider do arise
 as morphisms.
\begin{defn} Suppose \(C^*,D^*\) are objects of \(GMF_w(R)\), and that
  \(f: C^* \to D^*\) is a morphism.
 We say that \(f\) is an {\em weak equivalence} if  
\(f:C^*_+ \to D^*_+\) is a homotopy equivalence. 
\end{defn}

\begin{lem}
A weak equivalence is a quasi-isomorphism.
\end{lem}

\begin{proof} Suppose \(E^*\) is an object of \(GMF_{-w}(R)\). Then 
 \(f\otimes 1: C^*\otimes E^* \to
  D^* \otimes E^*\) is a morphism of \(GMF_0(R)\), so 
the induced map \((f \otimes 1)^*\) commutes with
  \(d_-^*\). On the other hand, the map \((f\otimes 1): C^*_+
  \otimes E^*_+ \to D^* \otimes E^*_+\) is a homotopy equivalence, so
  \((f \otimes 1)^*\) is an isomorphism. 
\end{proof}


A second source of quasi-isomorphisms is provided by a process we 
refer to as {\it twisting}. Suppose that \(C^*\) is a matrix
factorization of length \(3\), so that \(C^i\) is trivial for \(i\neq
0,1,2\). Given a homomorphism \(H:C^2 \to C^0\), we define  
a deformed version of \(d_-\) by the equation
\begin{equation*}
d_-(H) = d_- + d_+H - Hd_+.
\end{equation*}
The
twisted factorization \(C^*(H)\) is the triple \((C^*,d_+,d_-(H))\). 
\begin{lem}
 \(C^*(H)\) is a graded matrix
factorization with the same potential  as 
\(C^*\). 
\end{lem}

\begin{proof}
Suppose that \(C\) has potential \(w\). It is enough to check that
\begin{align*}
d_+d_-(H)+d_-(H)d_+ & = d_+d_- +d_+^2H -d_+Hd_+ + d_-d_+ + d_+Hd_+
-Hd_+^2 \\
& = d_-d_+ + d_+d_- = w
\end{align*}
and \(d_-(H)^2 = 0\). The latter expression contains nine terms. Five
of these (\(-d_-Hd_+\), \(d_+Hd_-\), \(d_+Hd_+H\), 
\(Hd_+Hd_+\) and \(-d_+H^2d_+\))
vanish for dimensional reasons. Two others (\(d_-^2\) and
\(-Hd_+d_+H\)) vanish because \(C^*\) is a matrix factorization. The
remaining two terms  \(d_-d_+H\) and \(-Hd_+d_-\) represent nontrivial
maps \(C^2 \to C^0\). On \(C^2\), \(d_-d_+\) vanishes for
dimensional reasons, so  \(d_+d_- = w \cdot \text{Id}\). Similarly, on
\(C^0\),  \(d_-d_+
= w\cdot \text{Id}\). Thus the final two terms cancel each other.
\end{proof}

\begin{lem}
\label{Lem:Twist}
If \(C^*\) and \(H\) are as above, the obvious
identification \(C^*_+ \cong C^*(H)_+\) is a quasi-isomorphism. 
\end{lem}

\begin{proof}
Suppose \(E^*\) is a matrix factorization with potential
\(-w\). Viewed as endomorphisms of the complex \(C^*_+\otimes E^*_+\), the
negative differentials on \(C^*\otimes E^*\) and \(C^*(H) \otimes
E^*\) have the form \(d_-\otimes 1 \pm 1 \otimes d_-\) and 
\(d_-(H)\otimes 1 \pm 1 \otimes d_-\). Their difference 
\((d_-(H)-d_-)\otimes 1= (d_+H-Hd_+) \otimes 1\) is null-homotopic. 
\end{proof}

\subsection{Koszul factorizations} Suppose \(R\) is a ring and that
\(a,b \in R\). The {\it short matrix factorization} \(\{a,b\}\)
is the rank two factorization given by the diagram
$$\xymatrixcolsep{4pc}  
\xymatrix{ R  \ar@<0.4ex>[r]^{b}&
  R \ar@<0.4ex>[l]^{a}}.$$
It has potential \(ab\). 

\begin{defn} \cite{KRI}
 Suppose \({\bf a} = (a_1,\ldots,a_n)\) and \( {\bf b} = 
(b_1,\ldots,b_n)\) are elements of \(R^n\).
The  {\it Koszul factorization} \(\{{\bf a},{\bf b}\}\) is the tensor
product of the short factorizations \(\{a_i,b_i\}\):
\begin{equation*}
\{{\bf a},{\bf b}\} = \bigotimes_{i=1}^n \{a_i,b_i\}.
\end{equation*}
It is a \(\Z\)--graded matrix 
factorization over \(R\), with potential \({\bf a \cdot b} = \sum_i
a_ib_i\). We say that the {\it order} of the factorization is \(n\). 
\end{defn}

\noindent
When we want to explicitly record the values of \(a_i\) and \(b_i\),
we represent \(\{{\bf a},{\bf b}\}\) by the {\it Koszul matrix}
\begin{equation*}
\begin{pmatrix}
a_1 & b_1 \\
a_2 & b_2 \\
\ldots & \ldots \\
a_n & b_n 
\end{pmatrix}
\end{equation*}

More intrinsically, we can view the underlying module of 
\(\{{\bf a},{\bf b}\}\) as the exterior algebra \(\Lambda^*R^n\),
where \({\bf b}\) as an element of \(R^n\)
and \({\bf a}\) is an element of the dual module \((R^n)^*\).  
The differentials are given by
\begin{equation*}
d_+({\bf x}) = {\bf x} \wedge {\bf b} \quad \text{and} \quad 
d_-({\bf x}) = {\bf x} \ts \neg \ts {\bf a}.
\end{equation*}
From this perspective, it's clear that if we express \({\bf b}\) and
\({\bf a}\) in terms of a new basis for \(R^n\) and its dual basis,
the resulting Koszul factorization will be isomorphic to 
\(\{{\bf a},{\bf b}\}\). In particular, consider the change of
basis operation which replaces the standard basis element 
 \({\bf e}_i\) of \(R^n\) with \({\bf e}_i + c {\bf
  e}_j\). At the level of Koszul matrices, this corresponds to the
{\it row operation} which sends
\begin{equation*}
\begin{pmatrix}
a_i & b_i \\  a_j & b_j
\end{pmatrix}
\mapsto
\begin{pmatrix}
a_i + c a_j & b_i \\  a_j & b_j - cb_i
\end{pmatrix}
\end{equation*}
and leaves the remaining rows of the Koszul matrix unchanged.


We now recall an important technical tool introduced in
  \cite{KRI}. This is the process of
``excluding a variable.'' Suppose that \(C = \{{\bf
  a},{\bf b}\}\) is a Koszul factorization over the ring \(R[x]\)
with potential \(w\) which happends to be contained in \(R\), 
and that \(b_1=f(x)\) is a
monic polynomial of positive degree in \(x\). Let
\({\bf a}',{\bf b}' \in R[x]^{n-1}\) be the vectors obtained from 
\({\bf a}\) and \({\bf b}\) by omitting the first component, and put
\(C' = \{{\bf a}',{\bf b}'\}.\)   At the level of modules,
\(C \cong C' \oplus C'\), and with respect to this decomposition the
 differentials on \(C\) are given by 
\begin{align*}
d_+(u,v) & = (d_+'u,d_+'v + f(x)u) \\
\quad d_-(u,v) & = (d_-'u+a_1 v,d_-'v).
\end{align*}

Next, we form the quotient ring \(R_1
= R[x]/(f(x))\), and let \(\pi: R[x] \to R_1\) be the projection. 
The factorization \(C'' =  \{\pi({\bf a}'),\pi({\bf b}')\}\) 
is a Koszul factorization over \(R_1\) with potential \(\pi(w) = w \in R\). 

\begin{lem}
\label{Lem:Exclude}
 The map \(\phi:C \to C'' \) defined by
\(\phi(u,v) = \pi(v)\) is a weak equivalence in the category \(GMF_w(R)\). 
\end{lem} 

\begin{proof}
Using the formulas above, it is easy to see that \(\phi\) defines a
morphism of matrix factorizations. Thus we need only verify that
\(\phi\) has a homotopy inverse with respect to \(d_+\). Since \(f\)
is monic, every \(r \in R_1\) may be written uniquely in the form
\(r = r_0 + r_1 x+ \ldots r_{k-1}x^{k-1}\), where \(r_i \in R\) and 
 \(k=\deg f\). The map
which sends \(r\in R_1\) to this representative defines an 
\(R\)--module  homomorphism  \(\iota: R_1 \to R[x]\).  \(C'\) and
\(C''\) are free  over  \(R[x]\) and \(R_1\),
respectively,  so \(\iota\) can be used to define an \(R\)--module
homomorphism \(\iota: C'' \to C'\). 
We define a map  
\(\psi: C'' \to C\) by  
\begin{equation*}
\psi(y) = \left( \frac{\iota (d_+''y) - 
d_+'\iota(y)}{f(x)},\iota(y)\right).
\end{equation*}
It is easy to see that
\(\psi\) commutes with \(d_+\) and that  \( \phi \psi = Id_{C''}\). Finally,
we define \(H:C \to C\) by 
\begin{equation*}
H(u,v) = \left( \frac{v - \iota(\pi(v))}{f(x)}, 0 \right ). 
\end{equation*}
We leave it as an exercise to the reader to check that \(H\) is a
homotopy between \(\psi \phi \) and \(Id_{C}\). 
\end{proof}

\begin{cor}
\label{Cor:QIRed}
Suppose \(C^* \sim D^*\) as objects of
 \(GMF_w(R[X])\), where \(w \in R\). 
 Then the quotients \(C^*/(X)\) and \( D^*/(X)\) are
 quasi-isomorphic as objects of  \(GMF_w(R)\). 
\end{cor}

\begin{proof}
By the lemma, the quotient \(C^*/(X)\) is quasi-isomorphic to 
\(C^* \otimes \{0,x\}\). This, in turn, is
quasi-isomorphic to  \(D^* \otimes  \{0,x\}\). 
\end{proof}
 

Now suppose that \(R\) is a polynomial ring, that \(w \in R\), and
that \({\bf b} \in R^n\). It is clear that we can choose \({\bf a} \in
(R^n)^*\) so that \(\{{\bf a},{\bf b}\}\) has
potential \(w\) if and only if \(w\) is in the ideal generated by the
\(b_i\). To what extent is the choice of  \({\bf a}\) unique? 
When \(n=1\), we have \(a_1b_1=w\), so \(a_1\) is uniquely determined
unless \(b_1=w = 0\). For \(n=2\), we have the following result.

\begin{lem}
\label{Lem:KTwist}
Suppose \(R\) is a UFD and that \(\{{\bf a},{\bf b}\}\) and 
 \(\{{\bf a}',{\bf b}'\}\) are two order two Koszul factorizations
 over \(R\) with potential \(w\). If \(b_1\) and \(b_2\) are
 relatively prime, the factorizations 
 \(\{{\bf a},{\bf b}\}\) and \(\{{\bf a}',{\bf b}'\}\) are related by
 a twist. 
\end{lem}

\begin{proof}
We have
\begin{equation*}
a_1b_1 + a_2b_2 = w = a_1'b_1+a_2'b_2,
\end{equation*}
which implies that  \((a_1-a_1')b_1+ (a_2-a_2')b_2 = 0.\) Since \(b_1\) 
and \(b_2\) are relatively prime, our two factorizations must be
represented by Koszul matrices of the form
\begin{equation*}
\begin{pmatrix}
a_1&b_1 \\ a_2& b_2
\end{pmatrix}
\quad \text{and} \quad 
\begin{pmatrix}
a_1-kb_2 &b_1 \\ a_2 +kb_1& b_2
\end{pmatrix}.
\end{equation*}
The second factorization is a twist of the first one, via the map
\(H:R \to R\) which sends \(x \mapsto kx\). 
\end{proof}
\noindent {\bf Remark:} In fact, it is not difficult to see that 
\(\{{\bf a},{\bf b}\}\) and \(\{{\bf a}',{\bf b}'\}\) are
 isomorphic as \(\Z/2\)--graded matrix factorizations. 

\subsection{Equivalence of definitions}
\label{SubSec:DefCom}

The ideas described above can be used to verify that the definitions of the
various KR--homologies given in section~\ref{Sec:Defs} agree with
those in \cite{KRI} and \cite{KRII}. 
We assume the
reader is already somewhat familiar with these papers, and only 
briefly recall their content.

To a planar diagram \(D\), we  associate the  ring
\(R'(D) = \Q[X_i]\), where \(i\) runs over the edges of \(D\). 
In \cite{KRI}, Khovanov and Rozansky assign to  \(D\) a complex of matrix
factorizations \(C_p'(D)\) defined over  \(R'(D)\) and with
potential \(w_p(D)\). 
(Although the definition in \cite{KRI} is only
stated for \(p(x)=x^{N+1}\), it works equally well for any \(p\), as
 implicitly noted by Gornik \cite{Gornik}.) 
\(C_p'(D)\) is a tensor product of factors, one for each internal
vertex of \(D\). These factors are as follows. To a mark
 with incoming and outgoing edges labeled \(i\) and \(j\), Khovanov and
Rozansky associate the short factorization
\begin{equation*}
\left\{\frac{p(X_j)-p(X_i)}{X_j-X_i}, X_j-X_i\right\}.
\end{equation*}
 To the singular diagram
\(D_s\), they associate an order \(2\) Koszul factorization given by
the Koszul matrix 
\begin{equation*}
C_p'(D_s) = \begin{pmatrix}
* & X_k+X_l-X_i-X_j \\ * & X_kX_l-X_iX_j
\end{pmatrix}.
\end{equation*}
According to the remark following  Lemma~\ref{Lem:KTwist}, any two such
factorizations are isomorphic as \(\Z/2\)--graded factorizations.
Thus the  entries in the left-hand column are more
or less immaterial, and we will simply mark them by \(*\)'s. 

Finally, the positive and negative crossings are associated to short
complexes of order 2 Koszul factorizations, as follows: 
\begin{align*}
C_p'(D_+) & = 
\begin{pmatrix}
* & X_k+X_l-X_i-X_j \\
* & X_kX_l-X_iX_j 
\end{pmatrix}
\xrightarrow{\phantom{X} \chi_1 \phantom{X}} 
\begin{pmatrix}
* & X_l-X_j \\
* & X_k-X_i 
\end{pmatrix} \\
C_p'(D_-) & = 
\begin{pmatrix}
* & X_l-X_j \\
* & X_k-X_i 
\end{pmatrix}
\xrightarrow{\phantom{X} \chi_0 \phantom{X}} 
\begin{pmatrix}
* & X_k+X_l-X_i-X_j \\
* & X_kX_l-X_iX_j 
\end{pmatrix}.
\end{align*}
The composition \(\chi_0\chi_1\) is given by multiplication by
\(X_k-X_j\). 
Applying  a row operation, we see that these complexes are isomorphic to
\begin{align*}
C_p'(D_+) & \cong 
\begin{pmatrix}
* & X_k+X_l-X_i-X_j \\
* & X_kX_l-X_iX_j 
\end{pmatrix}
\xrightarrow{\phantom{X} \chi_1 \phantom{X}} 
\begin{pmatrix}
* & X_k+X_l-X_i-X_j \\
* & X_k-X_i 
\end{pmatrix} \\
C_p'(D_-) & \cong 
\begin{pmatrix}
* & X_k+X_l -X_i-X_j \\
* & X_k-X_i 
\end{pmatrix}
\xrightarrow{\phantom{X} \chi_0 \phantom{X}} 
\begin{pmatrix}
* & X_k+X_l-X_i-X_j \\
* & X_kX_l-X_iX_j 
\end{pmatrix}.
\end{align*}

The matrix factorizations used by Khovanov and Rozansky are 
 \(\Z/2\)--graded, 
 rather than the \(\Z\)--graded factorizations that we
have been considering. One advantage of this approach is that there is
a good notion of homotopy equivalence for such factorizations. This
enables them to work in the homotopy category \(hmf_w(R)\) of
\(\Z/2\)--graded matrix factorizations with potential \(w\). There is
an obvious forgetful functor from \(Kom(GMF_w(R))\) to \(Kom(hmf_w(R))\), so we
can view both \(C_p'(D)\) and \(\Cu_p(D)\)  as  objects
 of the latter category.

\begin{lem}
\label{Lem:Exclusion}
 If \(D\) is a closed diagram, 
\(C_p'(D) \cong \widetilde{C}_p(D)\) in 
\(Kom(hmf_0(R(D)))\). 
\end{lem}

\begin{proof}
In both cases, the complex associated to a disconnected diagram is the
tensor product of the complexes associated to its components. Thus 
we may  assume that \(D\) is connected. We fix an edge \(i\) of \(D\) and
consider the diagram \(D(i)\) obtained by inserting a
bivalent vertex \(v_0\) into \(i\). 
For each vertex \(v\) of \(D(i)\), the linear relation \(\rho(v)\) 
appears as a
matrix entry of every Koszul factorization in the complex
\(C_p'(D(i))\). By 
Lemma~\ref{Lem:EdgeRing}, the relations \(\{\rho(v) \ts \vert \ts v
\neq v_0\}\) are all linearly
independent. Thus we can apply Proposition 10 of \cite{KRI} to
exclude them. The result is an isomorphic complex \(C_1\)
defined over the ring \(R'(D_0)/(\rho(v)) \cong R(D)\).

It is shown in \cite{KRI} that 
\(C_p'(D) \cong C_p'(D(i))\), so to prove the lemma, it is enough to
show that 
 \(C_1 \cong \widetilde{C}_p(D,i)\). To see
this, we examine each factor in the complex individually. For example,
consider the  factor associated to a singular crossing. 
 \(X_kX_l-X_iX_j = -(X_k-X_i)(X_k-X_j)\) in \(R(D)\), so \(C_p'(D_s)\)
 reduces to a short factorization of the form \(\{\beta,
 -(X_k-X_i)(X_k-X_j)\}\), where \(\beta\) is the image of some \(\beta '\in
 R'(D)\) which satisfies 
\begin{equation*}
\alpha' (X_k+X_l-X_i-X_j) + \beta' (X_kX_l-X_iX_j) = w_p(D_s) =W_p(X_i,X_j,X_k,X_l).
\end{equation*}
It follows that \(\beta = -W_p/(X_k-X_i)(X_k-X_j)=p_{ij}\) in
\(R(D)\). This is the factorization assigned to \(D_s\) in
section~\ref{SubSec:ETangle}. 

A similar argument shows that  \(C'_p(D_+)\) and
\(C'_p(D_-)\) reduce to complexes of the form
\begin{align*}
 \{p_{ij}, -(X_k-X_i)(X_k-X_j)\} &\xrightarrow{\phantom{X}
    \chi_1 \phantom{X}} \{p_i, X_k-X_i\} \\
 \{p_i, X_k-X_i\} &\xrightarrow{\phantom{X}
    \chi_0 \phantom{X}}\{p_{ij}, -(X_k-X_i)(X_k-X_j)\}
\end{align*}
Since the composition \(\chi_0\chi_1\) is given by multiplication by
\(X_k-X_j\), it is not difficult to see that \(\chi_0\) and \(\chi_1\)
agree with the corresponding maps defined in section~\ref{SubSec:ETangle}.
 
Finally, we consider the short factorization
\begin{equation*}
\left \{\frac{p(X_{i_+}) - p(X_{i_-})}{X_{i_+}-X_{i_-}},
X_{i_+}-X_{i_-}\right \}
\end{equation*} 
coming from the vertex \(v_0\). 
  \(X_{i_+} = X_{i_-}\) in \(R(D)\), so this 
 reduces to the short factorization \(\{p'(X_i), 0\}\) which appears in the
 definition of \(\widetilde{C}_p(D,i)\). 
\end{proof}

\begin{prop}
\(\H_N(L,i)\) and \( \Hu_N(L)\) are isomorphic (as doubly graded
groups) to the reduced and unreduced \(sl(N)\) homology of
\cite{KRI}. 
\end{prop}

\begin{proof}
Suppose \(L\) is represented by a planar diagram  \(D\).
The unreduced  \(sl(N)\) homology of \(L\) is defined to be
\(H(H(C_p'(D),d_{tot}),d_v^*)\), where \(p(x)=x^{N+1}\). From the
lemma, it follows that this is isomorphic to the group \(
H(H(\Cu_p(D,i),d_{tot}),d_v^*)\) which appears in
Definition~\ref{Def:Unreduced}. 

The argument for reduced homology is slightly more involved. In
 \cite{KRI}, the 
reduced homology of \(L\) with respect to an edge \(i\) is defined to be
\(H(H(C_p'(D),d_{tot})/X_i,d_v^*) \). Comparing with
 Definition~\ref{Def:Reduced}, we see that we must show that 
\begin{equation*}
H(C_p'(D),d_{tot})/X_i \cong H(C_p(D)/X_i,d_{tot}). 
\end{equation*}
The complex \(C_p'(D)\) is free  over \(\Q[X_i]\),
but it is shown in \cite{KRI} that
 \(H(C_p'(D),d_{tot})\) is a torsion module over
\(\Q[X_i]\). Applying the universal coefficient theorem, we see that
\begin{equation*}
H(C_p'(D)/X_i, d_{tot}) 
\cong H(C_p'(D),d_{tot})/X_i \otimes H^*(S^1).
\end{equation*}
On the other hand, the lemma tells us that the quotient
\(C_p'(D)/X_i\) is homotopy equivalent to
\(\widetilde{C}_p(D,i)/X_i \cong C_p(D)/X_i \otimes
U_p(i)/X_i\). \(U_p(i)/X_i\)
is a rank two factorization with trivial differentials, so 
\begin{equation*}
H(\widetilde{C}_p'(D)/X_i,d_{tot}) \cong H(C_p(D)/X_i,d_{tot})\otimes
H^*(S^1). 
\end{equation*}
Canceling out the extra factors of \(H^*(S^1)\), we obtain the
desired isomorphism

It remains to check that the bigradings 
 agree. For the second, homological
grading, this is clearly the case --- it is given by
  \(\gr_v\) in both cases.
To see that  \(\gr_N = i + (N-1)j/2\) coincides with the \(q\)--grading of
\cite{KRI}, first note that the complex \(C_p'(D)\) is set up so that
the right-hand group in each linear factor is unshifted with respect to
 the \(q\)--grading. If we exclude the linear term appearing in such a
 factor, the \(q\)--grading is unaffected.
Thus it suffices to check
that the gradings agree on quadratic factors. Consider the factorization
\(C_p(D_s)\) associated to a singular point. According to
section~\ref{SubSec:ETangle}, the two copies of \(R(D)\) used to define
this factorization have \((i,j)\) grading shifts of \(\{1,-2\}\) and
\(\{-1,0\}\). These correspond to shifts of \(\{2-n\}\) and \(\{-1\}\) in
\(\gr_N\), which precisely 
match the shifts in the \(q\)--grading which appear in the
definition of \(C_p'(D_s)\) on p. 48 of \cite{KRI}. The calculation
for \(C_p'(D_\pm)\) is similar, except it also uses the
grading shifts on p. 81 of \cite{KRI}, part of which goes into 
the shifts in \(C_p(D_{\pm})\), and part into the overall shift by
\(w(N-1)\) which appears in the definition of \(\Hu_N(L)\) and \(\H_N(L,i)\). 
\end{proof}

\begin{prop}
 The middle HOMFLY homology  \(H(L)\) is isomorphic to the HOMFLY
 homology of \cite{KRII}. The identification is such that an element
 with grading \((i,j,k)\) in our notation corresponds to an element
 with grading \((j/2,i-j/2,k/2)\) in the notation of \cite{KRII}.
\end{prop}

\begin{proof}
The homology of  \cite{KRII} is defined to be
\(H(H(C_a'(D),d_{tot}),d_v^*)\), where \(C_a'(D)\) is a certain complex of
matrix factorizations defined over the ring \(R'(D)[a]\). If we
substitute \(a=0\), \(C_a'(D)\) reduces to \(C_0'(D)\). On the other
hand, it is proved in \cite{KRII} that \(a\) acts by \(0\) on
\(H(C_a(D),d_{tot})\). Applying the universal coefficient theorem, we
find that 
\begin{equation*}
H(C_0'(D), d_{tot}) \cong H(C_a'(D),d_{tot}) \otimes
H^*(S^1).
\end{equation*}
On the other hand, the lemma implies that
\begin{equation*}
H(C_0'(D), d_{tot})  \cong H(\Cu_0(D), d_{tot}) 
 \cong H(C_0(D)\otimes U_0(i),d_{tot}).
\end{equation*}
Since
\(U_0(i)\) has trivial differential, the last group is isomorphic to 
\( H(C_0(D),d_{tot})\otimes H^*(S^1)\).
Cancelling the factors of \( H^*(S^1)\), we see that 
\begin{equation*}
H(H(C_a'(D),d_{tot}),d_v^*) \cong H(H(C_0(D),d_{tot}),d_v^*) \cong H(L).
\end{equation*}

It remains to compare the triple grading on the two theories.  The
 ring \(R'(D)[a]\) is bigraded, with an additional grading
 corresponding to the power of \(a\) as well as the usual \(q\)--grading.
 The first grading in \cite{KRII} is nominally given by the power of
 \(a\). Since \(a\) acts by \(0\) on homology, however, 
  any class is homologous to one represented by elements of \(R'(D)\). The
 \(a\)--grading  of such a class
comes entirely from the grading shifts introduced in the
 definition of \(C_a'(D)\). It is easily verified that these 
 shifts are the same as those for \(\gr_h\), so the first grading is
 \(\gr_h = j/2\). The second grading in \cite{KRII} corresponds to the usual
 \(q\)--grading on the ring \(R(D)\), but the grading shifts in 
 \(C_a'(D)\) differ from ours. Up to an overall shift, the grading shift in 
 \cite{KRII} corresponds to the difference between our shift in \(q\)
 and \(gr_h\). Thus the second grading is given by \(i-j/2\). 
Finally, 
 the third grading in \cite{KRII} is given by \(\gr_v = k/2\).
\end{proof}

As an further application of these techniques, we can now make good
on our claim  that the unreduced complex is
independent of the choice of the marked edge used to define it.

\begin{prop}
\label{Prop:MarkQI}
If \(i\) and \(j\) are two edges of a connected diagram
 \(D\), the unreduced complexes \(\Cu_p(D,i)\) and \(\Cu_p(D,j)\) are
 quasi-isomorphic. 
\end{prop}

\begin{proof}
Let \(D(i,j)\) be the diagram obtained by inserting bivalent vertices
\(v_i\) and \(v_j\) 
in edges \(i\) and \(j\). Consider the complex \(C_p'(D(i,j))\) as an
element of the category \(GMF_0(R(D(i,j))\). 
Arguing as in the proof of Lemma~\ref{Lem:Exclusion}, we use
Lemma~\ref{Lem:Exclude} to exclude the linear relations \(\{\rho(v)
\ts \vert \ts v \neq v_i\}\). The result is a new complex of matrix
factorizations \(C_i\) which is quasi-isomorphic to \(C_p'(D(i,j))\). 
The  same argument used in the proof of Lemma~\ref{Lem:Exclusion} shows
that \(C_i\cong \Cu_p(D,i)\). Thus \(C_p'(D(i,j))\) is
quasi-isomorphic to \(\Cu_p(D,i)\). Similarly, \(C_p'(D(i,j))\)  is
quasi-isomorphic to \(\Cu_p(D,j)\). This proves the claim.
\end{proof}

\section{Braid Graphs and MOY Relations}
\label{Sec:Singular} 

A tangle diagram all of whose crossings are
singular is called a {\it graph}; a braid diagram all of whose
crossings are singular is  a {\it braid graph.}  In \cite{MuOhYa}, 
Murakami, Ohtsuhki, and Yamada  explain how to assign
 a HOMFLY polynomial \(\wt{P}(D)\) to a closed graph \(D\). This assignment
 can be used to give a state model definition of the HOMFLY
 polynomial similar to the Kauffman state model \cite{Kauffman2} for the Jones
 polynomial.  (See \cite{JonesHOM,Turaev} for  related constructions.)
Murakami, Ohtsuki, and Yamada also show that the HOMFLY polynomial of
a graph satisfies certain relations, which we refer to as {\it MOY
  relations}. 

In this section, we briefly review these results and describe their
generalizations to KR-homology. 
In \cite{KRI,KRII}, Khovanov and Rozansky show that \(C_p(D)\) 
 satisfies relations analogous to the MOY relations for the
HOMFLY polynomial. The main technical result of this section is 
that these relations continue to hold in the context of
\(\Z\)--graded matrix factorizations. As an application, we show that
the HOMFLY homology of a braid graph is determined by its HOMFLY
polynomial. We then use the MOY state model to give a proof of 
  Theorem~\ref{Thm:Chi}  along the lines of the 
proof for the \(sl(N)\) homology given  in \cite{KRI}. 

\subsection{The MOY state model}
We begin by recalling the state model of Murakami, Ohtsuki,
and Yamada  \cite{MuOhYa}. Although their paper is phrased in terms
of the \(sl(N)\) polynomials, the results we want are easily translated
into the language of the HOMFLY polynomial, and we will state them in
this form.

\begin{figure}
\includegraphics{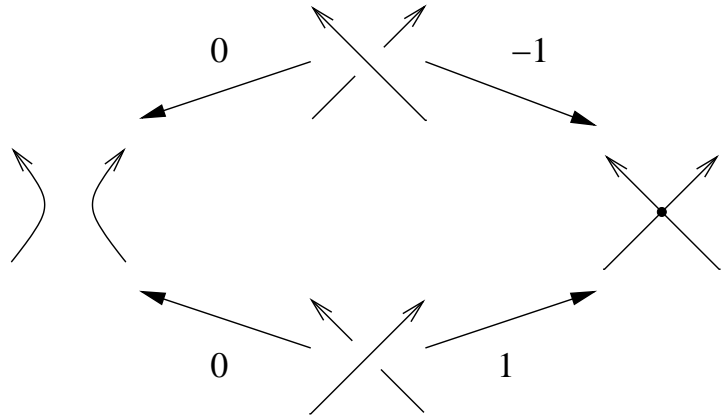}
\caption{\label{Fig:Resolutions} Resolutions and their weights.}
\end{figure}

Suppose \(D\) is a  diagram representing an oriented link
\(L\). We can ``resolve'' each crossing of \(D\) in one of two ways:
either into a pair of arcs (the oriented resolution) or into a
singular crossing. To each such resolution, we assign a weight \(\mu
\in \Z\), depending on whether the crossing is positive or
negative and on which resolution it receives. The possible
resolutions and their weights are illustrated in
Figure~\ref{Fig:Resolutions}.

A {\it state} of the diagram \(D\) is a choice of resolution for each
crossing of \(D\). If \(D\) has \(n\) crossings, it will have  \(2^n\)
different states. To a state \(\sigma\), we assign a weight
\(\mu(\sigma)\) given by the sum of the local weights at each
crossing. In addition, each  \(\sigma\) gives rise to a graph
\(D_\sigma\). In \cite{MuOhYa},
it is shown that the
unnormalized HOMFLY polynomial of \(L\) is given by the formula
\begin{equation}
\label{Eq:MOY}
\wt{P}(L) = (aq^{-1})^{w(D)} \sum_{\sigma} (-q)^{\mu(\sigma)}
\wt{P}(D_\sigma)
\end{equation}
where the  quantity  \(\wt{P}(D_\sigma)\) is an invariant of the 
 graph \(D_\sigma\). We think of
 \(\wt{P}(D_\sigma)\) as the HOMFLY polynomial of \(D_{\sigma}\), and
 view formula~\eqref{Eq:MOY} as generalizing the definition of \(\wt{P}\)
  to  closed tangle diagrams with an arbitrary number of singular
 crossings.

\subsection{Polynomials of braid graphs} 
In order to use formula~\eqref{Eq:MOY}, we need some way to determine
\(\wt{P}(D)\) when \(D\) is a graph.
In \cite{MuOhYa}, the authors give 
a direct geometric procedure for finding these polynomials --- or 
  rather, their specializations to \(a=q^N\). For our purposes,
  however, it is more convenient to characterize \( \wt{P}(D)\) in
  terms of certain relations given in \cite{MuOhYa}. 

\begin{figure}
\input{MOY.pstex_t}
\caption{\label{Fig:MOYMoves} MOY relations \(O,I\), and \(II\).} 
\end{figure}

\begin{figure}
\input{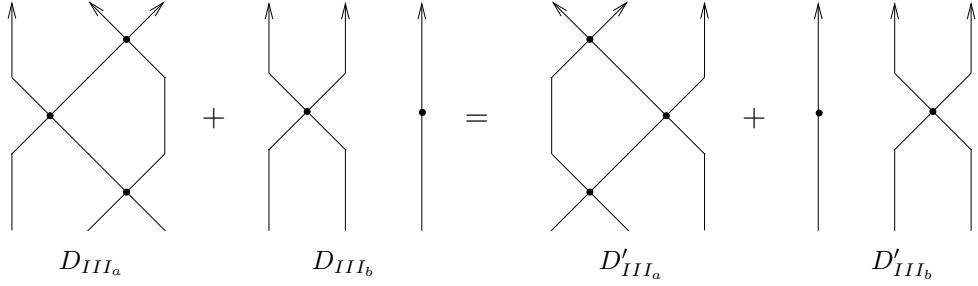}
\caption{\label{Fig:MOY3} MOY relation \(III\).} 
\end{figure}

Suppose \(D_{O}\), \(D_{I}\), \(D_{II}\), \(D_{IIIa}\) and \(D_{IIIb}\)
 are braid graphs containing 
regions like  those shown on the left-hand sides of
Figures~\ref{Fig:MOYMoves} and \ref{Fig:MOY3}, and let 
 \(D_{O'}\), \(D_{I}'\), \(D_{II}'\), \(D_{IIIa}'\)and \(D_{IIIb}'\)
be the graphs obtained by replacing this region with the
corresponding one on the right-hand side of the figure. It is
 shown in \cite{MuOhYa} that 
\(\wt{P}\) satisfies the following {\it MOY relations}:
\begin{align*}
O) \quad &    \wt{P}(D_{O})  = \frac{a-a^{-1}}{q-q^{-1}} \wt{P}(D_{O}') \\
 I) \quad &    
 \wt{P}(D_{I})  = \frac{aq^{-1}-a^{-1}q}{q-q^{-1}} \wt{P}(D_{I}') \\
II) \quad &  \wt{P}(D_{II})  = (q+q^{-1}) \wt{P}(D_{II}')  \\
 III) \quad &  \wt{P}(D_{IIIa})  + \wt{P}(D_{IIIb})  =
 \wt{P}(D_{IIIa}') + \wt{P}(D_{IIIb}').
\end{align*}

The HOMFLY polynomial of a braid graph is completely determined by
 these relations. To see this, we use an induction scheme
 introduced by Wu \cite{Wu}. Suppose \(D\) is a
braid graph on \(b\) strands. The crossings of \(D\) are naturally
arranged into \(b-1\) columns, which we number \(1,\ldots b-1\) going
from left to right. If \(c\) is a crossing of \(D\), let \(i(c)\)
be the number of the column containing it. Following Wu, we define the
{\it complexity} of \(D\) to be the sum

\begin{equation*}
i(D) = b + \sum_{c} i(c). 
\end{equation*}

The complexity of a diagram on the
left-hand side of Figure~\ref{Fig:MOYMoves} is strictly greater than the
complexity of the corresponding diagram on the right. Similarly, the
complexity of diagram \(D_{III_a}\) is greater than that of the other
three diagrams in Figure~\ref{Fig:MOY3}. 

\begin{lem}
\label{Lem:Hao}
\cite{Wu} Suppose \(D\) is a nonempty braid graph which is the
closure of an open braid graph \(D_o\). Then either \(D\) contains a
region of the form \(D_O\) or \(D_I\) or \(D_o\) contains a region of
the form \(D_{II}\) or \(D_{III_a}\).
\end{lem}

\noindent In other words, \(D\) can be 
related to braid graphs of lesser complexity by one of MOY moves
\(O-III\). Moreover, we may assume that moves of type \(II\) and
\(III\) take place in the open braid \(D_o\). 

\begin{cor}
If \(D\) is braid graph, \(\wt{P}(D)\) is determined by MOY
relations \(O-III\) and the fact that \(\wt{P}\) of the empty graph
is \(1\). 
\end{cor}

\subsection{Homology of braid graphs}
\label{SubSec:SingHom}
The KR-complex of a braid graph satisfies
 decomposition rules analogous to the MOY relations
\(O\)--\(III\). 
In the context of
 \(\Z/2\)--graded matrix factorizations, such rules were  introduced
 in \cite{KRI} and later applied to the HOMFLY
 homology in \cite{KRII, Wu}. 
Similar {MOY decompositions}  also hold  in the derived 
category of \(\Z\)--graded matrix factorizations.
We collect their statements here, but postpone the proofs to the
 end of this section. Note that although
\(\Cu_p(D)\) is generally a complex of matrix factorizations, when
\(D\) is a graph, the complex is supported in a single vertical
grading. Thus \(\Cu_P(D)\) is most naturally viewed as an object of
the \(GMF(R(D))\). It is doubly graded, with gradings  \((q,2\gr_h)\). 


\begin{prop}
\label{Prop:MOYO}
Let \(D_O\) and \(D_O'\) be two  braid graphs related as in the first
line of 
Figure~\ref{Fig:MOYMoves}, and let \(C_p(O)\) be the matrix
factorization
$$\xymatrixcolsep{5pc} 
\xymatrix{ \Q[X_1]\{0,-2\}  \ar@<0.4ex>[r]^{0}&
  \Q[X_1]\{0,0\}. \ar@<0.4ex>[l]^{p'(X_1)}}$$
Then \(C_p(D_O) \cong C_p(D_O') \otimes_{\Q} C_p(O) \) in 
\(GMF(R(D_O))\). 
\end{prop}

\begin{prop}
\label{Prop:MOYI}
Let \(D_I\) and \(D_I'\) be two braid graphs related as in
 the second line of 
 Figure~\ref{Fig:MOYMoves}, and let \(C_p(I)\) be the matrix
factorization
$$\xymatrixcolsep{5pc} 
\xymatrix{ \Q[X_1,X_2]\{1,-2\}  \ar@<0.4ex>[r]^{0}&
  \Q[X_1,X_2]\{-1,0\} \ar@<0.4ex>[l]^{p'_{12}}}$$
where \(p'_{12} = (p'(X_1) -
    p'(X_2))/(X_1-X_2)\). 
Then \(C_p(D_I) \cong C_p(D_I') \otimes_{\Q[X_1]} C_p(I) \) in 
\(GMF(R(D_I))\). 
\end{prop}


\begin{prop} 
\label{Prop:MOYII}
Let \(D_{II}\) and \(D_{II}'\) be two braid graphs formed by taking
the union of a fixed  graph \(D\) with the diagrams in the last line of
Figure~\ref{Fig:MOYMoves}. Then 
\begin{equation*}
 C_p(D_{II}) \sim C_p(D_{II}')\{-1,0\} \oplus
C_p(D_{II}')\{1,0\}. 
\end{equation*}
in  \(GMF(R(D))\). 
\end{prop}

\begin{prop} 
\label{Prop:MOYIII}
Let \(D_{IIIa}\), \(D_{IIIb}\),  \(D_{IIIa}'\), and 
\(D_{IIIb}'\) be braid graphs 
formed by taking
the union of a fixed  graph \(D\) with the diagrams in
Figure~\ref{Fig:MOY3}. 
Then 
\begin{align*}
 C_p(D_{IIIa}) \oplus C_p(D_{IIIb}) \sim C_p(D_{IIIa}') \oplus C_p(D_{IIIb}')
\end{align*}
in the category \(GMF(R(D))\). 
\end{prop}

As an  immediate consequence, we have relations
\begin{align*}
O)  \quad & \Hu(D_O) \cong (\wt{H}(D_O')\{0,-2\} \oplus
\wt{H}(D_O')) \otimes _\Q \Q[x] \\
I)  \quad & \wt{H}(D_I) \cong (\wt{H}(D_I')\{1,-2\} \oplus
\wt{H}(D_I')\{-1,0\}) \otimes _\Q \Q[x] \\
II) \quad & \wt{H}(D_{II}) \cong \wt{H}(D_{II}')\{-1,0\} \
\oplus \wt{H}(D_{II}')\{1,0\} \\
III) \quad & \wt{H}(D_{IIIa}) \oplus \wt{H}(D_{IIIb}) \cong 
\wt{H}(D_{IIIa}') \oplus \wt{H}(D_{IIIb}')
\end{align*}
which closely parallel the MOY relations for the HOMFLY polynomial. In
fact, these relations are proved by Khovanov and Rozansky in
\cite{KRII}, where they are used to show that \(\wt{H}\) is invariant
under braidlike Reidemeister moves. 

Like the HOMFLY polynomial, the
 HOMFLY homology of a braid graph is determined by the MOY
relations. In fact, the two carry precisely  the same
 information. More specifically, let 
\begin{equation*}
\wt{\mathcal{P}}(D) = \sum_{i,j}
(-1)^{j/2}a^jq^i \dim \Hu^{i,j}(D) 
\end{equation*}
 be the signed
Poincar{\'e} polynomial of \( \wt{H}(D)\). Then we have 

\begin{prop}
\label{Prop:a=grh}
If \(D\) is a closed braid graph on \(b\) strands, 
\(\wt{ \mathcal{P}}(D) = (-aq)^{-b}  \wt{P}(D)\). 
\end{prop}

\begin{proof}  We induct on the complexity of
  \(D\). The base case is the empty diagram, which has complexity
  \(0\), HOMFLY polynomial \(1\), and KR-homology \(\Q\) supported in
  bigrading \((0,0)\). For the induction step, we apply
  Lemma~\ref{Lem:Hao} to see that \(D\) is related to diagrams of
  lesser complexity by an MOY move. To complete the proof, we
  need only check that the MOY relations for \(\wt{P}\) are
  consistent with the corresponding MOY decompositions for \(\Hu\). 

For example, consider MOY move \(O\). By the induction hypothesis, we
know that \(\wt{\mathcal{P}}(D_O') = (-aq)^{1-b}\wt{P}(D_O')\). On the
other hand, relation \(O\) above shows that 
\begin{align*}
\wt{\mathcal{P}} (D_{O}) & = 
(1-a^{-2}) \left(\sum_{i=0}^\infty q^i \right)
\wt{\mathcal{P}}(D_O') \\
& = (-aq)^{-1} \left( \frac{a-a^{-1}}{q-q^{-1}} \right)
\wt{\mathcal{P}}(D_O') \\
& = (-aq)^{-b} \left( \frac{a-a^{-1}}{q-q^{-1}} \right) \wt{P}(D_O') \\
& = (-aq)^{-b} \wt{P}(D_O)
\end{align*}
so the claim holds for \(D_O\) as well. 

We leave it to the reader to check the remaining MOY moves.
The argument for move \(I\) is very similar to the one for move \(O\), 
and moves \(II\) and \(III\) are even easier, since
 all diagrams involved have the same number of strands. 
\end{proof}

Using the MOY relations, it is not difficult to see that if \(D\) is 
a braid graph on \(b\) strands, the denominator of \(\wt{P}(D)\) is 
\((q-q^{-1})^b\). This fact is nicely reflected  in the module
structure of \(\Hu(D)\). To see this, write \(D\) as the closure of an
open braid graph \(D_o\), and label the outgoing edges of 
\(D_o\) by \(1,2,\ldots,b\). (In \(D\), these are identified with the
incoming edges of \(D_o\).) The ring \(R_b = \Q[X_1,X_2,\ldots,X_b]\) is a
subring of \(R(D)\). 

\begin{prop}
\label{Prop:Free}
If \(D\) is a closed braid graph on \(b\) strands, 
\(\Hu(D)\) is a free module of finite rank over \(R_b \). 
\end{prop}

\begin{proof}
Again, we induct on the
complexity of \(D\). The base case is when \(D\) is the empty diagram, and
\(\wt{H}(D) \cong \Q \) is free of rank \(1\) over \(\Q\). For the induction
step, we use Lemma~\ref{Lem:Hao} to see that \(D\) can be simplified
either by a MOY \(O\) or \(I\) move, or by a MOY \(II\) or \(III\)
move which takes place in the open braid \(D_o\). We consider each of
these four possibilities separately. 

For move \(O\), it follows from Proposition~\ref{Prop:MOYO} that 
that \(\Hu(D_O)\) is a direct sum of two copies of \(\Hu(D_O')\)
tensored over \(\Q\)  with \(\Q[X_k]\), where the strand to be
eliminated  has label \(k\). 
By the induction hypothesis,
\(\Hu(D_O')\) is free of finite rank over
\(\Q[X_1,\ldots,X_{k-1},X_{k+1},\ldots,X_b]\), so \(C_p(D_O)\)
will be free of finite rank over \( R_b\). 
The argument for move \(I\) is similar. 

For moves of type \(II\) and \(III\), the fact that 
 the move takes place in  \(D_o\) implies that 
\(R_b\) is contained in the ring \(R(D)\) over which the relations  of
Propositions~\ref{Prop:MOYII} and \ref{Prop:MOYIII} hold. 
 Thus these decompositions also hold over \(R_b\).
The result for move \(II\) follows easily from this, 
since \(\Hu(D_{II})\) is a
direct sum of two copies of \(\Hu(D_{II}')\), which is free of finite
rank by the induction hypothesis. 

For move \(III\), 
the induction hypothesis implies that \(\Hu(D_{III_a}')\) and
 \(\Hu(D_{III_b}')\) 
are free. It follows that
 \(\Hu(D_{III_a})\oplus \Hu(D_{III_b})\) 
is free as well, so \(\Hu(D_{III_a})\) is a  projective module
 over the polynomial ring \(R_b\). By the theorem of Quillen and
 Suslin (see {\it e.g.} \cite{Lang}), any such  module is free.
  Finally, \(\Hu(D_{III_a}')\)
 and \(\Hu(D_{III_b}')\)  are of finite
rank, so the same must be true for \( \Hu(D_{III_a})\). 
\end{proof}


%
%

\subsection{States and the KR-complex}
\label{SubSec:GradedChi}
Now that we understand the relation between
 the MOY state model and  \(\Cu_p(D)\) for braid graphs, we
 consider what happens when
 \(D\) is an arbitrary  braid. 

\begin{lem}
\label{Lem:States}
Suppose \(D\) is a closed braid diagram. Then 
\begin{equation*}
\Hu^+(D) \cong 
\bigoplus_\sigma
\Hu (D_\sigma)\{\mu(\sigma),0,-2\mu(\sigma)\}.
\end{equation*}
where the sum runs over MOY states of \(D\). 
\end{lem}

\begin{figure}
\includegraphics{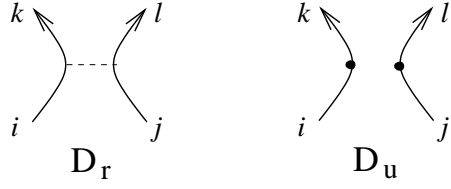}
\caption{\label{Fig:Dr} Diagrams \(D_r\), which represents a four-valent
  vertex, and \(D_u\), which represents a pair of two-valent vertices.}
\end{figure}

\begin{proof}
We temporarily enlarge our notion of a tangle diagram
to include a fourth sort of crossing
\(D_r\), represented by the diagram of Figure~\ref{Fig:Dr}. The
local factor associated to such a crossing is 
$$\xymatrixcolsep{4pc} C_p(D_r) = 
\xymatrix{ R\{0,-2,0\}  \ar@<0.4ex>[r]^{(X_k-X_i)}&
  R\{0,0,0\}. \ar@<0.4ex>[l]^{p_{i}}}$$
The definition of the  KR-complex  is otherwise unchanged. 

Referring to the diagrams in section~\ref{SubSec:ETangle}, we see that 
if we ignore the vertical differential, there are decompositions
\begin{align*}
C_p(D_+) &=C_p(D_s)\{1,0,-2\}  \oplus
  C_p(D_r) \\  C_p(D_-) & =C_p(D_r) \oplus 
  C_p(D_s)\{-1,0,2\}.
\end{align*}
\(C_p(D)\) is a tensor product of factors, one for each
crossing of \(D\). If we ignore \(d_v\),  \(\Cu_p(D)\) will
split into a direct sum of \(2^n\) summands, where \(n\) is the
number of ordinary crossings of \(D\). 
By assigning the summand \(C_p(D_r)\) to the oriented resolution of a
crossing and \(C_p(D_s)\) to its singular resolution, we get
 a bijection between summands and MOY states of \(D\). Comparing the
 grading shifts with the weights in Figure~\ref{Fig:Resolutions}, we
 see that 
\begin{equation*}
\Cu_p(D) \cong 
\bigoplus_\sigma
\Cu_p (D(\sigma))\{\mu(\sigma),0,-2\mu(\sigma)\},
\end{equation*}
where the diagram \(D(\sigma)\) is obtained by replacing each ordinary
crossing of \(D\) with either \(D_s\) or \(D_r\), depending on
\(\sigma\). Note that this is not quite the same the diagram \(D_\sigma\),
which is obtained by replacing each ordinary crossing with
 either \(D_s\) or the oriented resolution \(D_u\). 
 To remedy this discrepancy we use the following lemma, whose proof
 is given in the next section. 

\begin{lem}
\label{Lem:Marks}
Suppose \(D\) is a closed tangle diagram containing a crossing of type
\(D_r\), and let \(D'\) be the diagram obtained by replacing this
crossing by a pair of marks, as illustrated by the diagram \(D_u\) 
in Figure~\ref{Fig:Dr}. Then
\(\Cu_p(D)\) is quasi-isomorphic to \(\Cu_p(D')\) over \(R(D')\).
\end{lem}
\noindent
Applying the lemma repeatedly, we see that \(\Cu_p(D(\sigma))\) is
quasi-isomorphic to \(\Cu_p(D_\sigma)\). Thus
\(\Hu^+(D(\sigma)) \cong \Hu^+(D_\sigma)\), and the claim is proved.
\end{proof}

The similarity between the HOMFLY homology and the original Khovanov
homology \cite{Khovanov} is now evident. Like the chain complex used
to define the Khovanov homology, the summands of 
\(\Hu_p^+(D)\) naturally lie at the vertices of the ``cube of
resolutions'' of \(D\), each of whose vertices corresponds to a MOY
state. The components of the induced differential \(d_v^*\) correspond to
edges of the cube. 
This analogy can be used to give an alternate proof of
Theorem~\ref{Thm:Chi}, which is easily seen to be equivalent to the
statement below.

\begin{prop}
\( \displaystyle{\wt{P}(L) = 
\sum_{i,j,k} (-1)^{(k-j)/2} a^j q^i  \dim \Hu^{i,j,k}(L)}\) .
\end{prop} 

The argument is  similar to the proof that the
Euler characteristic of the Khovanov homology is the Jones polynomial, 
but with the MOY state model in place of the Kauffman state model.

\begin{proof} Recall that if 
 \(D\) is a closed braid diagram representing  \(L\),
\begin{equation*}
\Hu(L) = H(\Hu^+(D),d_v^*)\{-w+b,w+b,w-b\}.
\end{equation*}
Since \(d_v^*\) preserves both \(q\) and \(\gr_h\), the graded Euler
characteristic 
\begin{equation*}
\chi(\Hu(L)) = \sum_{i,j,k} (-1)^{(k-j)/2} a^j q^i  \dim \Hu^{i,j,k}(L)
\end{equation*}
can be computed from \(\Hu^+(D)\). We find
\begin{align*}
\chi(\Hu(L)) & = (-1)^{-b} a^{w+b} q^{-w+b}\sum_{i,j,k} 
(-1)^{(k-j)/2} a^j q^i  \dim
\Hu_p^{+ \ i,j,k}(D) \\
& = (aq^{-1})^w
\sum_{\sigma} (-q)^{\mu(\sigma)} (-aq)^b \sum_{i,j} (-1)^{j/2}a^jq^i \dim
\Hu^{i,j}(D_\sigma) \\
& = (aq^{-1})^w
\sum_{\sigma} (-q)^{\mu(\sigma)} \wt{P}(D_\sigma) \\
& = \wt{P}(L). 
\end{align*}
\end{proof}

\subsection{MOY decompositions}
\label{SubSec:MOYDecomp}

We now prove the various technical results
 used throughout the section. We begin with the proof of 
 Lemma~\ref{Lem:Marks}, which asserted that a ``crossing'' of type
 \(D_r\) was equivalent to its oriented resolution.
 
\begin{proof} (of Lemma~\ref{Lem:Marks}.)
Either \(D\) and \(D'\) have the same number of connected components,
or \(D'\) has one more component than \(D\). Suppose we are in the
first case. 
Then  \(X_i\) and \(X_k\) are
independent linear elements of the polynomial ring \(R(D)\), and we
can use Lemma~\ref{Lem:Exclude} to exclude the linear factor
\(X_k-X_i\) appearing in \(\Cu_p(D_r)\). We obtain a quasi-isomorphic
complex \(C'\) defined over the ring \(R(D)/(X_k-X_i)\). The ideal
generateds by \((X_k-X_i,X_l-X_j)\) and \((X_k-X_i,X_k+X_l-X_i-X_j)\)
are clearly equal, so  \(R(D)/(X_k-X_i) \cong
R(D')\). Then \(C'\) and \(\Cu_p(D')\) are Koszul factorizatiouns over
\(R(D')\) with the same Koszul matrices, so \(C'\cong \Cu_p(D)\).

Now suppose that replacing \(D_r\) with \(D_u\)
 increases the number of components in \(D\). In this case,
 \(X_i=X_k\) and \(X_j=X_l\)  in \(R(D)\), so \(R(D) \cong R(D')\). 
We compute 
\begin{align*}
p_i & = \frac{p(X_k)+p(X_l)-p(X_i)-p(X_j)}{X_k-X_i} \\
& = \frac{p(X_k)-p(X_i)}{X_k-X_i} + \frac{p(X_j+X_i-X_k)-p(X_j)}{X_k-X_i}
\\
& = p'(X_i) - p'(X_j)
\end{align*}
so  \(C_p(D_r)\) is given by the short factorization
 \(\{p'(X_i)-p'(X_j),0\}\). 

Recall that \(\Cu_p(D)\)  is obtained by tensoring \(C_p(D)\) with short
factorizations of the form \(\{p'(X_n),0\}\), where we pick one edge
\(n\) for each component of \(D\). By Proposition~\ref{Prop:MarkQI}, we
may assume  that the component containing
\(D_r\) has marked edge \(i\). Thus \(\Cu_p(D)\) has short
factors \(\{p'(X_i),0\}\) and \(\{p'(X_i)-p'(X_j),0\}\). Applying a Koszul
row operation, we see that this is isomorphic to a factorization with
short factors \(\{p'(X_i),0\}\) and \(\{p'(X_j),0\}\). If we choose
    \(j\) as the marked edge on the new component, this is the
factorization for \(\Cu_p(D')\).

\end{proof}

Next, we take up the task of proving the MOY decompositions stated in
Propositions \ref{Prop:MOYO}--\ref{Prop:MOYIII}. In each case, the
argument follows the proofs of the corresponding results in
\cite{KRI,KRII}, although some additional care is required for the
 MOY \(III\) move. The proof for the MOY \(O\) move is easiest. 

\begin{proof} (of Proposition~\ref{Prop:MOYO}). The diagram \(D_O'\)
    is obtained from \(D_O\) by deleting a small loop consisting of a
    single edge, labeled \(1\), attached at both ends to a single
    mark. The relation \(\rho(v)\) associated to this mark is \(0\),
    so \(R(D_O) \cong R(D_O') \otimes_{\Q} \Q[X_1]\). Both diagrams
    have the same set of crossings, so the only difference between
    \(\Cu_p(D_O)\) and \(\Cu_p(D_O')\) comes from the factor
     associated to the deleted component. This is precisely the
    factorization \(C_p(O)\) from the statement of the
    proposition.
\end{proof}

The argument for the MOY \(I\) move is not much harder. 

\begin{proof} (Of Proposition~\ref{Prop:MOYI}.) We start by
  considering the case when \(D_I\) is  the open diagram shown
  in Figure~\ref{Fig:MOYMoves}. Then \(R(D_I) \cong
  \Q[X_1,X_2,X_3]/(X_3+X_2-X_2-X_1) \cong \Q[X_1,X_2]\), while
  \(R(D_I') = \Q[X_1]\), so \(R(D_I)
  \cong  R(D_I')\otimes \Q[X_2]\). \(D_I'\) has no crossings, so
\(C_p(D_I') = R(D_I')\), while \(C_p(D_I)\) is the short factorization
$$\xymatrixcolsep{7pc} 
\xymatrix{ R\{1,-2,0\}  \ar@<0.4ex>[r]^{(X_3-X_1)(X_3-X_2)}&
  R\{-1,0,0\} \ar@<0.4ex>[l]^{p_{12}}},$$
where \(R=R(D_I)\). \(X_1=X_3\) in \(R(D_1)\), so the entry on
the upper arrow is \(0\).  To compute \(p_{12}\), we first take the
quotient 
 \(
p(X_3)+p(X_4)-p(X_2)-p(X_1)/(X_3-X_1)(X_3-X_2)\) in the ring
\(\Q[X_1,X_2,X_3,X_4]/(X_3+X_4-X_1-X_2)\)
 and then set \(X_2=X_4\). In other words, 

\begin{equation*}
p_{12}  = \frac{1}{X_1-X_4} \left[ 
\frac{p(X_1+X_2-X_4) - p(X_1)}{X_2-X_4} + \frac{p(X_4) -
  p(X_2)}{X_2-X_4} \right]\Biggl\lvert_{X_2=X_4} 
\end{equation*}
which  reduces to \((p'(X_1)-p'(X_2))/(X_1-X_2)\). Thus
\(C_p(D_I)\) is exactly the factorization \(C_p(I)\) described in the
statement of the proposition, and  \(C_p(D_I) \cong
C_p(D_I') \otimes_{\Q[X_1]} C_p(I)\). 

More generally,  
suppose that \(\overline{D}_I\) and \(\overline{D}_I'\) 
are formed by gluing a fixed graph
\(D\) to \(D_I\) and \(D_I'\). Then Lemma~\ref{Lem:Gluing} 
tells us that  \( \Cu_p(\overline{D}_I')  \cong
\Cu_p(D)\otimes_{\Q[X_1]} C_p(D_I')\), so 
\begin{align*}
\Cu_p(\overline{D}_I) & \cong \Cu_p(D)\otimes_{\Q[X_1]} C_p(D_I) \\
& \cong \Cu_p(D)\otimes_{\Q[X_1]} C_p(D_I') \otimes_{\Q[X_1]} C_p(I)
\\ & \cong \Cu_p(\overline{D}_I') \otimes  C_p(I)
\end{align*}
and the general case follows from the local one. 
\end{proof}

The proof of  the MOY \(II\) relation follows its counterpart in
\cite{KRII} almost verbatim.

\begin{proof} (Of Proposition~\ref{Prop:MOYII}.)
As before, we start by assuming that \(D_{II}\) and \(D_{II}'\) are
the open graphs shown in Figure~\ref{Fig:MOY2}. We label their edges
as shown in the figure. 
\(C_p(D_{II})\) is an order 2 Koszul factorization
 over the ring
\begin{equation*}
R  = \Q[X_1,\ldots,X_6]/(X_5+X_6-X_3-X_4,X_3+X_4-X_1-X_2)   \cong R_0[X_3]
\end{equation*}
where \(R_0 = \Q[X_1,X_2,X_5,X_6]/(X_5+X_6-X_1-X_2) \cong
\Q[X_1,X_2,X_5]\) is isomorphic to both \(R_e(D_{II})\) and
\(R_e(D_{II}')\). It is given by a Koszul matrix of the form
\begin{equation*}
\begin{pmatrix}
* & -(X_3-X_1)(X_3-X_2) \\
* & (X_3-X_5)(X_3-X_6)
\end{pmatrix}. 
\end{equation*}

\begin{figure}
\input{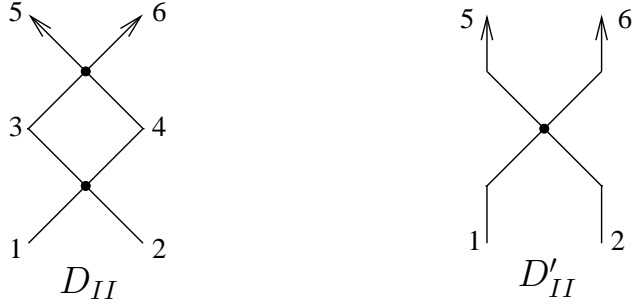}
\caption{\label{Fig:MOY2} Diagrams for the MOY \(II\) move.}
\end{figure}


We use the entry \(-(X_3-X_1)(X_3-X_2)\) in the first row to
exclude the internal variable \(X_3\). The result is  a new 
 factorization \(C_1\) which  is quasi-isomorphic 
 to \(C_p(D_{II})\) over the ring \(R_0\).
\(C_1\) is an order one Koszul factorization
defined over the ring 
\begin{equation*}
R_1 = R/(X_3^2-(X_1+X_2)X_3+X_1X_2).
\end{equation*} We can write
\(C_1 = \{P,Q\}\), where \(Q\) is obtained by
substituting \(X_3^2 = (X_1+X_2)X_3 - X_1X_2\) into the lower right entry
of the factorization above. We have
\begin{align*}
Q & = (X_1+X_2)X_3 -X_1X_2 -(X_5+X_6)X_3 +X_5X_6  \\
& = - X_1X_2 + X_5X_6 
 = -(X_5-X_2)(X_5-X_1).
\end{align*}
Thus although \(Q\) is {\it a priori} an element of \(R_1\), we find
that actually \(Q\in R_0\). Since  the product \(PQ =
w_p(D_{II})\) is also contained in \(R_0\),  \(P \in R_0\) as well. 
Thus \(C_1 \cong C_2 \otimes_{R_0} R_1\), where \(C_2\) is
the short factorization over \(R_0\) defined by the pair \(\{P,Q\}\). 
In other words, \(C_2 = C_p(D_{II}')\). 
Viewed as a module over \(R_0\), \(R_1 \cong R_0
 \oplus X_3 R_0\), so over \(R_0\)
\begin{equation*}
C_1  \cong C_2 \oplus X_3 C_2 
 = C_p(D_{II}') \oplus X_3 C_p(D_{II}').
\end{equation*}

Next, we check the grading shift of the two
summands. \(C_p(D_{II})\) is a direct sum of \(4\) copies of \(R\), with
grading shifts \(\{-2,0\},\{0,-2\},\{0,-2\}\), and \(\{2,-4\}\). When
we exclude \(X_3\) to get \(C_1\), we are left with two copies of
\(R_1\), with grading shifts  \(\{-2,0\}\) and \(\{0,-2\}\). Since
 \(R_1 \cong R_0 \oplus X_3 R_0 = R_0 \oplus
R_0\{2,0\}\),  \(C_1\) is a direct sum of \(4\) copies of
\(R_0\), with grading shifts \(\{-2,0\},\{0,0\},\{0,-2\}\),and
\(\{2,-2\}\). On the other hand, \(C_p(D_{II}')\) is a direct sum of
two copies of \(R_0\) with grading shifts \(\{-1,0\}\) and
\(\{1,-2\}\). Thus \(C_1\) must decompose as \(C_p(D_{II}')\{-1,0\}\oplus
C_p(D_{II}')\{1,0\}\). 

Finally, we consider the general situation, in which \(D_{II}\) and
\(D_{II}'\) are formed by attaching the diagrams shown in the figure
to an arbitrary graph \(D\). In this case the result follows from the
special case considered above, the local nature of the KR-complex
(Lemma~\ref{Lem:Gluing}), and the 
 fact that if \(A\sim B\) over \(R\), then  \(A\otimes_R C \sim B
 \otimes_R C\). 

\end{proof}

\begin{figure}
\includegraphics{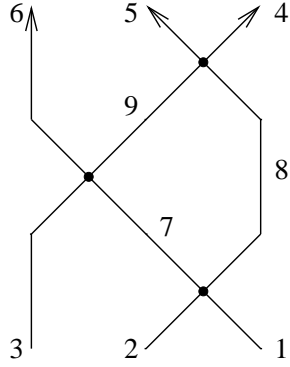}
\caption{\label{Fig:D3a} The diagram \(D_{IIIa}\)} 
\end{figure}
Lastly, we turn to the MOY \(III\) move. As usual, it suffices to
prove the statement of Proposition~\ref{Prop:MOYIII} for the graphs shown in
Figure~\ref{Fig:MOY3} and then appeal to the local nature of the
KR-complex to show that it holds in  general. We number the edges of
the diagram \(D_{IIIa}\) as shown in Figure~\ref{Fig:D3a}, and  label
the external edges of 
\(D_{IIIb}\),  \(D_{IIIa}'\), and  \(D_{IIIb}'\) to match.
All four diagrams  share the same potential
\begin{equation*}
W = p(X_4)+p(X_5) + p(X_6) - p(X_1)-p(X_2)-p(X_3)
\end{equation*}
and the same exterior ring
\begin{equation*}
R_0 = \Q[X_1,\ldots X_6]/(X_4+X_5+X_6-X_1-X_2-X_3). 
\end{equation*}
In \cite{KRI,KRII}, Khovanov and Rozansky
 introduce an additional factorization \(\Upsilon\) defined over
 \(R_0\) and  with potential \(W\). 
\(\Upsilon\) is a order two Koszul factorization given by the
 Koszul matrix
\begin{equation*}
\begin{pmatrix}
P_2 & s_2(X_1,X_2,X_3) - s_2(X_4,X_5,X_6) \\
P_3 & X_1X_2X_3 - X_4X_5X_6
\end{pmatrix}
\end{equation*}
where \(s_2(x,y,z) = xy+xz+yz\) is the degree two symmetric
polynomial. Since  \(W\) is symmetric in \(X_1,X_2,X_3\)
and \(X_4,X_5,X_6\), it is not difficult to see that it is in the
ideal generated by the symmetric differences
\(X_1+X_2+X_3-X_4-X_5-X_6\), \(s_2(X_1,X_2,X_3) - s_2(X_4,X_5,X_6)\),
and \(X_1X_2X_3 - X_4X_5X_6\). The first symmetric difference vanishes in
\(R_0\), so we can find 
 \(P_2,P_3\in R_0\) so that 
\begin{equation*}
(s_2(X_1,X_2,X_3) - s_2(X_4,X_5,X_6))P_2 + (X_1X_2X_3 - X_4X_5X_6)P_3
= W.
\end{equation*}
in \(R_0\). 
 The choice of \(P_2\)
and \(P_3\) is not unique, but since \(s_2(X_1,X_2,X_3) -
s_2(X_4,X_5,X_6)\) and \(X_1X_2X_3 - X_4X_5X_6\) are relatively
prime, Lemma~\ref{Lem:KTwist} implies that any two choices are related
by a twist. For definiteness, we fix some  values of  \(P_2\)
and \(P_3\) which are symmetric in \(X_1,X_2,X_3\) and (separately) in
\(X_4,X_5,X_6\). 

In \cite{KRII}, Khovanov and Rozansky exhibit  homotopy equivalences
\begin{align*}
f& : C_p(D_{IIIa})_+ \to \Upsilon_+ \oplus C_p(D_{IIIb}')_+ \\
g&: C_p(D_{IIIb})_+ \to \Upsilon_+ \oplus C_p(D_{IIIa}')_+
\end{align*}
Our goal is to show that \(f\) and \(g\) are 
quasi-isomorphisms. We will prove

\begin{prop}
\label{Prop:MOYIIIB} 
\(C_p(D_{IIIa}) \sim \Upsilon \oplus
C_p(D_{IIIb}')\) over \(R_0\).
\end{prop}
For the moment, let us  assume that this proposition is true.
Observe that \(D_{IIIa}'\) is essentially the same graph as \(D_{IIIa}\),
but with the labels on edges \(1\) and \(3\) and \(4\) and \(6\)
reversed.  Since  \(\Upsilon\) is symmetric in both \(X_1,X_2,X_3\) and
\(X_4,X_5,X_6\),  \(C_p(D_{IIIa}')\) will be
quasi-isomorphic to \(\Upsilon \oplus C_p(D_{IIIb})\). It
follows that \(C_p(D_{IIIa}) \oplus C_p(D_{IIIb})\) and \(
C_p(D_{IIIa}') \oplus C_p(D_{IIIb}')\) are both quasi-isomorphic to 
\( \Upsilon \oplus C_p(D_{IIIb}) \oplus C_p(D_{IIIb}') \). Thus
Proposition~\ref{Prop:MOYIII} is implied by  Proposition~\ref{Prop:MOYIIIB}.  
To prove the latter,  we follow step-by step the argument given in 
\cite{KRII}.

The factorization \(C_p(D_{IIIa})\) is defined over the ring 
\(R=R(D_{IIIa}) = \Q[X_1,\ldots,X_9]/I \), where 
\begin{equation*}
I = (X_7 +X_8 -  X_1 - X_2, \ X_6+X_9 - X_3-X_7, \ X_4+X_5 -
X_8-X_9). 
\end{equation*}
We use these relations to eliminate \(X_1\), \(X_7\), and \(X_8\),
thus expressing \(R\) as a polynomial ring in  variables
\(X_2,X_3,X_4,X_5,X_6,X_9\). In this ring, \(C_p(D_{IIIa})\) is an order
three Koszul factorization, with Koszul matrix
\begin{equation*}
\begin{pmatrix}
* & (X_2-X_8)(X_2-X_7) \\
* & (X_3-X_9)(X_3-X_6) \\
* & (X_9-X_4)(X_9-X_5)
\end{pmatrix}.
\end{equation*}
Eliminating \(X_7\) and \(X_8\), this becomes
\begin{equation*}
\begin{pmatrix}
* & (X_2+X_9-X_4-X_5)(X_2+X_3-X_6-X_9) \\
* & (X_3-X_9)(X_3-X_6) \\
* & (X_9-X_4)(X_9-X_5) 
\end{pmatrix}.
\end{equation*}
We use the right-hand entry of the last row to exclude the internal
variable \(X_9\). The result is an order two Koszul factorization
\(C_1\) over
the ring \(R_1 = R/(X_9-X_4)(X_9-X_5)\), with Koszul matrix
\begin{equation*}
\begin{pmatrix}
* & X_9(X_3-X_6) +X_4X_5 +(X_2-X_4-X_5)(X_2+X_3-X_6) \\
* & (X_3-X_9)(X_3-X_6) \\
\end{pmatrix}.
\end{equation*}
After a row operation in which we add the 
 bottom entry in the right-hand row to the top one, we get
\begin{equation*}
\begin{pmatrix}
* & (X_2-X_4)(X_2-X_5) +(X_2+X_3-X_4-X_5)(X_3-X_6) \\
* & (X_3-X_9)(X_3-X_6) \\
\end{pmatrix}. 
\end{equation*} 
More explicitly, this factorization is represented by the following
diagram
\begin{equation*}
\xymatrixcolsep{4pc}
\xymatrixrowsep{4pc}
\xymatrix{R_1 \ar@<0.4ex>[r]^{b_1} \ar@<0.4ex>[d]^{a_2}  & R_1
  \ar@<0.4ex>[l]^{a_1} \ar@<0.4ex>[d]^{a_2} 
\\
R_1 \ar@<0.4ex>[r]^{-b_1} \ar@<0.4ex>[u]^{b_2}
 & R_1 \ar@<0.4ex>[l]^{-a_1}
\ar@<0.4ex>[u]^{b_2}}
\end{equation*}
where \(a_1\) and \(a_2\) are unknown, \(b_1 = (X_2-X_4)(X_2-X_5)
+(X_2+X_3-X_4-X_5)(X_3-X_6)\), and \(b_2 = (X_3-X_9)(X_3-X_6)\).

We now think of \(C_1\) as an object of \(GMF(R_0)\), and 
\(R_1\) as a free module of rank \(2\) over \(R_0\). 
 Following \cite{KRII}, we choose an explicit basis
\(\{1,X_9+X_3-X_4-X_5\}\) for the two copies of \(R_1\) in the lower
row of the diagram,  and the basis \(\{1,X_9-X_3\}\) for the two
copies of \(R_1\) in the upper row. 
With respect to these bases, \(C_1\) takes the form
\begin{equation*}
\xymatrixcolsep{4pc}
\xymatrixrowsep{4pc}
\xymatrix{R_0\oplus R_0 \ar@<0.4ex>[r]^{B_1'} \ar@<0.4ex>[d]^{A_2}  & R_0\oplus R_0
  \ar@<0.4ex>[l]^{A_1'} \ar@<0.4ex>[d]^{A_2} 
\\
R_0\oplus R_0 \ar@<0.4ex>[r]^{-B_1} \ar@<0.4ex>[u]^{B_2}
 & R_0\oplus R_0 \ar@<0.4ex>[l]^{-A_1}
\ar@<0.4ex>[u]^{B_2}}
\end{equation*}
where
 \(A_1,A_1',B_1,B_1',A_2,B_2\) are  \(2\times2\) matrices over \(R_0\)
representing multiplication by \(a_1,b_1,a_2\) and \(b_2\). 
The pairs \(A_1\) and \(A_1'\) and
\(B_1\) and \(B_1'\) represent the same linear maps with respect to
 different bases, so they are conjugate.
For \(B_1\), this is 
irrelevant ---  \(X_9\) does not appear in \(b_1\), 
so \(B_1\) is a multiple of the identity map:
\begin{equation*}
B_1= B_1' = \begin{pmatrix} x & 0 \\ 0 & x
\end{pmatrix}
\end{equation*}
where \(x = (X_2-X_4)(X_2-X_5) +(X_2+X_3-X_4-X_5)(X_3-X_6)\).
Direct computation shows that 
\begin{equation*}
B_2 = \begin{pmatrix} 0 & y \\ z & 0 
\end{pmatrix}
\end{equation*}
where \(y = (X_3-X_4)(X_3-X_5)(X_3-X_6)\) and \(z = X_6-X_3\).

\begin{lem}
\(A_1\) and \(A_1'\) may be expressed in the form
\begin{equation*}
A_1 = \begin{pmatrix}  a& yc/z  \\ c & a-qc
\end{pmatrix}
\quad 
A_1' = \begin{pmatrix}  a-qc& yc/z  \\ c & a
\end{pmatrix}
\end{equation*}
where \(q= X_4+X_5-2X_3\). 
\end{lem}

\begin{proof}
Suppose that \(A_1 = \left( \begin{smallmatrix} a&b  \\ c & d
\end{smallmatrix} \right)\). 
Changing basis, we see that 
\begin{equation*}
A_1' = \begin{pmatrix}  a - qc & b + q(a-d)-q^2c \\ c & d + qc 
\end{pmatrix}
\end{equation*}
where \(q = X_4+X_5-2X_3\). The component of
\(d_+d_-+d_-d_+\) which maps from the bottom right corner of the
square to the upper left must vanish, so \(A_1'B_2 = B_2A_1\), or,
more explicitly 
\begin{equation*}
\begin{pmatrix}  a - qc & b + q(a-d)-q^2c \\ c & d + qc 
\end{pmatrix}
\begin{pmatrix} 0 & y \\ z & 0 
\end{pmatrix} =
\begin{pmatrix} 0 & y \\ z & 0 
\end{pmatrix}
\begin{pmatrix} a&b  \\ c & d
\end{pmatrix}.
\end{equation*}
Multiplying out and equating terms, we find that we must have \(d = a -qc
\) and \(yc = zb\). 
\end{proof}

\begin{lem}
\( {\displaystyle
A_2 = \begin{pmatrix} -xc' & \beta \\ \gamma & -xc'
\end{pmatrix}},\)
where \(z c' = c\). 
\end{lem}

\begin{proof}
Suppose that \(A_2 = \left( \begin{smallmatrix} \alpha & \beta   \\
  \gamma & \delta 
\end{smallmatrix} \right)\). Inspecting the component of \(d_-d_+
  +d_+d_-\) which maps the lower left-hand corner of the diagram to
  itself, we see that \(A_2B_2 + B_1A_1 = W \cdot \text{Id}\), or
\begin{equation*}
\begin{pmatrix} \alpha & \beta   \\
  \gamma & \delta 
\end{pmatrix}
\begin{pmatrix} 0 & y \\ z & 0 
\end{pmatrix} +
\begin{pmatrix} x & 0 \\ 0 & x
\end{pmatrix}
\begin{pmatrix}  a& yc/z  \\ c & a-qc
\end{pmatrix} =
\begin{pmatrix} W & 0 \\ 0 & W
\end{pmatrix}
\end{equation*}
Inspecting the off-diagonal elements, we find that \(z \delta + x c =
y \alpha + xyc/z = 0 \). Since \(x\) and \(z\) are relatively prime,
we must have \(c = z c'\), \(\delta = -x c'\) for some \(c' \in
R_0\). Substituting into the second equation, we see that \(\alpha =
\delta\). 
\end{proof}
Thus we can write \(A\) and \(A'\) in the form
\begin{equation*}
A_1 = \begin{pmatrix}  a& yc'  \\ zc' & a-qzc'
\end{pmatrix}
\quad 
A_1' = \begin{pmatrix}  a-qzc'& yc'   \\ zc' & a
\end{pmatrix}. 
\end{equation*}
Consider the map \( H\) from the upper right-hand copy of \(R_0\oplus
R_0\) to the lower left given by the
matrix \(H = \left( \begin{smallmatrix} -c' & 0   \\
  0 & -c' \end{smallmatrix} \right)\). The twisted factorization
\(C_2(H)\) has the same positive differentials as \(C_2\), but the
negative differentials are given by matrices
\begin{equation*}
A_1(H) = \begin{pmatrix}  a& 0  \\ 0 & a-qzc'
\end{pmatrix}
\quad 
A_1'(H) = \begin{pmatrix}  a-qzc'& 0   \\ 0 & a
\end{pmatrix} \quad 
A_2(H) = \begin{pmatrix} 0 & \beta \\ \gamma & 0
\end{pmatrix}. 
\end{equation*}
It is now clear that \(C_2(H)\) decomposes as a direct sum: one
summand  consists of the first copies of \(R_0\) in the top row and the
second copies in the bottom, and the other of the second copies in the top
 and the first in the bottom. Both summands are order two Koszul 
factorizations over \(R_0\), with Koszul matrices 
\begin{equation*}
\begin{pmatrix}
a-qzc' & x \\ \gamma & z
\end{pmatrix}
\quad \text{and} \quad 
\begin{pmatrix}
a & x \\ \beta & y
\end{pmatrix}
\end{equation*}
Recalling that \(y = (X_3-X_4)(X_3-X_5)(X_3-X_6)\),
 \(z = X_3-X_6\), and \(x = z(X_4+X_5-X_2-X_3) + (X_2-X_4)(X_2-X_5)\),
 and using row operations to simplify the right-hand columns, we see
 that these are equivalent to Koszul matrices of the form
\begin{equation*}
\begin{pmatrix}
* & (X_2-X_4)(X_2-X_5) \\ * & X_3-X_6
\end{pmatrix} \quad \text{and} \quad
\begin{pmatrix}
* & s_2(X_1,X_2,X_3)-s_2(X_4,X_5,X_6) \\ * & X_1X_2X_3 - X_4X_5X_6
\end{pmatrix}.
\end{equation*}
By Lemma~\ref{Lem:KTwist}, the first factorization is a twisted version
of 
\begin{equation*}
\begin{pmatrix}
p_{1245} & (X_2-X_4)(X_2-X_5) \\ p_{36} & X_3-X_6
\end{pmatrix} 
\end{equation*}
where 
\(p_{1245}  = \frac{p(X_4)+p(X_5)-p(X_1)-p(X_2)}{X_4+X_5-X_1-X_2} \)
and 
\(p_{36}  = \frac{p(X_6)-p(X_3)}{X_6-X_3}.\)
Arguing as in the proof of Lemma~\ref{Lem:Marks}, we see that this
factorization is quasi-isomorphic to \(C_p(D_{IIIb}')\). 
Likewise, Lemma~\ref{Lem:KTwist} shows that the second
factorization is a twisted version of \(\Upsilon\). 

To recap, we have shown that 
\(C_p(D_{IIIa}) \sim C_1\). By Lemma~\ref{Lem:Twist},
 \(C_1 \sim C_1(H)\), which decomposes
into a direct sum of  two order two Koszul complexes. Finally, a further 
application of Lemmas~\ref{Lem:Twist} and \ref{Lem:KTwist} shows that
these are quasi-isomorphic to \(C_p(D_{IIIb}')\) and
\(\Upsilon\). We leave it as an exercise for the reader
to check that the two summands have the correct bigrading. 
 \qed

\section{Relation between \(\H\) and \(\H_N\)}
\label{Sec:d+}

We are now in a position to address
 the relation between  the HOMFLY and \(sl(N)\)
homologies. Here is our main result.

\begin{thrm}
\label{Thm:Hp}
Suppose \(L \subset S^3\) is a link, and let \(i\) be a marked
 component of \(L\). For each \(p \in \Q[x]\), there is
 a spectral sequence \(E_k(p)\) with \(E_1(p) \cong \H(L)\) and 
\(E_{\infty}(p) \cong \H_p(L,i)\). For all \(k>0\), the isomorphism
 type of \(E_k(p)\) is an invariant of the pair \((L,i)\). 
\end{thrm}

\begin{cor}
The isomorphism type of \(\H_p(L,i)\) is an invariant of  \((L,i)\). 
\end{cor}

The relation of these sequences with the various gradings may be
summarized as follows.
Let \(d_k(p):E_k(p) \to E_k(p)\) be the \(k\)th differential in the
 sequence. If \(p(x) = x^{N+1}\), then \(d_k(p)\) 
is homogenous of degree \((2Nk ,-2k,2-2k)\) with respect to the
triple grading on \(\H(K)\). In particular, each \(d_k(p)\)
preserves the  grading \(\gr_N' = q + 2N \gr_h\). The
grading on \(E_\infty(p)\) induced by \(\gr_N'\) 
 is equal to the polynomial grading \(\gr_N = 
q + (N-1) \gr_h\) on \(\H_N(L,i)\).

For general values of \(p\), \(d_k(p)\) is no longer homogenous with respect
to the \(q\)-grading, but it is still the case that \(d_k(p)\) shifts
\(\gr_h\) by \(-k\) and  \(\gr_v\) by \(1-k\). Thus the \(d_k(p)\)
are all homogenous of degree \(1\) with respect to the grading
\(\gr_-\). The grading induced by \(\gr_-\) on the \(E_{\infty}\) term
is equal to the homological grading on \(\H_p(k)\).

A few other remarks on the theorem are in order. First
it is possible to prove an analogous result for the unreduced
homology. The argument is very similar to the one in the reduced case,
except that we don't need to worry about keeping track of a 
marked edge. Second, in terms of invariance, the spectral sequence suffers
from the same drawback as the HOMFLY homology --- we can show that any
two diagrams representing the same link  give rise to isomorphic spectral
sequences, but not that the isomorphism is canonical. Finally,
we expect that  \(\H_p(L,i)\)
should be
 determined by the order of vanishing of \(p'(x)\) at \(x=0\), and that
\(\Hu_p(L)\) should be determined by the multiplicities of the roots of
\(p'(x)\).
(This idea has its source in the work of Gornik \cite{Gornik}.
Although we will not pursue   it here, some
 supporting evidence has been provided by 
Mackaay and Vaz \cite{MackaayVaz}.)
In  particular, it seems unlikely that the set of all homologies
\(\H_p(L,i)\), \(p \in \Q[x]\) contains more information than is
present in the \(sl(N)\) homologies. 

\subsection{Definition and basic properties}
\label{Subsec:DiffDefs}
We now construct the spectral sequence \(E_k(p)\). Given
a link \(L\) with a marked component, we fix a braid diagram
\(D\) representing \(L\), and an edge \(i\)  belonging to the marked
component. The complex \(\Cr_p(D,i)\) is endowed with differentials
\(d_+,d_-\), and \(d_v\). Since \(D\) is a closed diagram, all three
differentials anticommute. It follows that \(\H^+_p(D,i)\) inherits a pair of
 anticommuting differentials \(d_-^*\) and \(d_v^*\).
 \(d_-^*\) lowers \(\gr_h\) by \(1\) and preserves
 \(\gr_v\), while \(d_v^*\) raises \(\gr_v\) by \(1\) and preserves
 \(\gr_h\). Thus the triple \((\H^+_p(D,i),d_v^*,d_-^*)\) defines a double
 complex with total differential \(d_{v-} = d_v^*+d_-^*\) and total grading
 \(\gr_- = \gr_v-\gr_h\). 
Like any double complex,  \((\H^+_p(D,i),d_{v-})\) comes with two
natural filtrations: a {\it horizontal filtration} induced by
\(\gr_h\), and a {\it vertical filtration} induced by
\(\gr_v\). 
\begin{defn}
\(E_k(p)\) is the spectral sequence induced by the
horizontal filtration on the complex \((\H^+_p(D,i),d_{v-})\). 
\end{defn}
As in the definition of \(\H\), we shift the triple grading on
\(E_k(p)\) by a factor of \(\{-w+b-1,w+b-1,w-b+1\}\), where \(w\) and
\(b\) are the writhe and number of strands in the diagram \(D\). 
With this normalization, the first claim of Theorem~\ref{Thm:Hp} is 
 easily verified. \(E_0(p) =
\H^+_p(D,i)\{-w+b-1,w+b-1,w-b+1\}\), and the differential 
\(d_0(p):E_0(p) \to E_0(p)\) is  the part of \(d_v^*+d_-^*\) which
preserves \(\gr_h\). In other words, \(d_0(p) = d_v^*\), so
\begin{equation*}
E_1(p)  = H(\H^+_p(D,i), d_v^*) \{-w+b-1,w+b-1,w-b+1\} 
 \cong \H(L,i). 
\end{equation*}
To complete the proof of the theorem, we must  show that the total
homology 
\begin{equation*}
H(\H^+_p(D,i),d_{v-}) \cong \H_p(D,i)
\end{equation*}
 and  that the sequence is
an invariant of the pair \((L,i)\). Before we doing this, we pause to
discuss some  elementary properties of \(E_k(p)\). First, note that 
when \(p\) is a linear polynomial, the differential \(d_-\) is
identically zero, and the spectral sequence converges
trivially to \(\H_p(D,i) \cong \H(D)\). Thus the sequence is only
interesting when \(\deg p > 1\). For the rest of the section,
we will assume that this is the case. 

Next, we address the issue of gradings. 

\begin{lem} The differential \(d_k(p)\) is homogenous of degree
  \(-k\) with respect to \(\gr_h\) and degree \(1-k\) with respect to
  \(\gr_v\). In addition, 
if \(p(x) = x^{N+1}\), then
  \(d_k(p)\) is homogenous of degree \(2Nk\) with respect to the
  \(q\)--grading.  
\end{lem}

\begin{proof}
When \(p(x) = x^{N+1}\), this follows immediately from the fact
 \(d_-\) and \(d_v\) are
 homogenous of degrees \((2N,-2,0)\) and  \((0,0,2)\)
with respect to the triple grading on \(\H_p^+(D,i)\). For general
 values of \(p\),  \(d_-\) is no longer  homogenous with
 respect to the \(q\)-grading, but its behavior with respect to the
 homological gradings remains unchanged. 
\end{proof}

When \(p=x^{N+1}\), the differentials \(d_k(p)\) all preserve the
quantity \(\gr'_N = q + 2N \gr_h\), 
so the
graded Euler characteristic of 
 \(H(\H_p^+(D,i),d_{v-})\) with respect to \(\gr_N'\) will be the same as
 the graded Euler characteristic of \(E_1(p)\).  
 Using Theorem~\ref{Thm:Chi} we compute that
\begin{align*}
\chi(E_1(p)) & = \sum (-1)^{\gr_-} q^{\gr_N'}\dim \H(L,i) \\
& = \sum_{i,j,k} (-1)^{k-j} q^{i+Nj} \dim \H^{i,j,k}(L,i) \\
& =  \sum_{i,j,k} (-1)^{k-j} a^j q^{i} \dim \H^{i,j,k}(L,i)
   \bigl \vert_{a=q^N} \\
& = P_L(q^N,q).
\end{align*}
Since the graded Euler characteristic of \(H(\H_p^+(D,i),d_{v-})\)
 is given  by the \(sl(N)\) polynomial, it's at
 least plausible that the homology should
agree with \(\H_N\).

Next, we consider the relation between \(E_k(p)\) for different values
of \(p\).
In the original complex \(\Cr_p(D,i)\), the underlying group and the
differentials \(d_+\)
and \(d_v\) are  independent of \(p\). Thus
 \(E_0(p) = \H_p^+(D,i)\) is  independent of \(p\), as is 
 \(d_0 = d_v^*\). It follows that we can view \(E_1(p) = \H(L,i)\) as being
 equipped with an infinite dimensional family of differentials
 \(d_1(p)\) --- one for each \(p \in \Q[x]\). 

\begin{lem} \( d_1(ap+bq)  = ad_1(p) + b d_1(q)\).  
\end{lem}

\begin{proof}
Denote the differential \(d_-\) on \(\Cr_p(D)\) by \(d_-(p)\). We
claim that 
\begin{equation*}
d_-(ap+bq) = ad_-(p) + bd_-(q).
\end{equation*}
 Since \(d_1(p)\) is 
the map induced by \(d_-(p)\), the claim implies the statement of the lemma.
To prove it, observe that
for an elementary tangle \(D\), the potential 
\begin{equation*}
W_p = p(X_3)+p(X_4) - p(X_1)-p(X_2)
\end{equation*}
 satisfies \(W_{ap+bq}
= aW_p + b W_q\). 
The coefficients of \(d_-(p)\) are quotients of \(W_p\)
by fixed polynomials, so they are also linear in \(p\). Finally, it is
easy to see that the linearity property is preserved under tensor
product, so the claim holds.
\end{proof}

\begin{cor}
\label{Cor:AntiCom}
For all \(p,q \in \Q[x] \), \(d_1(p)\) and \(d_1(q)\) anticommute.
\end{cor}

\begin{proof}
We have 
\( d_1(p)d_1(q) + d_1(q)d_1(p) = (d_1(p) + d_1(q))^2 = d_1(p+q)^2 = 0\).
\end{proof}

\subsection{The total homology}
\label{SubSec:+HomCalc}
Our goal in this section is to calculate the homology group
\(H(\H_p^+(D,i),d_{v-})\) to which \(E_k(p)\) converges.
(Throughout, we continue to assume that \(\deg p > 1\).)
 To do this, we use two more spectral
sequences --- one which converges to \( H(\H_p^+(D,i),d_{v-})\), and
another which can be used to calculate \(\H_p(D,i)\). The key point is
to show that both of these sequences converge at the \(E_2\) term, and
that the \(E_2\) terms agree. 

We start off with some notation.
 Suppose that \(C^*\) is graded matrix factorization
with potential \(0\), so that \(d_+\) anticommutes with \(d_-\). Then
we define
\begin{equation*}
\Hpm(C^*) = H(H(C^*,d_+),d_-^*).
\end{equation*}
When \(C^*=C_p(D,i)\), we abbreviate this to \(\Hpm(D,i)\). Both
\(d_+\) and \(d_-\) are homogenous with respect to the homological grading on
\(C\) (albeit with different degrees), so this grading descends to a
well-defined grading on \(\Hpm(C^*)\). 

\begin{lem}
\label{Lem:WeakE}
If \(C^*\) and \(D^*\) are quasi-isomorphic factorizations
with potential \(0\), then \(\Hpm(C^*) \cong \Hpm(D^*)\).
\end{lem}

\begin{proof}
From the definition of a quasi-isomorphism, we know that the complex
\((H(C^*,d_+),d_-^*)\)  is isomorphic to \( (H(D^*,d_+),d_-^*)\)
\end{proof}

\begin{lem} 
\label{Lem:OneGrade}
Suppose \(D\) is a closed braid graph on \(b\) strands. If,
 then \(\H^\pm_p(D,i)\) is
  supported in horizontal grading \(\gr_h = 1-b\).
\end{lem}

\begin{proof}
As in section~\ref{SubSec:SingHom}, 
we use the MOY relations to induct on the complexity of \(D\). 
In the base case, \(D\) is a single
circle, and  \(\H^\pm_p(D,i) \cong \Cr_p(D,i) \cong
\Q\) is supported in grading \(\gr_h = 0\). 

For the induction step, we use Lemma~\ref{Lem:Hao} to see that \(D\)
can be related to diagrams of lesser complexity using the MOY
moves. In fact, we claim that \(D\) can be simplified by a MOY move
which has the marked edge \(i\) as an external edge. To see this, write
\(D\) as the closure of an open braid in such a way that \(i\) is one
of the edges which appear in the closure. Then \(i\) will be an
external edge for any MOY \(II\) or \(III\) move provided by the
lemma. If \(i\) lies on a loop which could be eliminated by a MOY
\(O\) move, we ignore it and simplify using some other MOY
move. Finally, if \(i\) is a small loop which is about to be eliminated by
a MOY \(I\) move, \(i\) must be the rightmost strand in \(D\). 
In this case, we consider the mirror image  \(\overline{D}\)
of \(D\). It's easy
to see that \(\Cr_p(D,i) \cong \Cr_p(\overline{D},i)\), and 
the marked edge in \(\overline{D}\) is on the leftmost strand. 
We now  use  Lemma~\ref{Lem:Hao} to simplify \(\overline{D}\) as before.

Thus in order to prove the lemma, it is enough to check that if the
statement holds for the less complex diagram(s) in each of the
four MOY
moves, it also holds for the more complex one. 
For example, suppose \(D_O\) and \(D_{O}'\) are related by a MOY
\(O\) move, so that \(C_p(D_O) \cong C_p(D_O') \otimes_\Q C_p(O)\). Then
\(\Cr_p(D_O,i) \cong \Cr_p(D_O',i) \otimes_\Q C_p(O)\), and  \(d_+ = 0 \) on 
\(C_p(O)\), so \(\H^+_p(D_O,i) \cong \H^+_p(D_O',i) \otimes_\Q
C_p(O)\). By the Kunneth formula, 
\begin{equation*}
\H^\pm_p(D_O,i) \cong \H_p^\pm (D_O',i) \otimes_{\Q} H(C_p(O),d_-).
\end{equation*}
\(D_O'\) has one fewer strand than \(H_p(D)\), so \(H_p(D_O',i)\) is
supported in \(\gr_h = 2-b \) by the induction hypothesis. When
\(\deg(p) > 1\), \(  H(C_p(O),d_-)\) is
supported in \(\gr_h = -1\). Thus  \( \H^\pm_p(D_O,i)\) is supported in
\(\gr_h = 1-b\) as claimed.

Similarly, if \(D_I\) and \(D_I'\) are related by a MOY \(I\) move, 
\(C_p(D_I,i) \cong C_p(D_I',i) \otimes_{\Q[X_1]} C_p(I)\)
and \(\H^+_p(D_I,i) \cong \H^+_p(D_I',i) \otimes_{\Q[X_1]}
C_p(I)\). The complex \((C_p(I),d_-)\)  has the form
$$\xymatrixcolsep{3pc} 
\xymatrix{  \Q[X_1,X_2]\{0,-2\}  &
  \Q[X_1,X_2]\{0,-0\}\ar[l]_{p_{12}'}} $$
where \( p_{12}' = (p'(X_1) -p'(X_2))/(X_1-X_2)\). When \(\deg p > 1\),
its homology 
is a free module over \(\Q[X_1]\), supported in grading
\(\gr_h = -1\). As in the previous case, we apply the Kunneth formula
to conclude that \( \H^\pm_p(D_I,i)\) is supported in \(\gr_h =
1-b\).

Next, suppose that \(D_{II}\) and \(D_{II}'\) are related by a MOY
\(II\) move which takes place away from the marked edge \(i\). 
 Proposition~\ref{Prop:MOYII} tells us
that 
\begin{equation*}
C_p(D_{II}) \sim C_p(D_{II}')\{-1,0\} \oplus C_p(D_{II}')\{1,0\}
\end{equation*} over a ring \(R \) which contains \(\Q[X_i]\) as a
subring. 
By Corollary~\ref{Cor:QIRed}, it follows that 
\begin{equation*}
C_p(D_{II},i) \sim C_p(D_{II}',i)\{-1,0\} \oplus
C_p(D_{II}',i)\{1,0\}.
\end{equation*} 
Applying  Lemma~\ref{Lem:WeakE}, we see that
\begin{equation*}
\H^\pm_p(D_{II},i) \sim \H^\pm_p(D_{II'},i)\{-1,0\} \oplus  
\H^\pm_p(D_{II'},i)\{1,0\}.
\end{equation*}
Since both diagrams have the same number of strands, the  result
 follows from the induction hypothesis. 
The argument for the MOY \(III\) move is very similar, and is left to
the reader. 
\end{proof}

\begin{cor}
\label{Cor:OneGrade} If  \(D\) is a closed braid on \(b\) strands,
 \(\H^\pm_p(D,i)\) is supported in horizontal grading
 \(\gr_h = 1-b  \). 
\end{cor}

\begin{proof} 
 Since we are only taking homology with respect to \(d_+\) and
 \(d_-\), \(\H^\pm_p(D,i)\) decomposes as a direct sum over MOY states
 of \(D\) ({\it cf.} Lemma~\ref{Lem:States}): 
\begin{equation*}
\H^\pm_p(D,i)  \cong \bigoplus_{\sigma}
\H_p^\pm(D_\sigma,i)\{\mu(\sigma),0,-\mu(\sigma)\}.
\end{equation*}
If \(D\) is a braid on \(b\) strands, 
each diagram \(D_\sigma\) will be a braid graph on \(b\) strands.  
There are no shifts in \(\gr_h\), so each summand 
 \( \H_p^\pm(D_\sigma,i)\) is supported in \(\gr_h = 1-b\). 
\end{proof}

\begin{prop}
\label{Prop:Seq1}
If \(D\) is a closed braid, then 
\begin{equation*}
 H(\H^+_p(D,i),d_{v-}) \cong H(\H^\pm_p(D,i),d_v^*).
\end{equation*} 
\end{prop}

\begin{proof}
We compute the total homology using 
the spectral sequence induced by the vertical 
filtration on \((\H_p^+(D,i),d_v^*,d_-^*)\).
In this sequence, \(d_0 = d_-^*\),
so the \(E_1\) term is 
\begin{equation*}
H(\H^+_p(D,i),d_-^*)  =  \H^\pm_p(D,i)  .
\end{equation*}
By Corollary~\ref{Cor:OneGrade}, this group is supported in a
single horizontal grading. 

The differential \(d_1\) is the induced map \(d_v^*: \H^\pm_p(D,i)  \to
\H^\pm_p(D,i)\). Thus the \(E_2\) term of the sequence is the group
\(H(\H^\pm_p(D,i),d_v^*)\). 
For \(n>1\), the differential 
 \(d_n\) raises \(\gr_v\) by \(n\) and \(\gr_h\) by
\(n-1\). Since the \(E_1\) term (and thus the \(E_2\) term) is supported in
a single horizontal grading, \(d_n \equiv 0\) for all \(n>1\), and
 the sequence converges at the \(E_2\) term. 
{\it A priori}, this
implies that the graded group
\(H(\H^\pm_p(D,i),d_v^*)\) is isomorphic to the
associated graded group of \( H(\H^+_p(D,i),d_{v-})\).
In fact, for each
value of the homological grading \(\gr_-\), the former group is
supported in a unique value of the filtration grading \(\gr_v\). Thus the
two groups are canonically isomorphic.
\end{proof}

Our next task is to relate the group \(H(\H^\pm_p(D,i),d_v^*)\) to
\(\H_p(D,i)\). To do so, we use a  slightly different
spectral sequence. Recall that if \((C^*,d_{\pm})\) is a matrix
  factorization with potential \(0\), we can form the total
  differential \(d_{tot} = d_++d_-\). 

\begin{lem}
\label{Lem:TotSeq}
If \(C^*\) is a matrix factorization with potential
\(0\), there is a spectral sequence with \(E_2\) term
\(\Hpm(C^*)\) which converges to \(H(C^*,d_{tot})\). 
\end{lem}

\begin{proof} We define an increasing filtration of \((C^*,d_{tot})\)
  by \(F^n = \oplus_{i<n} C^i \oplus \ker d_+^{n}\), where \(d_+^n\)
  denotes the component of \(d_+\) which maps   \(C^n\)  to \(C^{n+1}\). 
 The \(E_0\) term of the spectral sequence induced by this filtration
is the associated  graded complex
\begin{equation*}
E_0^n = 
\frac{F^n}{F^{n-1}} \cong \frac{C^{n-1}}{\ker d_+^{n-1}} \oplus \ker
d_+^{n}. 
\end{equation*}

The differential \(d_0:E_0^n \to E_0^n\) is given by \(d_0(x,y) =
(d_-y,d_+x) \). If \(d_0(x,y) = 0 \), then \(d_+x = 0\), so \(x = 0\) as an
element of \(C^{n-1}/\ker  d_+^{n-1}\). Conversely, \(d_+d_-y = -
d_-d_+y = 0\) for any \(y \in  \ker  d_+^{n}\), so \(d_- y = 0\) as an
element of \(C^{n-1}/\ker  d_+^{n-1}\). Thus 
\begin{equation*}
\ker d_0 = \{(0,y) \ts
\vert \ts y \in \ker  d_+^{n}\} \cong \ker  d_+^{n}.
\end{equation*}
Similarly, \(\im d_0 \cong \im d_+^{n-1}\), so 
  \(E_1^n =  H(E_0^n,d_0) \cong  H^n(C^*,d^+)\).

Next, we consider the differential \(d_1:E_1^n\to E_1^{n-1}\). 
An element of \(E_1^n\) can be  represented by \(x \in \ker
d_+^n\), and  \(d_1x \) is the image of \(d_{tot}x = d_- x \) in
\(E_1^{n-1}\). In other words, \(d_1\) is given by
\(d_-^*:H^n(C^*,d^+)\to H^{n-1}(C^*,d^+) \), and  the \(E_2\) term is  
\(H(H(C^*,d_+),d_-^*) = \Hpm(C^*)\).
\end{proof}

\begin{lem}
If \(D\) is a closed braid graph, then 
\(H(\Cr_p(D,i),d_{tot}) \cong \H^\pm_p(D,i).\)
\end{lem}

\begin{proof}
We apply the sequence of the preceding lemma to 
 the complex \((\Cr_p(D,i),d_{tot})\). By
Lemma~\ref{Lem:OneGrade}, the \(E_2\) term 
is supported in a single homological grading, and thus in a
single filtration grading as well. 
It follows that the sequence has converged at the \(E_2\)
term. As in the proof of Proposition~\ref{Prop:Seq1}, the fact that
 the \(E_{\infty}\) term is supported in a single filtration grading
 implies that it is canonically isomorphic to the total homology. 
\end{proof}


\begin{prop}
\label{Prop:Seq2}
If \(D\) is a closed braid, then 
\(\H_p(D,i) \cong H(\H^\pm_p(D,i),d_v^*).\)
\end{prop}

\begin{proof}

By definition, \(\H_p(D,i)\) is the homology of the complex
\(( H(\Cr_p(D,i),d_{tot}),d_v^*)\). 
To prove the proposition, it suffices to 
 show that this complex is isomorphic to 
 \((\H^\pm_p(D,i),d_v^*)\).  
The  complex \((\Cr_p(D,i),d_{tot})\) splits as
a direct sum over MOY states of \(D\), so by  the previous lemma, the
underlying group  \(H(\Cr_p(D,i),d_{tot})\) is isomorphic to
\(\H^\pm_p(D,i)\). 


To see that the differentials are identified under this isomorphism, 
we must check that  the following diagram commutes:
$$ \xymatrixrowsep{1pc}
\xymatrix{ H(\Cr_p(D,i),d_{tot}) \ar^{d_v^*}[r] \ar@{=}[d]&
   H(\Cr_p(D,i),d_{tot}) \ar@{=}[d] \\
\H^\pm_p(D,i)\ar^{d_v^*}[r] & \H^\pm_p(D,i)}
$$
To show this, we filter \(C_{tot} = (\Cr_p(D,i),d_{tot})\) 
as in the proof of Lemma~\ref{Lem:TotSeq}. Since \(d_v\) commutes with
\(d_+\), the map \(d_v:C_{tot} \to C_{tot}\) preserves the
filtration. It thus induces a morphism of spectral sequences
\((d_v)_k: E_k \to E_k\) which converges to  \(d_v^*:H_{tot}
\to H_{tot}\). In particular, the map \((d_v)_2\) is the induced map \( d_v^*:
\H^\pm_p(D,i) \to \H^\pm_p(D,i)\). Since both sequences converge at the
\(E_2\) term,  \((d_v)_2\) is also equal to the associated
graded map of \(d_v^*: H_{tot} \to H_{tot}\).  \(H_{tot}\) is
supported in a single filtration grading, so the two maps are actually
equal. 
\end{proof}
\noindent
To sum up, we have
\begin{prop}
The spectral sequence \(E_k(p)\) converges to \(\H_p(D,i)\). The
grading \(\gr_-\) on \(E_k(p)\) corresponds to the homological grading on
\(\H_p(D,i)\), and if \(p = x^{N+1}\), the grading \(\gr_N'\) on
\(E_k(p)\) corresponds to the polynomial grading \(\gr_N\) on
\(\H_p(D,i)\). 
\end{prop}

\begin{proof}
The first statement is immediate from 
Propositions~\ref{Prop:Seq1} and \ref{Prop:Seq2}, so we just need to 
check that the gradings agree. The triple grading on \(E_1(p)\) is the
 grading on \(\Cr_p(D,i)\), shifted by
 \(\{-w+b-1,w+b-1,w-b+1\}\), where \(w\) and \(b\) are the writhe and
 number of strands in \(D\). We write \(\gr_-(E)\) and \(\gr_-(C)\)
 for the shifted and unshifted gradings, so that
\begin{align*}
\gr_-(E)  & = \gr_v(E) -\gr_h(E) \\ &= \gr_v(C) -\gr_h(C) -b + 1
\quad \text{and} \\
\gr_N'(E) & = q(E) + 2N \gr_h(E)  \\ & = 
q(C) + 2N \gr_h(C) + (N+1)(b-1) + (N-1)w
\end{align*}
On \(E_{\infty}(p)\), \(\gr_-(C)\) and \(\gr_N'(C)\) agree with the
gradings 
\(\gr_-\) and \(\gr_N'\) on the total homology
\(H(\H^+_p(D,i),d_{v-})\). The gradings on the latter group can be computed
using the spectral sequence of Proposition~\ref{Prop:Seq1}, whose 
\(E_\infty\) term is supported in \(\gr_h = 1-b\). 
Substituting, we find that
\begin{align*}
\gr_-(E) & = \gr_v(C) -(1-b)-b+1 \\ &= \gr_v(C),
\end{align*}
which is the homological grading on \(\H_p(D,i)\), and 
\begin{align*}
\gr_N'(E) & = q(C) - (N-1)(b-1) +(N-1)w \\ & = q(C) + (N-1) \gr_h(C) +
(N-1)w \\  & = \gr_N(C) + (N-1)w,
\end{align*} 
which is the polynomial grading on \(\H_N(D,i)\). 
\end{proof}

\subsection{Change of marked edge}
We now turn to the last part of Theorem~\ref{Thm:Hp} --- the invariance
of \(E_k(p)\). The first step is to show that for a fixed braid
diagram \(D\) representing \(L\), \(E_k(p)\) depends only on  the
component of \(L\) containing the marked edge \(i\), and not on \(i\)
itself. 

Suppose for the moment that \(D\) is an arbitrary tangle diagram. 
If \(i\) is an internal edge of \(D\), multiplication by \(X_i\) defines an
endomorphism \(X_i:C_p(D) \to C_p(D).\) Viewing \(C_p(D)\) as an
object of \(Kom(GMF_w(R(D)))\), we can form the mapping cone
\(Cone(X_i)\), which is also an object of \(Kom(GMF_w(R(D)))\). 
We would like to view \(Cone(X_i)\) as a factorization over the ring  
\(R_i(D) = R(D)/(X_i)\). To do so, we observe that since \(X_i\) is a
generator of the polynomial ring \(R(D)\), there is an
inclusion \(R_i(D) \subset R(D)\) with the property that \(R_i(D)[X_i]
\cong R(D)\). 

\begin{lem}
\(Cone(X_i)\) is homotopy equivalent to \(\Cr_p(D,i)\) 
in  \(Kom(GMF_w(R_{i}(D)))\). 
\end{lem}

\begin{proof}
As a matrix factorization, \(Cone(X_i) = C_p(D) \oplus C_p(D)\). 
The vertical 
differential is given by \(d_v(x,y) = (d_vx,d_vy +(-1)^{\gr_v x} X_i
\cdot x)\). Let \(\pi_0:C_p(D) \to C_p(D,i)\) be the projection, and
define \(\pi:Cone(X_i) \to \Cr_p(D,i)\) by \(\pi(x,y) =
\pi_0(y)\). Since \(C_p(D)\) is a free module over \(R(D) =
R_i(D)[X_i]\), we have an injection \(\iota_0: \Cr_p(D,i) \to
C_p(D)\) modeled on the inclusion \(R_i(D) \subset R(D)\). We extend
this to an inclusion \(\iota:\Cr_p(D,i) \to Cone(X_i) \) given by
\(\iota(z) = (0, \iota_0(z))\). The composition \(\pi\iota\) is the
identity map, and \(\iota \pi\) is homotopic to the identity via a
homotopy \(H:Cone(X_i) \to Cone(X_i)\) given by 
\(H(x,y) = ( (y - \iota \pi y)/X_i,0)\).   
\end{proof}

\begin{lem}
\label{Lem:EdgeHom}
Suppose \(D\) is a tangle diagram, and let \(j\) and \(k\) be two
edges of \(D\) which belong to the same component of the underlying
tangle and are separated by a single ordinary crossing. (For example,
the edges labeled \(j\) and \(k\)
 in Figure~\ref{Fig:Tangles}.) Then 
\(X_j \) and \(X_k\) are homotopic morphisms from \(C_p(D)\) to
itself. 
\end{lem}

\begin{proof}
Suppose that \(D\) is the elementary diagram \(D_+\). The complex
\(C_p(D_+)\) has the form

$$\xymatrixcolsep{7pc} 
\xymatrixrowsep{4pc}
\xymatrix{ 
R\{0,-2,0\}  \ar@<0.4ex>[r]^{(X_k-X_i)}&
  R\{0,0,0\} \ar@<0.4ex>[l]^{p_{1}} \\
R\{2,-2,-2\}  \ar@<0.4ex>[r]^{-(X_k-X_i)(X_k-X_j)} 
\ar[u]^{(X_j-X_k)}&
  R\{0,0,-2\} \ar@<0.4ex>[l]^{p_{12}} \ar[u]_{1}} $$
where \(R=R(D_+). \)
 The  homotopy \(H\) is given by the vertical arrows in the
diagram below.
$$\xymatrixcolsep{7pc} 
\xymatrixrowsep{4pc}
\xymatrix{ 
R\{0,-2,0\}  \ar@<0.4ex>[r]^{(X_k-X_i)}   \ar[d]_{1} &
  R\{0,0,0\} \ar@<0.4ex>[l]^{p_{1}} \ar[d]^{(X_j-X_k)} \\
R\{2,-2,-2\}  \ar@<0.4ex>[r]^{-(X_k-X_i)(X_k-X_j)} &
  R\{0,0,-2\}. \ar@<0.4ex>[l]^{p_{12}} } $$
It's easy to see that \(H\) commutes with  \(d_+\) and
\(d_-\) and that \(d_vH+Hd_v = X_j-X_k\). 
Since \(X_j-X_k = X_l-X_i\) in \(R(D_+)\), multiplication by \(X_i\) and
\(X_l\) are also homotopic. The reader can easily check that there is
a similar homotopy when \(D\) is the elementary diagram \(D_-\).

For a general diagram \(D\), the local nature of the KR-complex
implies that we can write \(C_p(D) \cong
C_p(D_{\pm})\otimes_{R(D_{\pm})} C_p(D')\), where \(D_{\pm}\) is the
crossing separating \(j\) from \(k\), and \(D'\) is the rest of the diagram. 
 In terms of this decomposition,
 \(X_j:C_p(D) \to C_p(D)\) can be written as \(X_j \otimes 1\) and
similarly for \(X_k\). Clearly
\(f\sim g\) implies \(f\otimes 1 \sim g \otimes 1\), so we are done.
\end{proof}

Next, we need a general result from homological algebra.

\begin{lem}
Suppose \(\mathcal{A}\) is an additive category.
If \(f:A \to B\) and \(g:A\to B\) are homotopic in
\(Kom(\mathcal{A})\), then \(Cone(f) \cong Cone(g)\).  
\end{lem}

\begin{proof}
If \(H:A \to B\) is the homotopy from \(f\) to \(g\), then the map
\(H: Cone(f) \to Cone(g)\) defined by \(h(a,b) = (a,b-Ha)\) is an
isomorphism with inverse \(h^{-1}(a,b) = (a,b+Ha)\). 
\end{proof}
\noindent

By repeatedly applying the lemmas, we see that if \(i\) and \(j\)
belong to the same  component of \(L\), then 
\(\Cr_p(D,i)\) and \(\Cr_p(D,j)\) are homotopy equivalent as objects of
 \(Kom(GMF_0(\Q))\).

\begin{prop} 
\label{Prop:MoveMark}
Suppose \(D\) is a diagram representing a link \(L\)
 and that \(i\) and \(j\) are edges
of \(D\) which belong to the same component of \(L\). If we denote by
 \(E_k(p,D,i)\)  the spectral sequence associated to the
 pair \((D,i)\), then \(E_k(p,D,i) \cong E_k(p,D,j)\) for all \(k>0\). 
\end{prop}

\begin{proof}
Let  \(f:\Cr_p(D,i) \to \Cr_p(D,j)\)
be  a homotopy equivalence. Recall the 
functor \(\mathcal{H}^+: GMF_0(\Q) \to Kom(\Q)\), which takes a
factorization \(C^*\) to the complex \((H^+(C^*),d_-^*)\).
We apply \(\mathcal{H}^+\) to \(f\) and get a 
 homotopy equivalence \(f^+: \H^+_p(D,i) \to \H^+_p(D,j)\)
 in the category \(Kom(Kom(\Q))\) ({\it i.e.} double complexes over
 \(\Q\)). \(f^+\) respects the
 horizontal filtration on \(\H^+_p\), so it induces a map of spectral
 sequences \(f^+_k: E_k(p,i) \to E_k(p,j)\). The map \(f_1^+: E_1(p,i)
 \to E_1(p,j)\) is the induced map 
\begin{equation*}
(f^+)^*:  H(\H^+_p(D,i),d_v^*) \to H(\H^+_p(D,i),d_v^*).
\end{equation*}
Since \(f^+\) is a homotopy
 equivalence with respect to  \(d_v^*\), this map is an isomorphism.
This proves the claim when \(k=1\).   Finally, it is a well-known proverty
 of spectral sequences that if \(f_r^+\) is an isomorphism,
then \(f_k^+\) is an isomorphism for all \(k>r\) as well.  (See {\it
  e.g.} Theorem 3.4 of \cite{UsersGuide}.) 
\end{proof}

\subsection{Invariance under Reidemeister moves}

\begin{figure}
\input{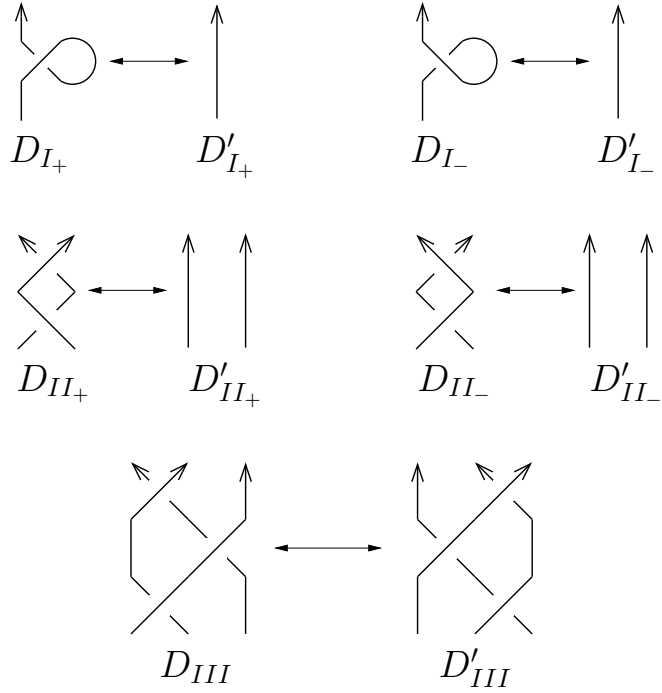}
\caption{\label{Fig:RMoves} Braidlike Reidemeister moves.}
\end{figure}

The final step in the proof of Theorem~\ref{Thm:Hp} is to show that
\(E_k(p)\) remains invariant when we vary the diagram \(D\). 
Following \cite{KRII}, we make some preliminary simplications of the
problem. By assumption, the diagram representing 
 \(L\) is a braid diagram. Any two braid diagrams representing
the same link \(L\) can be joined by a sequence of the five moves shown
in Figure~\ref{Fig:RMoves}, so it is enough to prove that \(E_k(p)\)
is invariant under these moves. Also, 
using Proposition~\ref{Prop:MoveMark}, we can
assume that the marked edge \(i\) does not participate in the move.

In what follows, it will be important to keep track of the
category we are working in. To help with this, we introduce the
following notation, writing
\begin{align*}
\cc_j & = Kom(GMF_{w(D_j)}(R_e(D_j))) & \ccp_j & = Kom(Kom(R_e(D_j))) \\
\ck_j & = Kom(hmf(R_e(D_j))) & \ckp_j & = Kom(Mod(R_e(D_j)))
\end{align*}
where \(j=I,II,III\). The same symbols without subscripts 
 indicate the corresponding category for a closed diagram. 
({\it e.g.} \(\cc =Kom(GMF_0(\Q))\).) There is a commutative square of
 functors 
$$\xymatrixcolsep{2pc} 
\xymatrixrowsep{2pc}
\xymatrix{ 
\cc  \ar[r]^{\mathcal{H}^+} \ar[d]^{\mathcal{F}} & \ccp \ar[d]^{\mathcal{F}} \\
\ck  \ar[r]^{\mathcal{H}^+} & \ckp} $$
where \(\mathcal{F}\) is the forgetful functor which
corresponds to ignoring the differential
\(d_-\) (or \(d_-^*\)) and \(\mathcal{H}^+\) is the functor
which corresponds to taking homology with respect to \(d_+\).

Suppose that \(\od_j\) and \(\od_j'\) are closed diagrams related by the
\(j\)th Reidemeister move. Below, we will show that there are morphisms 
\begin{equation*}
\sigma_{j}:  \H^+_p(\od_{j},i) \to  
\H^+_p(\od_{j}',i) \quad (j=I_\pm,II_\pm,III)
\end{equation*}
in the category \(\ccp\) with the property that
\(\mathcal{F}(\sigma_j)\) is a homotopy equivalence in \(\ckp\). 
This is sufficient to prove the theorem. Indeed, arguing as in 
proof of Proposition~\ref{Prop:MoveMark}, we see that \(\sigma_j\)
induces a morphism of spectral sequences \((\sigma_j)_k:E_k(P,\od_j,i)
\to E_k(P,\od_j',i)\) which is an isomorphism for \(k>0\). 

Most of the work involved in constructing the \(\sigma_j\) and
showing they are homotopy equivalences has
already been done by Khovanov and Rozansky. 
In \cite{KRII}, they prove
invariance of the HOMFLY homology by exhibiting homotopy equivalences 
\begin{equation*}
\rho_{j}: \mathcal{F}(C_p(D_{j})) \to \mathcal{F}(C_p(D_{j}'))
\end{equation*}
in the category \(\ck_j\). From the local nature of the KR-complex, it
follows that there are homotopy equivalences 
\begin{equation*}
\rho_j\otimes 1:  \mathcal{F}(\Cr_p(\od_{j},i)) \to  \mathcal{F}(
\Cr_p(\od_{j}',i)) 
\end{equation*}
in \(\ck\). The morphism 
\(\sigma_j\) will be derived from \(\rho_j\), in
the sense that \( \mathcal{F}(\sigma_j) = \mathcal{H}^+(\rho_j\otimes
1)\). 
\vskip0.07in
\noindent{\bf Reidemeister \(I\) move:} In this case,
 we can work directly in the category \(\cc_I\). 

\begin{lem}
There are morphisms
\(\rho_{I_{\pm}}: C_p(D_{I_\pm}) \to C_p(D_{I_\pm}')\)
in  \(\cc_{I}\) with the property that
\(\mathcal{F}(\rho_{I_\pm})\) is a homotopy equivalence. 
\end{lem}

\begin{proof}
The ring \(R_e(R_{I_{\pm}})\) is isomorphic to \( \Q[X_1]\). 
Since \(R_I'\) has no crossings,
 \(C_p(R_I') = \Q[X_1]\) as well. 
Arguing as in the proof of Proposition~\ref{Prop:MOYI},
we see that the complex  \(C_p(I_+)\) has the form 
$$\xymatrixcolsep{7pc} 
\xymatrixrowsep{4pc}
\xymatrix{ 
\Q[X_1,X_2]\{0,-2,0\}  \ar@<0.4ex>[r]^{0}&
  \Q[X_1,X_2] \{0,0,0\} \ar@<0.4ex>[l]^{p'(X_1)-p'(X_2)} \\
\Q[X_1,X_2] \{2,-2,-2\}  \ar@<0.4ex>[r]^{0} 
\ar[u]^{X_2-X_1}&
  \Q[X_1,X_2] \{0,0,-2\}. \ar@<0.4ex>[l]^{p'_{12}} \ar[u]_{1}} $$

The morphism \(\rho_{I_+}\) takes the copy of \(\Q[X_1,X_2]\) in the
top left to \(C_p(R_I') = \Q[X_1]\) by subsituting \(X_2=X_1\), and
is zero elsewhere. 
Similarly, \(C_p(I_-)\) is the  complex
$$\xymatrixcolsep{7pc} 
\xymatrixrowsep{4pc}
\xymatrix{ 
\Q[X_1,X_2]\{0,-2,2\}  \ar@<0.4ex>[r]^{0}&
\Q[X_1,X_2] \{-2,0,2\} \ar@<0.4ex>[l]^{p'_{12}} \\
\Q[X_1,X_2]\{0,-2,0\}  \ar@<0.4ex>[r]^{0} 
\ar[u]^1&
 \Q[X_1,X_2] \{0,0,0\}, \ar@<0.4ex>[l]^{p'(X_1)-p'(X_2)}
 \ar[u]_{X_2-X_1}} $$
and the morphism \(\rho_{I_-}\) takes the copy of \(\Q[X_1,X_2]\) in the
top right to \(C_p(R_I') = \Q[X_1]\) by subsituting \(X_2=X_1\), and
is zero elsewhere.
The reader can easily verify that both \( \phi_{I_+}\) and \( \phi_{I_-}\)
are morphisms in \(\cc_{I}\) and that their
restrictions to  \(\ck_I\)  are homotopy
equivalences. (In fact, they are  the morphisms \(\rho_{I_{\pm}}\)
  defined in  \cite{KRII}.)
\end{proof}

By the local nature of the KR-complex, there are morphisms
\(\rho_{I_{\pm}}\otimes 1: \Cr_p(\od,i) \to \Cr_p(\od,i)\) 
in \(\cc\) which restrict to homotopy
equivalences in \(\ck\). Finally, the morphism
\begin{equation*}
\sigma_{I_\pm}:  \H^+_p(\od_{I_\pm},i) \to  \H^+_p(\od_{I_\pm}',i)
\end{equation*}
is defined to be \(\mathcal{H}^+(\rho_{I_\pm}\otimes1)\). The fact that
\(\mathcal{F}(\sigma_{I_\pm})\) is a homotopy equivalence follows from
the relation \(\mathcal{F}\mathcal{H}^+ = \mathcal{H}^+
\mathcal{F}\). \qed
\begin{figure}
\input{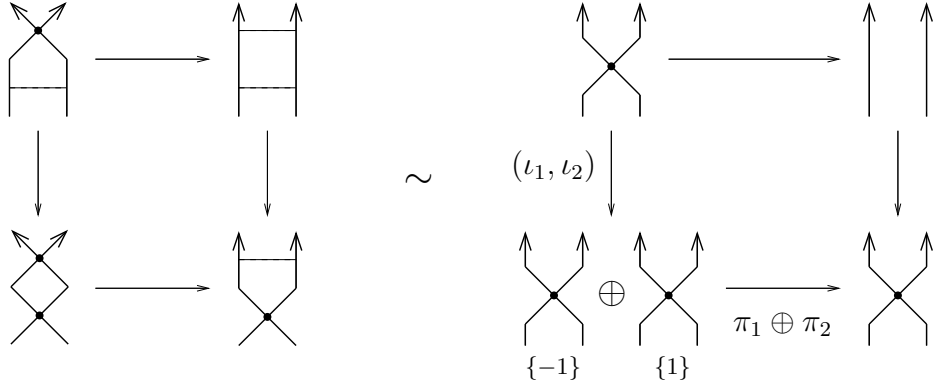}
\caption{\label{Fig:RII} The complex \(\H_p^+(D_{II_+},i).\) }
\end{figure}
\vskip0.05in 
\noindent{\bf Reidemeister \(II\) move:} 
As shown in Figure~\ref{Fig:RMoves}, there are two different versions of
the oriented Reidemeister \(II\) move. We  discuss only  the first
one --- the proof for the other is virtually identical.  
We can represent \(\Cr_p(\od_{II_+},i)\)  by the
diagram on the left-hand side of Figure~\ref{Fig:RII}. Each corner of
the square is an object of \(\cc\), and the edges represent 
additional components of \(d_v\) going between them. We apply the functor
\(\mathcal{H}^+\) to get \(\H^+_p(\od_{II_+},i)\).
Using the MOY \(II\) decomposition on each
factorization in the complex at the bottom left of the square (and
Lemma~\ref{Lem:Marks} on the other corners)  we see that
\(\H^+_p(\od_{II_+},i)\) has the form shown on the right-hand side of the
figure. 
\begin{figure}
\input{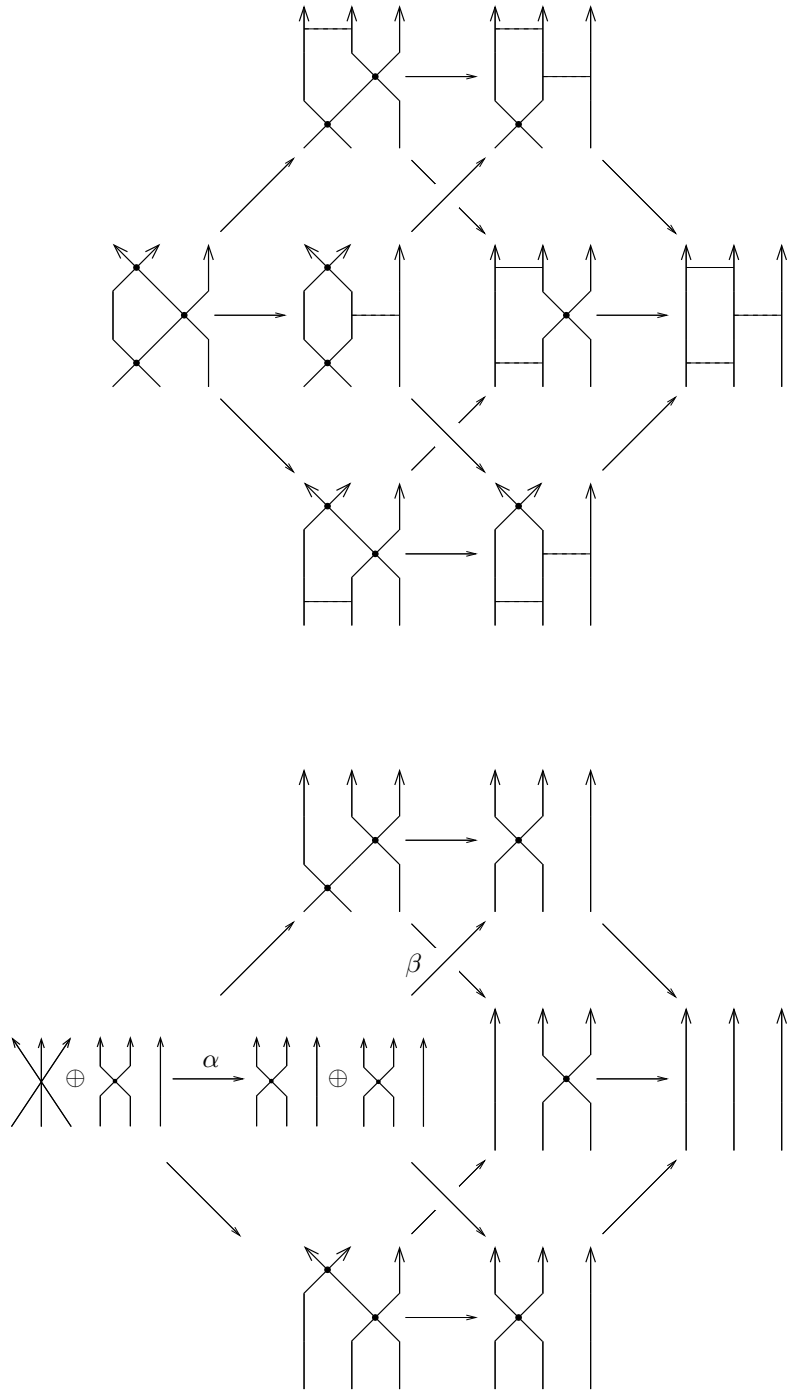}
\caption{\label{Fig:RIII} The complex \(\H_p^+(D_{III},i)\).}
\end{figure} 
Consider the maps \(\iota_1: \H^+_p(\doublepoint,i) \to 
\H^+_p(\doublepoint,i) \) (from the upper left-hand corner to the first
summand on the bottom left) and \(\pi_2: \H^+_p(\doublepoint,i) \to 
\H^+_p(\doublepoint,i) \) (from the second summand to the bottom
right). 
We claim that both  \(\iota_1\) and \(\pi_2\) are isomorphisms in
\(\ccp\). To see
this, we return to \(\Cr_p(D_{II},i)\), and apply the functor
\(\mathcal{F}\). In \(\ck\), we can use the MOY \(II\) decomposition
directly, without applying \(\mathcal{H}^+\) first. We get a diagram
like that on the right-hand side of Figure~\ref{Fig:RII}, with
corresponding morphisms \(\wt{\iota}_1\) and \(\wt{\pi}_2\). The main
ingredient in the
proof of invariance under the Reidemeister \(II\) move in \cite{KRII}
is to show that \(\wt{\iota}_1\) and \(\wt{\pi}_2\) are isomorphisms. 
From this,
it follows that \(\mathcal{H}^+(\wt{\iota_1}) = \mathcal{F}(\iota_1)\) and 
\(\mathcal{H}^+(\wt{\pi_2}) = \mathcal{F}(\pi_2)\) are isomorphisms. But 
if \(f\) is a morphism in \(\ccp \) with the property that
\(\mathcal{F}(f)\) is an isomorphism in \(\ckp\), then \(f\) is an isomorphism
in \(\ccp\). (In plainer language, a chain map that is an
isomorphism at the level of modules is an isomorphism.) This proves
the claim. 

At this point,
a standard cancellation argument like that used in
 the proof of invariance under
the Reidemeister \(II\) move in \cite{Khovanov} or \cite{KRI}
 shows that \(\H^+_p(\od_{II},i)\) is homotopy equivalent to \(
\H^+_p(\od_{II}',i)\) in \(\ccp\). \qed


\vskip0.05in 

\begin{figure}
\includegraphics{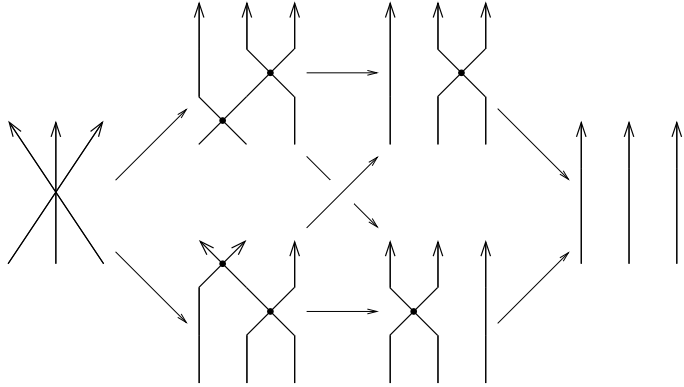}
\caption{\label{Fig:RIII2} Simplified version of  \(\H_p^+(D_{III},i)\).}
\end{figure}

\noindent{\bf Reidemeister \(III\) move:}
The argument in this case is similar to the one for the Reidemeister
\(II\) move. We start
out with the complex \(\Cr_p(\od_{III},i)\), which has the form
illustrated in the top half of Figure~\ref{Fig:RIII}. After
applying the functor \(\mathcal{H}^+\) and using the MOY \(II\) and
\(III\) decompositions, we get the diagram for \(\H^+_p(\od_{III})\)
shown in the bottom half
 of the figure. Our first claim is that the maps labeled
\(\alpha\) and \(\beta\) are isomorphisms. The proof is the same
as it was for the Reidemeister \(II\) move --- we consider the
analogous decomposition of \(\Cr_p(\od_{III},i)\) in the category
\(\ck\), where Khovanov and Rozansky
 proved that the corresponding maps \(\wt{\alpha}\) and
\(\wt{\beta}\) are isomorphisms. (This is the first part of the proof of 
Proposition 8 in \cite{KRII}). Canceling the summands connected by
\(\alpha\) and \(\beta\), we see that \(\H^+_p(\od_{III},i)\) is
homotopy equivalent (in \(\ccp\)) to a complex \(C\) of the form shown in
Figure~\ref{Fig:RIII2}.

An analogous simplification of  \(\H^+_p(D_{III}')\) shows
that it is homotopy equivalent  to a complex \(C'\)
which also has the form shown in Figure~\ref{Fig:RIII2}. To be
precise, it has the same six subquotients as \(C\). {\it A priori},
however, the morphisms between them may be different. We claim that in
reality, this is not the case  --- the maps in \(C\) and \(C'\)
corresponding to a fixed arrow in the diagram are same up to
multiplication by a nonzero element of \(\Q\). It follows
that \(C \cong C'\) in \(\ccp\), which gives us the desired homotopy
equivalence between \(\H^+_p(D_{III},i)\) and \(\H^+_p(D_{III}',i)\).

To prove the claim, we go back to the proof of the Reidemeister
\(III\) move in \cite{KRII}. There, Khovanov and Rozansky consider the
complexes of the open diagrams \(C_p(D_{III})\) and \(C_p(D_{III}')\)
in the category \(\ck_{III}\). As we have described above, they show
that they are homotopy equivalent to complexes \(C_o\) and \(C_o'\) of
the form illustrated in Figure~\ref{Fig:RIII2}. Moreover, they show
that the morphisms in \(C_o\) and \(C_o'\) corresponding to a fixed
edge in the figure are  nonzero multiples of each other. Going from
this statement to our claim is just a matter of applying
functoriality. More precisely, suppose \(f\) and \(f'\) are morphisms in
\(C\) and \(C'\) associated to some edge of the diagram, and that
\(f_o\) and \(f_o'\) are the corresponding morphisms in \(C_o\) and
\(C_o'\). Then 
\begin{equation*}
\mathcal{F}(f) = \mathcal{H}^+(f_o\otimes1) \quad \text{and} \quad 
 \mathcal{F}(f') = \mathcal{H}^+(f_o'\otimes1)
\end{equation*}
It follows that \(\mathcal{F}(f)\) is a nonzero multiple of
\(\mathcal{F}(f')\). Since a morphism of complexes is determined
by its action on modules, \(f\) is a nonzero multiple
of \(f'\). \qed


\section{An additional  sequence}
\label{Sec:d-1}

We now turn our attention to the spectral sequence described in 
Theorem~\ref{Thm:d-1}, which  is a special case of the following 

\begin{thrm}
Suppose \(L \subset S^3\) is an \(\ell\)--component link, and let
\(U^\ell\) be the \(\ell\)--component unlink. 
There's a spectral sequence \(E_k(-1)\) which has 
\(E_2\) term \(\H(L)\) and converges to \(\H(U^{\ell})\). 
The differentials 
in this sequence raise the cohomological grading \(\gr_+\) by
\(1\) and preserve the polynomial grading \(\gr_{-1}' = q - 2
\gr_h\).
\end{thrm}
\noindent  More precisely, the differential \(d_k\) is homogenous of
degree \((2-2k,2-2k,2k)\) with respect to the triple grading on
\(\H(L)\).

Compared with the sequences of the
preceding section, this construction of \(E_k(-1)\) is quite
simple. Indeed, the fact that such a sequence should exist 
 is well-known to
experts in the field. It is more surprising, however, that the
behavior of this sequence with respect to the triple grading should
so closely match the behavior predicted in \cite{superpolynomial} for the
``cancelling differential'' \(d_{-1}\).

\begin{proof}
We represent \(L\) by a braidlike diagram \(D\), and
consider the globally reduced complex \(C_r(D)\) defined in
section~\ref{SubSec:UnRed}. The triple \((C_r(D),d_+,d_v)\) is
a double complex with respect to the bigrading \((\gr_h,\gr_v)\). The
total differential on this complex is \(d_++d_v\), and the total
grading is \(\gr_h+\gr_v = \gr_+\). In addition, since \(d_+\) and
\(d_v\) have 
degrees \((2,2,0)\) and \((0,0,2)\) 
with respect to the triple grading on \(C_r(D)\), \(d_++d_v\)
preserves the polynomial grading \(\gr_{-1}' = q - 2 \gr_h\). 

The spectral sequence of the theorem is induced by the
horizontal filtration on the complex \((C_r(D),d_++d_v)\). 
This sequence has \(E_2\) term \(H(H(C_r(D),d_+),d_v^*) = \H(L)\), and
converges to the total homology \(H(C_r(D),d_v+d_+)\). 

To compute the latter group, we consider the spectral sequence induced
by the vertical filtration on \(C_r(D)\).
 The \(E_1\) term of this sequence is
\(H(C_r(D),d_v)\). 
Recall that \(C_r(D)\) is a tensor product of
factors, one for each crossing of \(D\):
\begin{equation*}
C_r(D) = \bigotimes_c C_r(D_c). 
\end{equation*}
 If we ignore the differential
\(d_+\), the factor associated to a crossing with sign \(\pm1\) 
has the form
\begin{equation*}
C_r(D_c) = M_c\{\pm1,\mp1,\mp1\} \oplus N_c\{0,\pm1,\mp1\}
\end{equation*} 
where
\begin{align*}
M_c & = R_r(D)\{1,-1,-1\}\xrightarrow{(X_k-X_j)} R_r(D)\{-1,-1,1\} \\
N_c & = R_r(D)\{0,-1,-1\}\xrightarrow{\phantom{XX}1\phantom{XX}}
R_r(D)\{0,-1,1\}.  
\end{align*}
The complex \(N_c\) is contractible, so \(C_r(D)\) is the direct sum of 
\begin{equation*}
M =  \bigotimes_c M_c\{\pm1,\mp1,\mp1\}
\end{equation*}
and a contractible complex. It follows that \(H(C_r(D),d_v) \cong H(M,
d_v)\). 
Each factor \(M_c\) is supported in a single value of \(\gr_h\), so
 \(M\) and \(H(M,d_v)\) are supported in a single value of \(\gr_h\) as
 well. Thus the spectral sequence has converged at the
 \(E_1\) term, and 
\begin{equation*}
H(C_r(D),d_v+d_+) \cong H(M,d_v).
\end{equation*}

To evaluate  \(H(M,d_v)\), we observe that  \(M_c\) has
 the same form as the complex \(C_p(D_r)\) introduced in
the proof of Lemma~\ref{Lem:States}, but with \(d_v\) in place of
 \(d_h\). Thus \(H(M,d_v) \cong \H^+(D')\), where \(D'\) is the abstract graph
 obtained by replacing each crossing of \(D\) with a ``crossing'' of type
 \(D_r\). Since the
 differential in \(M_c\) is multiplication by \(X_k-X_j\), the ends
 labeled \(j\) and \(k\) lie on a solid segment of \(D_r\),
 as do the ends labeled \(i\) and \(l\). 
By
Lemma~\ref{Lem:Marks}, \(\H^+(D') \cong \H^+(D'')\), where \(D''\) is
 the abstract graph obtained by erasing the dashed lines in each copy
 of \(D_r\).  In other words, \(D''\) is  obtained
 by thinking of  \(L\) as a topological space and entirely 
forgetting its embedding in \(R\R^3\).
Thus \(D''\) is a disjoint union of \(\ell\) circles and  
 \(H(M,d_v) \cong \H(U^\ell)\). 
\end{proof}

When \(\ell=1\),  \(H(M)\) is one-dimensional, and is supported in the
top homological grading of \(M\).  If we compute the gradings \(\gr_+\)
and \(\gr_{-1}'\) for this generator, we find that they are given by
\(\gr_{+} = -2w\) and \(\gr_{-1}' = 2w\), where \(w\) is the writhe of
\(D\). Together with the overall shift of \(\{-w+b-1,w+b-1,w-b+1\}\)
in the triple grading, this means that the total homology 
\(H(C_r(D),d_v+d_+)\) is supported in gradings \(\gr_+ = \gr_{-1}' =
0\). More generally, the total homology will have Poincar{\'e}
polynomial
\begin{equation*}
\sum_{\substack{i = \gr_+  \\
j = \gr_{-1}'}} t^iq^j \dim  H^{i,j}(C_r(D),d_v+d_+) = \left(\frac{1+t^{-1}}{1-q}
\right)^{\ell-1} 
\end{equation*}

We can also prove an analog of Corollary~\ref{Cor:AntiCom}.

\begin{lem}
The differential \(d_k(-1)\) anticommutes with \(d_1(p)\) for any
value of \(p\). 
\end{lem}

\begin{proof} In this situation, it's more  convenient to work with
 \(\Cr_p(D,i)\) than \(C_r(D)\).  The isomorphism between the two described in
  Lemma~\ref{Lem:GlobalReduced} clearly respects their structure as
  double complexes, so we can think of \(E_k(-1)\) as being the
  spectral sequence induced by the horizontal filtration on
  \((\Cr_p(D,i),d_+,d_v)\). \(d_-\) anticommutes with both \(d_+\) and
  \(d_v\), so it defines
 a morphism \(d_-:  \Cr_p(D,i) \to
  \Cr_p(D,i)\) in the category \(Kom(Kom(R(D)))\). There is an induced 
  morphism  \((d_-)_k: E_k(-1) \to
  E_k(-1)\) which anticommutes with \(d_k(-1)\). The map 
\((d_-)_k\) is induced by the action of \(d_-\) on \(E_k(-1)\). 
  In particular, 
\((d_-)_2\) is the map induced by \(d_-\) on \(E_2(-1) \cong
  \H(L)\) --- in other words, \((d_-)_2 = d_1(p)\). This proves the claim.
\end{proof}

\section{Examples}
\label{Sec:Misc}
In this section, we compute the
KR-homology of some simple knots.  We begin by giving a quick proof of
Theorem~\ref{Thm:Limit}. Next, we discuss the notion of a KR-thin knot and
show that two-bridge knots are KR-thin. We then derive a skein exact
sequence  which is useful for making calculations. Combining this with
some computations of Webster \cite{Webster}, we are able to  determine
the KR-homology of all knots with \(9\) crossings or fewer.

\subsection{Homology of knots}

The reduced homology of a knot has the following important property:

\begin{prop}
\label{Prop:FGen}
If \(K\) is a knot, then \(\H(K)\) is a finite dimensional vector
space over \(\Q\).  
\end{prop}

\begin{proof}
Suppose \(D\) is a diagram representing a link \(L\). 
The complex \(\Cr_p(D,i)\) is a finitely generated module over the
ring \(R(D)\), so \(\H(D,i)\) is finitely generated over the ring
\(\Q[X_j]\), where \(j\) runs over the edges of \(D\). According to
Lemma~\ref{Lem:EdgeHom}, multiplication by \(X_j\) and \(X_k\) are
homotopic as morphisms of \(\Cr_p(D,i)\) whenever \(j\) and
\(k\) belong to the same component of \(L\). In particular, if \(L=K\)
is a knot, multiplication by any \(X_j\) is homotopic to multiplication by
\(X_i\), which is the zero map on \(\Cr_p(D,i)\). It follows that all
\(X_j\) act trivially on \(\H(K,i)\cong \H(K)\), so this group is 
finitely generated over \(\Q\). 
\end{proof}

\vskip0.05in
\noindent
{\it Proof of Theorem~\ref{Thm:Limit}:} Since \(\H(K)\) is finite
  dimensional, it is supported in finitely many \(q\)--gradings.
Consider the spectral sequence \(E_k(N)\) which relates \(\H(K)\) to
\(\H_N(K)\). The \(k\)-th differential in this sequence raises the
\(q\)--grading by \(2kN\). Thus when \(N\) is sufficiently large, all
the differentials beyond \(d_0\)
must vanish, and the sequence has converged at the
\(E_1\) term. \qed

\subsection{KR-thin knots}  

In both Khovanov homology and knot Floer homology, the simplest knots
exhibit a very similar pattern of behavior, in which there is a linear
relation between the two gradings and the signature of the knot.
  Such knots are said to be  {\it thin}. 
An analogous definition of
thinness in the context of KR-homology was proposed in \cite{superpolynomial}.
 In terms of our current normalizations, it is

\begin{figure}
\includegraphics{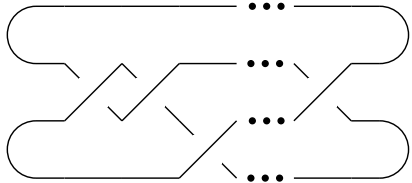}
\caption{\label{Fig:Twobridge} Plat diagram of a two-bridge knot.}
\end{figure}

\begin{defn}
A knot \(K\subset S^3\) is {\it KR-thin} if \(\H^{i,j,k}(K) = 0 \)
whenever \(i+j+k \neq \sigma (K)\). 
\end{defn}

Our sign convention for \(\sigma\)  is that positive knots have positive
signature. The quantity \(\delta = i + j +k \) which appears in the
definition occurs  frequently, and it is often convenient to think of
the grading on \(\H\) as being determined by the triple
\((i,j,\delta)\), rather than \(i,j,k\). 
From this point of view, it is clear that 
the HOMFLY homology of a KR-thin knot is completely determined by its
signature and HOMFLY polynomial. The same statement holds for the \(sl(N)\) 
homology as well:

\begin{prop}
\label{Prop:Thin}
If \(K\) is KR-thin, then the isomorphism of Theorem~\ref{Thm:Limit}
 holds  for all \(N>1\).
\end{prop}

\begin{proof} Consider the spectral sequence \(E_k(N)\) which relates
  \(\H(K)\) to \(\H_N(K)\). 
  The differential \(d_k(N)\) 
   is triply-graded of degree \((2kN,-2k,2-2k)\), so it raises
  \(\delta=i+j+k\) by \(2+2k(N-2)\). This quantity is positive whenever
  \(N>1\). Since \(E_1(N) \cong \H(K)\) is supported in a
  single value of \(\delta\), it follows that \(d_k \equiv 0\) for all
  \(k>0\), and the sequence converges at the \(E_1\) term. 
\end{proof}

In \cite{khthin}, knots for
which \(\H_N\) took this form were called \(N\)--thin. In this
language, the proposition says that if a knot is KR-thin, then it is
\(N\)--thin for all \(N>1\). 
Conversely, we  have the following result, which is an immediate
consequence of Theorem~\ref{Thm:Limit}.
\begin{prop}
If \(K\) is
\(N\)-thin for all sufficiently large \(N\), then \(K\) is
KR-thin.
\end{prop}

The  main result of \cite{khthin}  says that
two-bridge knots are \(N\)-thin for all \(N>4\), so they are KR-thin as well. 
We thus arrive at the statement of 
Corollary~\ref{Cor:TwoBridge} from the introduction. 

Curiously, it seems difficult to prove this result  without
appealing to the \(sl(N)\) homology. The issue is that \(\H_N(K)\) can be
computed using any planar diagram of \(K\), whereas the definition of
\(\H(K)\) requires that we use a braid diagram. Any two-bridge knot
admits a simple plat diagram of the form shown in
Figure~\ref{Fig:Twobridge}, which can be used to compute
\(\H_N(K)\). In contrast, the minimal braid diagram of such a knot can
be quite complicated, and there does not seem to be an easy way
 to compute \(\H(K)\) from it.

By a well-known theorem of Lee \cite{Lee}, the Khovanov homology of
any alternating knot is thin. In \cite{OS8} Ozsv{\'a}th and Szab{\'o}
proved a similar
result for the knot Floer homology. As observed in
\cite{superpolynomial}, however, it cannot be the case that all
alternating knots are KR-thin. Indeed, the HOMFLY polynomial of a
KR-thin knot must be {\it alternating}, in the sense that 
\begin{equation*}
P_K(a,q) = (-1)^{\sigma(K)} \sum_{i,j} c_{ij} a^{2j}(-q^2)^{i}
\end{equation*}
with \(c_{ij} \geq 0 \). On the other hand, it is not difficult to
find alternating knots whose HOMFLY polynomials are not alternating. 

We conclude our discussion of KR-thin knots by considering their
behavior with respect to the spectral sequences \(E_k(1)\) and
\(E_k(-1)\). We have already seen that the sequences \(E_k(N)\)  are
essentially trivial for \(N>1\). This cannot be true for \(E_k(\pm1)\),
however, since they converge to \(\Q\). Instead, we have

\begin{lem}
If \(K\) is KR-thin, then the spectral sequences \(E_k(1)\)
and \(E_k(-1)\) converge after the first differential on
\(\H(K)\). (That is, at \(E_2(1)\) and \(E_3(-1)\).)
\end{lem}

\begin{proof}
As we saw in the proof of Proposition~\ref{Prop:Thin}, the
differential \(d_k(1)\) shifts \(\delta\) by \(2-2k\). Since
\(E_1(1)\cong \H(K)\) 
is supported in a single \(\delta\)-grading, \(d_k(1)\) is
trivial for all \(k>1\). Similarly, the differential \(d_k(-1)\) is
triply graded of degree \((2-2k,2-2k,2k)\), so it shifts \(\delta\) by 
\(4-2k\). Thus it vanishes for all \(k>2\). 
\end{proof}
\noindent
It follows that the spectral sequence of a KR-thin knot behaves as
conjectured in \cite{superpolynomial}.

\subsection{A skein exact sequence}
Suppose \(D\) is a planar diagram  representing a two-component
link \(L\), and that \(i\) and \(j\) are edges of \(D\) belonging to the two
components of \(L\). Let \(\Cr_N(D,i) = (H(\Cr_p(D,i),d_{tot}),d_v^*)\) be
the \(sl(N)\) chain complex, and form the mapping cone
\begin{equation*}
\Cr_N(D,i,j) = \Cr_N(D,i)\{1,0,-1\}
\xrightarrow{\phantom{X}X_j\phantom{X}}\Cr_N(D,i)\{-1,0,1\}
\end{equation*}
 where \(X_j\) denotes the map induced by multiplication by \(X_j\).
The grading shifts are chosen so that \(X_j\) --- like \(d_v^*\) ---
is homogenous of degree \((0,0,2)\) with respect to the triple
grading. 
We call the homology of \(\Cr_N(D,i,j)\)  the {\it totally
  reduced homology} of \(L\) and  denote it by
\(\overline{\H}_N(L)\). Using Lemma~\ref{Lem:EdgeHom}, it is not
difficult to see that \(\overline{\H}_N(L)\) is independent of the
choice of \(i\) and \(j\), although we will not use this fact here. 

The group \(\H_N(L,i)\) can naturally be viewed as a \(\Q[X]\) module,
where \(X\) acts as multiplication by \(X_j\). 
If we understand the module structure of 
\(\H_N(L,i)\), we can easily determine the 
totally reduced homology from the long exact sequence 
\begin{equation*}
\begin{CD}
@>>> \overline{\H}_N(L) @>>> \H_N(L,i) @>{X_j}>> \H_N(L,i) @>>>
\overline{\H}_N(L) @>>> .
\end{CD}
\end{equation*}
One such case is when \(L\) is a two-bridge link. In \cite{khthin}, it
was shown that 
\(\H_N(L,i)\) is composed of a number of summands on which \(X\) acts
trivially, together with a single summand isomorphic to
\(\Q[X]/X^{N-1}\). The generators of each summand have
\(\delta\)--grading congruent to \(\sigma(L) \mod (2N-4)\). (On
\(\H_N\), we can't tell the difference between \(a^2\) and \(q^{2N}\),
so the \(\delta\)--grading is  only defined modulo
\((2N-4)\).) From this, it is not difficult to see that
\(\overline{\H}_N(L)\) is also thin, in the sense that it is supported
in \(\delta\)--gradings congruent to \(\sigma(L)\). In analogy with
the case of knots, we say that \(L\) is KR-thin if
\(\overline{\H}_N(L)\) is thin for all \(N \gg 0 \). 

For the moment, our interest in the group \(\overline{\H}(L)\) arises from 
the following skein exact sequence, which generalizes the oriented
skein relation for the \(sl(N)\) polynomial. 

\begin{prop}
\label{Prop:Skein}
Suppose \(L_+\) and \(L_-\) are two knots related by a crossing
change, and \(L_0\) is the two-component link obtained by resolving
the crossing. Then there is a long exact sequence 
\begin{equation*}
\begin{CD}
@>{(0,0)}>> \H_N(L_-) @>{(N,1)}>> \overline{\H}_N(L_0) @>{(N,1)}>>
\H_N(L_+) @>{(-2N,0)}>> \H_N(L_-) @>{(N,1)}>>.
\end{CD}
\end{equation*}
\end{prop}

The numbers over each arrow indicate the degree of the corresponding
map with respect to the \((q,\delta)\) bigrading on \(\H_N\). For
example, the map \(\overline{\H}_N(L_0) \to \H_N(L_+)\) raises the
\(q\)--grading by \(N\) and \(\delta\) by \(1\). 

\begin{proof}
The complex \(\Cr_N(L_-)\) is the mapping cone of the map \(\chi_0:
\Cr_N(L_0) \to \Cr_N(L_s)\), where \(L_s\) is the diagram obtained by
replacing the crossing in question with a singular point. Similarly,
\(\Cr_N(L_+)\) is the mapping cone of \(\chi_1: \Cr_N(L_s) \to
\Cr_N(L_0)\), from which it follows that \(\Cr_N(L_s)\) is homotopy equivalent to
the mapping cone of the inclusion \(i: \Cr_N(L_0) \to \Cr_N(L_+)\). An explicit homotopy equivalence is given by the
map 
\begin{equation*}
\iota: \Cr_N(L_s) \to  \Cr_N(L_0) \oplus \Cr_N(L_s) \oplus \Cr_N(L_0)
\end{equation*}
 which sends \(a \in  \Cr_N(L_s)\) to \((-\chi_1(a),a,0)\). It is easy
 to see that \(\iota\) is the inclusion in a strong
 deformation retract, in the sense of Bar-Natan \cite{DBN2}. By Lemma
 4.5 of \cite{DBN2}, it follows that \(\Cr_N(L_-)\) is homotopy equivalent
 to the mapping cone of \(\iota \chi_0: \Cr_N(L_0) \to Cone(i)\). 
This complex has a three-step filtration, as illustrated below:
\vskip0.07in
$$
\xymatrix{\Cr_N(L_0) \ar[r]^{-\chi_1\chi_0} \ar@/^1.8pc/[rr]  & \Cr_N(L_0) \ar[r] & \Cr_N(L_+).
}
$$
It follows that there is a short exact sequence
\begin{equation*}
\begin{CD}
0 @>>> \Cr_N(L_+) @>>> Cone(\iota \chi_0) @>>> Cone(\chi_1 \chi_0)
@>>> 0 .
\end{CD}
\end{equation*}
Considering the associated long exact sequence, we see that to prove
the lemma, it suffices to show that  \(Cone(\chi_1 \chi_0) \cong
\Cr_N(D,i,j)\). To show this, 
recall that 
the composition \(\chi_1 \chi_0\) is multiplication by \(X_j - X_i\),
where \(i\) and \(j\) are the edges of \(L_0\) adjacent to the
resolution. Taking these two edges to be the edges \(i\) and \(j\)
which appear in the definition of \(\overline{\H}(L)\) gives
 the desired isomorphism. Finally, the bigrading of each map in the
 sequence can easily be determined from Lemma 3.3 of \cite{khthin}. This
 is left as an exercise to the reader. 
\end{proof}

As an application, we have the following criterion for showing that a
knot is KR-thin. It is a slight generalization of criterion 5.4 from
\cite{khthin}. 

\begin{cor}
\label{Cor:ThinCrit}
Suppose that \(L_-\), \(L_0\), and \(L_+\) are as above,
that \(L_-\) and \(L_0\) are both KR-thin, and that 
\(\det L_- + 2 \det L_0 = \det L_+\). Then \(L_+\) is KR-thin as well. 
\end{cor}

\begin{proof}
Suppose  \(N\) is very large. Then all three terms in the exact sequence
stabilize --- \(L_0\) by hypothesis, and \(L_-\) and \(L_+\) by
Theorem~\ref{Thm:Limit}. We have
\begin{equation*}
\text{rank} \ts \H_N(L_+)  \geq \det L_+ 
 = \text{rank} \ts \H_N(L_-) + \text{rank} \ts \overline{\H}_N(L_0)
\end{equation*}
since both \(L_-\) and \(L_0\) are thin. For this to happen, the map \(
\H_N(L_-) \to  \overline{\H}_N(L_0)\) in the exact sequence must
vanish. To show that \(L_+\) is \(N\)-thin,  it is enough to check that
\(\sigma(L_+) = \sigma (L_-) = \sigma (L_0)+1.\) This follows 
from the usual skein-theoretic constraint on the signature.
 (See the proof of criterion
5.4 in \cite{khthin} for details.) Finally, 
since \(L_+\) is \(N\)-thin for all
large \(N\), it is KR-thin as well. 
\end{proof}

\noindent The analogous statement with the roles  of
\(L_-\) and \(L_+\) reversed also holds. 



\subsection{Connected sums}
In applying the skein exact sequence of the previous section, one
often encounters non-prime knots and links. For this reason, it is
convenient to understand the behavior of the KR-homology under
connected sum. 

Suppose  \(L_1\) and
\(L_2\) are oriented links with marked components \(i_1\) and \(i_2\).
Up to isotopy, there is a unique was to form their 
 orientation-preserving connected sum along \(i_1\)
and \(i_2\). We denote the resulting link by \(L_1 \#_{i_1=i_2}L_2\). 

\begin{lem}
There are isomorphisms
\begin{align*}
\H(L_1 \#_{i_1=i_2} L_2) & \cong \H(L_1) \otimes \H(L_2) \\
\H_N(L_1 \#_{i_1=i_2} L_2,i_1) & \cong \H_N(L_1,i_1) \otimes
\H_N(L_2,i_2).
\end{align*}
\end{lem}

\begin{figure}
\includegraphics{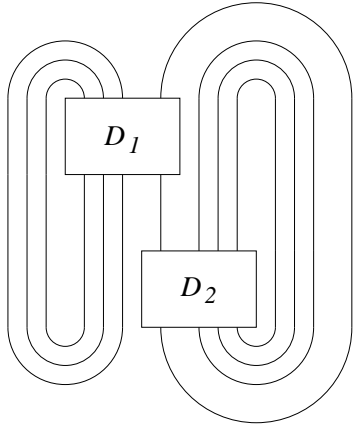}
\caption{\label{Fig:ConnectSum} The connected sum of braids \(D_1\)
  and \(D_2\).}
\end{figure}

\begin{proof}
Suppose \(L_i\) is represented by a braid diagram \(D_i\) on
\(b_i\) strands. Without loss of generality, we may arrange  the
diagrams \(D_i\) so \(i_1\) is the rightmost strand 
of \(D_1\) and  \(i_2\) is the  leftmost strand of \(D_2\). Then the
connected sum \(L_1\#_{i_1=i_2}L_2\) can be represented by a braid
diagram \(D_1\#D_2\) on \(b_1+b_2-1\) strands, as illustrated in
Figure~\ref{Fig:ConnectSum}. Let \(D_1^o\) be the open diagram  obtained by removing a neighborhood of the connected sum point in \(D_1\), and let  \(i_1^+\) and \(i_1^-\) be its free ends. 
Then \(X_{i_1^+}=X_{i_1^-}\) in \(R(D_1^o)\), and
\begin{equation*}
C_p(D_1) \cong  C_p(D_1^o)\vert_ {X_{i_1^+}=X_{i_1^-}} \cong C_p(D_1^o). 
\end{equation*}
From the local nature of the KR-complex, we see that 
\begin{align*}
C_p(D_1\#D_2) & \cong C_p(D_1^o)\otimes C_p(D_2^o) \vert_{X_{i_1^+} = X_{i_2^-},
X_{i_1^-}=X_{i_2^+}} \\
& \cong C_p(D_1) \otimes  C_p(D_2) \vert_{X_{i_1}=X_{i_2}}.
\end{align*}
It follows that the reduced KR-complex has the form
\begin{equation*}
C_p(D_1\#D_2,i_1) \cong  C_p(D_1,i_1) \otimes  C_p(D_2,i_2).
\end{equation*}
Applying the Kunneth formula twice gives the statement of the lemma. 
\end{proof}

\begin{cor}
The connected sum of two KR-thin knots is KR-thin. 
\end{cor}

\noindent Similarly, it is not difficult to see that the connected sum
of a 
KR-thin knot and a KR-thin link is also KR-thin.

\subsection{Small knots}
\begin{figure}
\includegraphics{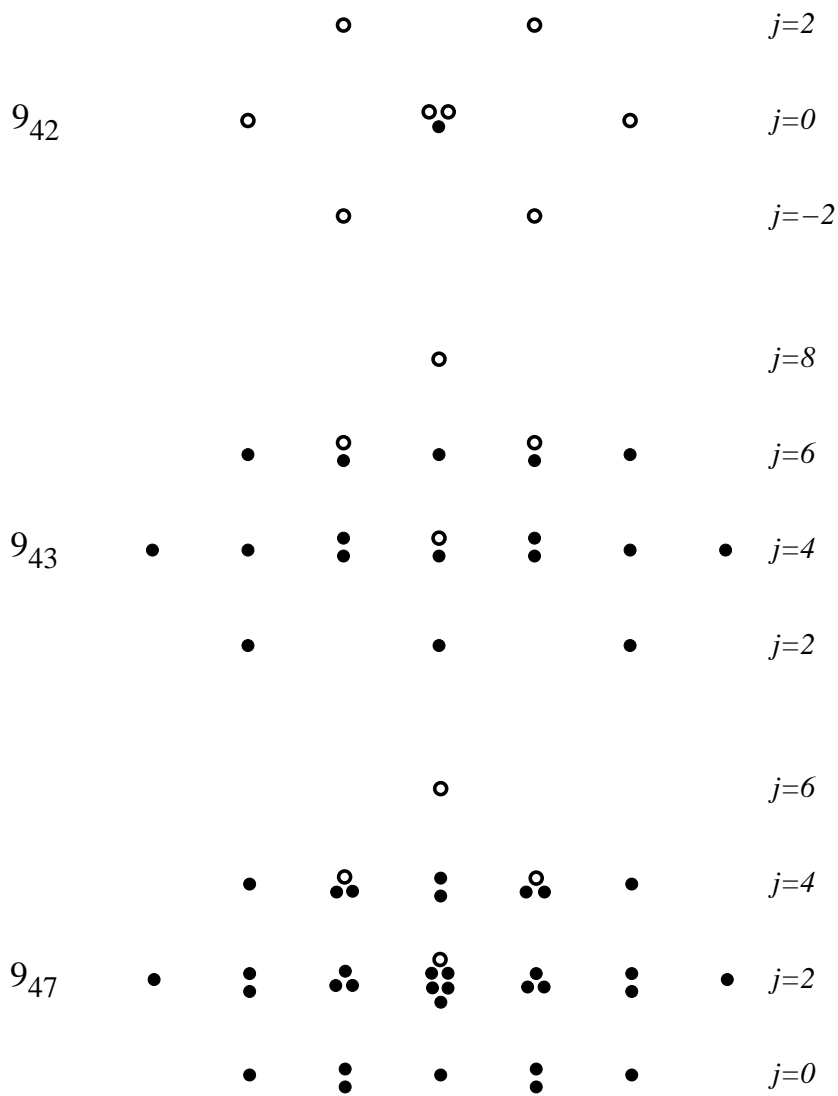}
\caption{\label{Fig:NonThin} HOMFLY homology of the knots 
\(9_{42}\), \(9_{43}\) and \(9_{47}\), represented by ``dot diagrams.''
 Each dot represents a
generator of the homology. The \(i\) and \(j\) gradings are indicated
by the positions in the horizontal and vertical directions,
respectively. (\(i=0\) corresponds to the axis of symmetry.) The solid
and hollow dots have different \(\delta\)--gradings. For \(9_{42}\),
hollow dots have \(\delta = -2\), and the solid dot has \(\delta =
0\). For \(9_{43}\), the values are \(\delta = 2\) and \(\delta = 4\),
respectively, and for \(9_{47}\), they are \(\delta = 0\) and \(\delta
= 2\).}
\end{figure}
We conclude by describing the KR-homology of knots with 9 or fewer
crossings. Previous computations of \(\H\) have been made by 
 Khovanov
\cite{KhSoergl}, who showed that the \((2,n)\) torus knots
 are KR-thin, and by Webster \cite{Webster},
who wrote a computer program for this purpose. Using it, he computed
 \(\H\) for knots up through \(7\) crossings, all
of which are KR-thin. For larger knots, the program is very effective
at computing the homology of knots which can be represented as
 closures of three-strand braids, but less useful in other
cases. Fortunately, many  of the small knots with large braid index
 are two-bridge, and thus covered by Corollary~\ref{Cor:TwoBridge}.  
The remainder can be analyzed  using the skein exact sequence
of Proposition~\ref{Prop:Skein}. Combining the information from these various
sources, we have

\begin{prop}
The only knots with \(9\) crossings or fewer which are not KR-thin are
\(8_{19}\), \(9_{42}\), \(9_{43}\), and \(9_{47}\). (Numbering as in
 Rolfsen   \cite{Rolfsen}.)
\end{prop}

\noindent{\bf Remarks:} The homology of 
\(8_{19}\) (the \((3,4)\) torus knot) was computed by Webster
\cite{Webster}. The homology of  the remaining three knots 
is illustrated in Figure~\ref{Fig:NonThin}. In all four cases, the homologies
are symmetric; the sequences \(E_k(-1),E_k(1)\), and \(E_k(2)\)
converge after the first differential on \(\H(K)\); and \(\H_N(K)
\cong \H(K)\) for \(N>2\). In addition, the calculated values of
\(\H(8_{19})\) and \(\H(9_{42})\) agree with the predictions made in
\cite{superpolynomial}.

\begin{figure}
\includegraphics{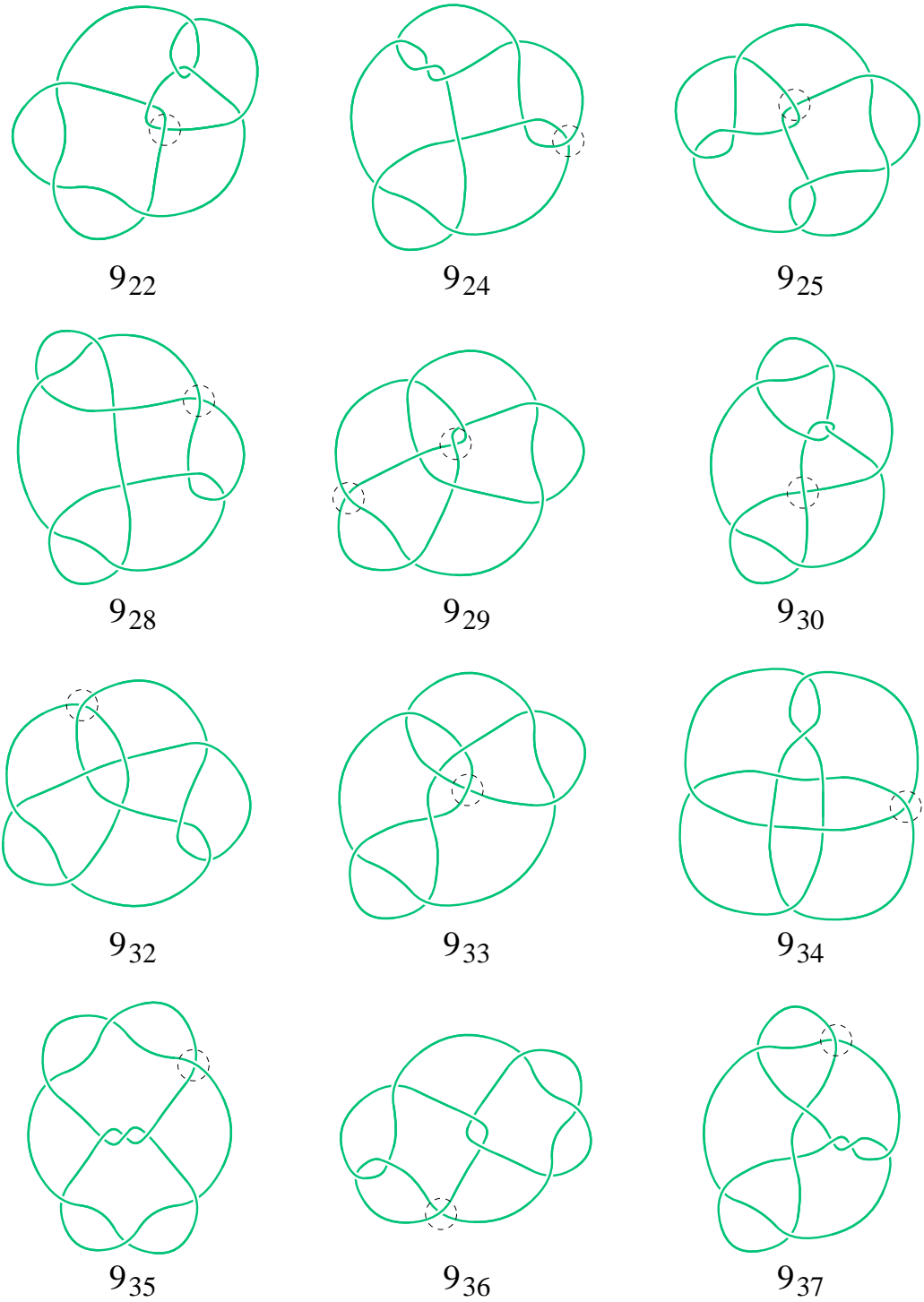}
\caption{\label{Fig:NineCrossing} 9 crossing knots (I). Figures drawn
by {\it Knotscape} \cite{knotscape}.}
\end{figure}

\begin{figure}
\includegraphics{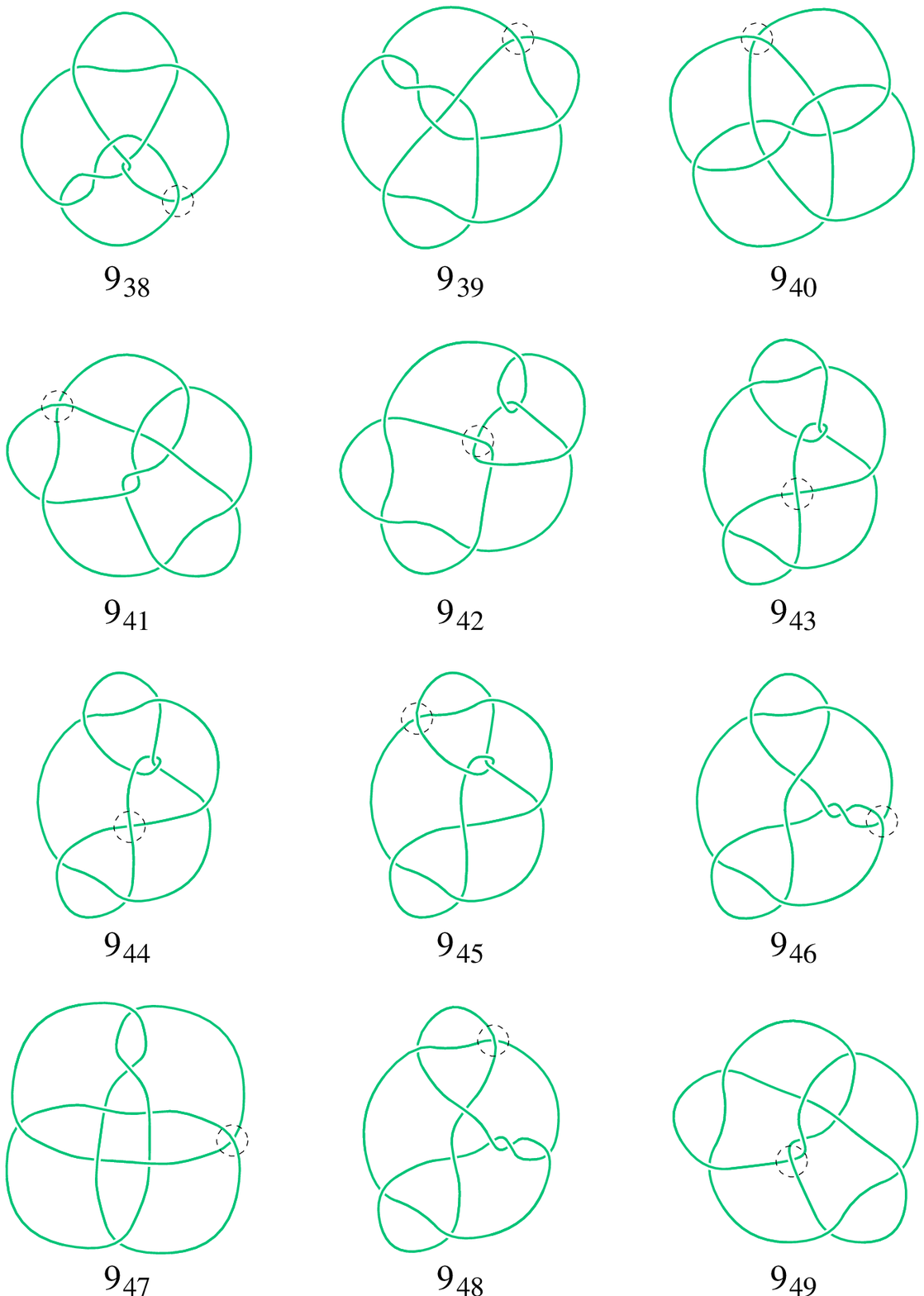}
\caption{\label{Fig:NineCrossing2} 9 crossing knots (II). Figures
  drawn by {\it Knotscape} \cite{knotscape}.}
\end{figure}

\begin{proof}

The only knots with \(8\) or fewer crossings which are not two-bridge
are \(8_{5},8_{10}\), and \(8_{15}\)--\(8_{21}\). Among these, all but
\(8_{15}\) have braid index \(3\) and were computed by Webster
\cite{Webster}. The only one which is not thin is the \((3,4)\)
torus knot \(8_{19}\). In \cite{khthin}, it was shown that \(8_{15}\)
is \(N\)-thin for all \(N>4\). Thus it is also KR-thin.

For the \(9\)-crossing knots, we have to work a little harder. The knots
\(9_{16}\), \(9_{22}\), \(9_{24}\), \(9_{25}\), \(9_{28}\),
\(9_{29}\), \(9_{30}\), and \(9_{32}\)--\(9_{49}\) are not
two-bridge. Only one -- \(9_{16}\)-- is the closure of a three-strand
braid, and Webster's program shows that it is KR-thin. Of the rest,
all but five --- \(9_{29}\),\(9_{42}\),\(9_{43}\), \(9_{46}\), and
\(9_{47}\) ---  can be shown to be KR-thin using the criterion of
Corollary~\ref{Cor:ThinCrit}. 
 These knots, and the crossing to which the criterion can be applied,
 are shown in Figures~\ref{Fig:NineCrossing} and
 \ref{Fig:NineCrossing2} at the end of the paper.

The knot \(9_{29}\) can be seen to be KR-thin by a similar, but
slightly more elaborate argument. If we change the marked crossing on
the right side of the knot in Figure~\ref{Fig:NineCrossing}, we get
the two-bridge knot \(7_6\). Resolving the crossing gives the link
\(7_4^2\). We claim that this link is KR-thin. To see this, we
consider the second marked crossing in the figure. Changing this
crossing gives the connected sum of the Hopf link and the trefoil
knot, which is KR-thin and has determinant \(6\). Resolving the
crossing gives the knot \(5_1\), which has determinant \(5\). Since
the determinant of \(7_4^2\) is \(16=6+2\cdot5\), it's not difficult to
see that \(7_4^2\) is KR-thin. Then \(9_{29}\) must be
KR-thin as well. 

\begin{figure}
\includegraphics{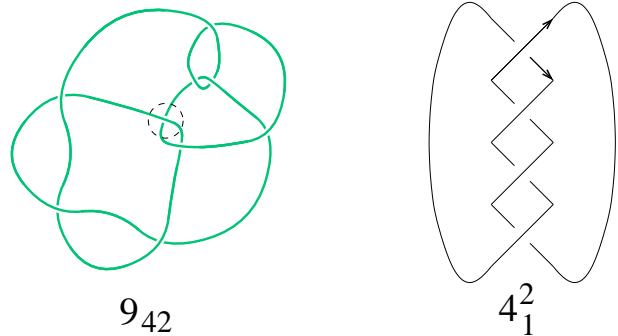}
\caption{\label{Fig:942} We consider the skein exact sequence
  associated to the circled crossing in diagram of \(9_{42}\) shown
  above. \(L_0\) is the link \(4^2_1\), oriented as shown.}
\end{figure}

To analyze the four remaining knots, we  resort to a more
detailed study of the skein exact sequence.
We illustrate this process in the case of  
 the knot \(9_{42}\), which is shown in Figure~\ref{Fig:942}. 
If we change the circled crossing in the figure 
 from positive to negative, the result is the connected
sum of the negative trefoil and the figure-eight knot. Resolving the
crossing, we get the two-bridge link \(4^2_1\) shown on the right-hand side of
the figure.  Thus we have a long exact sequence 
\begin{equation*}
\begin{CD}
@>>> \H_N(3_1\#4_1) @>>> \overline{\H}_N(4_1^2) @>>> \H_N(9_{42}) @>>>
  \H_N(3_1\#4_1) @>>>.
\end{CD}
\end{equation*}
When \(N\) is large, all three terms in this sequence stabilize. Both
\(3_1\#4_1\) and \(4_1^2\) are KR-thin, so their homologies are
determined by their HOMFLY polynomials. 
In Figure~\ref{Fig:942HOM},  we have superimposed diagrams
representing the homology of \(3_1\#4_1\) (hollow dots) 
and \(4_1^2\)
(solid dots). 
The \(j\)-gradings are shifted so that they correspond to the power of \(a\) in
the HOMFLY polynomial of \(9_{42}\). Under the assumption that \(N\)
is large, nontrivial  components of the map 
 \(\H_N(3_1\#4_1) \to \overline{\H}_N(4_1^2)\) must preserve the position
 of the generators. In other words, a generator corresponding to a
 hollow dot at any of the lettered positions
 can map nontrivially to a solid dot at
 the same position, but not to anything else.

From the figure, we can deduce some constraints on  the group 
\(\H_N(9_{42})\). For example, the group at position \(c\) must
have rank either \(2\) or \(0\), depending on whether the map from the
hollow to the solid generator at that position is trivial or
nontrivial. Since \(\H_N(9_{42})\cong \H(9_{42})\) when \(N\) is
large, the same is true for \(\H\) as well. 

We can now use Theorems~\ref{Thm:d+} and \ref{Thm:d-1} to deduce the
 exact value of the homology. For example, suppose
 the two generators labeled \(a\) on the right-hand side
of the figure  survive in
\(\H(9_{42})\). They must die in the spectral sequence
\(E_k(1)\), but there is nothing to kill them.
 We conclude that these
generators could not have survived. A similar argument using
\(E_k(-1)\) shows that the two left-hand generators labeled \(a\) do
not survive either. It is now easy to see that generators labeled \(b\) must
kill each other too. 
\begin{figure}
\includegraphics{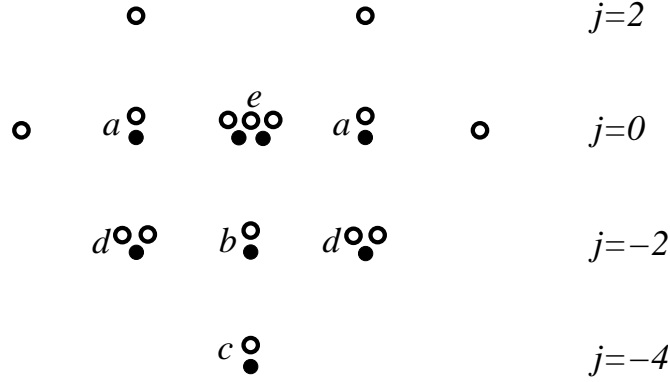}
\caption{\label{Fig:942HOM} Possible generators of \(\H(9_{42})\).}
\end{figure}
\end{proof}
To eliminate the generators labeled \(c\), we
consider the spectral sequence \(E_k(2)\), which converges to the
usual Khovanov homology. Clearly, if these generators survive in
\(\H(K)\), they will also appear in \(\H_2(K)\), where they will have
\(q\)--grading \(-8\). On the other hand, it is well known that
\(\H_2(9_{42})\) has Poincar{\'e} polynomial
\begin{equation*}
\mathcal{P}_2(9_{42}) = q^{-6}t^{-4} + q^{-4}t^{-3} + q^2t^{-2} +
2t^{-1} + 1 + q^2 + q^4 t + q^6 t^2. 
\end{equation*}
 There is no term with \(q^{-8}\), so the generators in position \(c\)
 must die. Next, we consider the positions
labeled \(d\), where we have a map from a two-dimensional space generated by
the hollow dots to a one-dimensional space generated by the solid
dot. Now that we know that there is nothing in position \(c\),
considering \(E_k(\pm1)\) shows that both maps must be
surjective. Finally, in position \(e\), we have a map from a space of
dimension 3 to a space of dimension 2. Examining the sequence
  \(E_k(2)\) shows that this map must have rank
\(1\). Thus the homology is as shown in Figure~\ref{Fig:NonThin}.

Similar considerations may be applied to compute the homology of the
knots \(9_{43},9_{46}\), and \(9_{47}\). 
 Rather than go into details, we simply indicate an
appropriate crossing for each knot in Figure~\ref{Fig:NineCrossing2},
and leave it to the interested reader to check the rest.

\bibliographystyle{plain}
\bibliography{../mybib}

\end{document}